\documentclass{article}
\usepackage{amsmath,amssymb,amsthm,amsfonts,mathabx}
\usepackage[all]{xy}

\usepackage[margin=1in]{geometry}
\usepackage{titling}
\usepackage{fancyhdr,indentfirst,color,hyperref}

\newtheorem{lemma}{Lemma}

\newtheorem{theorem}{Theorem}
\newtheorem{proposition}{Proposition}
\newtheorem{definition}{Definition}

\newtheorem{corollary}{Corollary}

\newtheorem{manualtheoreminner}{Theorem}
\newenvironment{manualtheorem}[1]{%
  \renewcommand\themanualtheoreminner{#1}%
  \manualtheoreminner
}{\endmanualtheoreminner}
\newtheorem{manualcorollaryinner}{Corollary}
\newenvironment{manualcorollary}[1]{%
  \renewcommand\themanualcorollaryinner{#1}%
  \manualcorollaryinner
}{\endmanualcorollaryinner}

\def\Hom{\operatorname{Hom}}
\def\Eu{\operatorname{Eu}}
\def\Res{\operatorname{Res}}
\def\det{\operatorname{det}}
\def\ev{\operatorname{ev}}
\def\ft{\operatorname{ft}}
\def\tr{\operatorname{tr}}
\def\ch{\operatorname{ch}}
\def\id{\operatorname{id}}

\def\rank{\operatorname{rank}}
\def\mod{\operatorname{mod}}

\def\virt{\operatorname{virt}}

\def\Rep{\operatorname{Rep}}
\def\Fix{\operatorname{Fix}}
\def\C{\mathbb{C}}
\def\J{\mathcal{J}}

\def\K{\mathcal{K}}

\def\L{\mathcal{L}}

\def\M{\widebar{\mathcal{M}}}
\def\t{\mathbf{t}}
\def\f{\mathbf{f}}
\def\g{\mathbf{g}}
\def\e{\mathbf{e}}

\title{Loop Space Formalism and K-Theoretic Quantum Serre Duality}
\date{}
\author{Xiaohan Yan}

\begin{document}

\maketitle

\begin{abstract}
In this paper, we prove the quantum Serre duality for genus-zero K-theoretic permutation-invariant Gromov-Witten theory. The formulation of the theorem relies on an extension to the formalism of loop spaces and big $\J$-functions more intrinsic to quantum K-theory. With the extended formalism, we also arrive at a re-interpretation of the level structures in terms of twisted quantum K-theories. We discuss the torus-equivariant theory in the end, and as an application generalize the K-theoretic quantum Serre duality to non-primitive vector bundles over flag varieties.
\end{abstract}

\noindent {\bf\small Keywords:} {\it\small Gromov-Witten invariants, K-theory, quantum Serre duality, loop space formalism, level structures}

\section{Introduction}

\subsection{Some background}
\par In Gromov-Witten (GW) theory, one practical idea is to express the counting problem on more complicated varieties in terms of such problem on better-understood ones, possibly with some modification on the deformation-obstruction theory of the moduli spaces. In Coates-Givental \cite{Coates-Givental}, the authors considered the so-called twisted GW theory of a smooth projective variety $X$, under the very general settings where the virtual fundamental cycles on the moduli spaces of stable maps to $X$ are capped with a cohomological class determined by $(c,E)$, with $c$ an invertible characteristic class and $E$ a vector bundle over $X$, and expressed the twisted GW invariants (of all genera) in terms of untwisted ones. By taking $c$ as the Euler class and $E$ as a convex bundle, they were able to prove the quantum Lefschetz theorem which recovers the GW theory of complete intersections (cut out by generic sections of $E$) on $X$, providing thus a new proof to the mirror formula \cite{CDGP} for curve-counting on quintic threefolds and extending several earlier results. Moreover, with a fiberwise $C^\times$-action, the authors also found out that the GW invariants twisted by $(\Eu^{-1}, E^\vee)$ are closely related to those of $(\Eu, E)$. In this way, they identified the GW theory of a subvariety $Y$ in $X$ defined by a global section of $E$ with the GW theory of the total space of $E^\vee$, which generalizes an observation made in \cite{Givental:EGWI,Givental:elliptic} on projective spaces. In the genus-zero case, such relation specializes to an identification of the two twisted Givental cones up to scaling.

\par Such correspondence, which is named quantum Serre duality (qSD), or non-linear Serre duality, exists in various generalized settings as well. The closest formulation came in Tseng \cite{Tseng} where the author studied extensively the twisted GW theory for orbifolds (smooth DM stacks) and proved qSD for such targets. In Iritani-Mann-Mignon \cite{IMM}, the authors formulated qSD as a duality between the quantum D-modules. Under such interpretation, the authors proved that the duality admits a non-equivariant limit, and that it is compatible with the integral structures in quantum cohomology introduced in Iritani \cite{Iritani}. Such construction was generalized to non-compact orbifolds by Shoemaker in \cite{Shoemaker:narrow} through consideration of the narrow cohomology. In quasimap (defined for GIT quotients in \cite{CKM}) settings, Heath-Shoemaker \cite{HS} proved that qSD-type correspondence holds true directly on the level of virtual fundamental cycles and thus for individual two-point quasimap invariants. 

\par The quantum Serre duality has led to various applications in the field. In Lee-Priddis-Shoemaker \cite{LPS}, the authors combined qSD with the so-called MLK correspondence, which relates certain LG models with orbifold GW theory, and thus provided a new proof of the LG/CY correspondence (introduced in \cite{CR}) in terms of crepant transformation conjecture. A $D$-module formulation of the LG/CY correspondence \cite{CIR} is recovered by Shoemaker \cite{Shoemaker:correspondence} using the narrow cohomology construction mentioned above. Moreover, in their proof of the crepant transformation conjecture for toric complete intersections, Coates-Iritani-Jiang \cite{CIJ} used qSD to convert the quantum cohomology of such targets into the quantum cohomology of the total spaces of their dual bundles, which are toric objects and thus have simpler structure. More recently, Mi-Shoemaker \cite{MS} studied the the GW theory of varieties related by extremal transitions, further demonstrating the utility of qSD in birational geometry. 

\par Due to the relation known as the adelic characterization \cite{Givental-Tonita,Givental:perm9} between quantum K-theory and quantum cohomology, a K-theoretic version of the above story should exist as well. Quantum K-theory was introduced by Givental and Lee \cite{Givental:WDVV,Lee,Givental-Lee} to study a K-theoretic variation of the GW invariants. The counterpart in quantum K-theory of the virtual fundamental cycles are the virtual structure sheaves, defined as elements in the K-groups of the moduli spaces. More recently, by considering stable envelopes on Nakajima quiver varieties and of hyperplane arrangements, Aganagic, Okounkov, Smirnov and collaborators \cite{Okounkov,AO,OS,PSZ,KPSZ,AFO} have revealed more relations of quantum K-theory with representation theory (of Lie groups and quantum groups) and finite-difference equations (e.g. the qKZ equations). 

\par Twisted K-theoretic GW invariants of smooth projective varieties were studied by Tonita \cite{Tonita} and under permutation-invariant settings by Givental \cite{Givental:perm11}. A K-theoretic version of the quantum Lefschetz theorem was established as a special case of the so-called quantum Adams-Riemann-Roch (qARR) theorem. Moreover, in genus-zero case and under the assumption that $K(X)$ is generated by line bundles, one may recover from limited information the entire twisted Givental cone using the $D_q$-module structure studied in \cite{Givental:perm4} and even further the big quantum K-ring \cite{IMT} using a different system of difference operators. Such assumption is slightly relaxed in \cite{Givental-Yan} and \cite{Yan:flag} to include grassmannians and partial flag varieties, through the idea of abelian/non-abelian correspondence. In a different direction, however, a K-theoretic analogue of the qSD is still missing.

\subsection{Main theorems}
\par In this paper, we derive the K-theoretic quantum Serre duality (KqSD) for genus-zero permutation-invariant quantum K-theory, using the language of Givental cones. We hope that the results here may serve as a starting point to extend the story in quantum cohomology as above into the realm of quantum K-theory. 

\par Consider smooth projective variety $X$ and primitive vector bundle $E$ over $X$ (i.e. generated by line bundles in $K(X)$). We study the relation between the $(\Eu,E)$-twisted quantum K-theory and the $(\Eu^{-1},E^\vee)$-twisted theory, where $\Eu$ is the equivariant (thus invertible) K-theoretic Euler class. As before, the former may roughly be regarded as a counting problem on the subvariety $Y$ of $X$ cut out by a generic section of $E$, while the latter as on the total space of $E^\vee$. 

\par One key observation of this paper is that in K-theoretic settings, however, the most natural connection is no longer directly from the $(\Eu,E)$-twisted Givental cone to the $(\Eu^{-1},E^\vee)$-twisted one, but rather involves an extra level structure. Level structures as defined in Ruan-Zhang \cite{R-Z} as a feature of quantum K-theory, modify the virtual structure sheaves in a different way from the more classical twistings of the type $(c,E)$ for $c$ an invertible (K-theoretic) characteristic class. For integer $l$, we denote by $\L^{X,(\Eu,E,l)}$ the $(\Eu,E)$-twisted Givental cone further with level-$l$ structure of $E$, and similarly by $\L^{X,(\Eu^{-1},E^\vee,l+1)}$ the $(\Eu^{-1},E^\vee)$-twisted Givental cone further with level-$(l+1)$ structure of $E^\vee$. Then, our correspondence is in between these two cones. To ensure invertibility, we introduce an auxiliary fiberwise $\C^\times$-action on $E$, and denote its equivariant parameter by $\mu^{-1}$. We state the main theorem in the following two forms, in terms of either $\mu$ or $\mu^{-1}$.

\begin{manualtheorem}{A}[K-Theoretic Quantum Serre Duality] \label{MainA}
For any point
\[
J^1 = \sum_d \ \prod_{i=1}^n \left( Q_i' \right)^{d_i} \cdot J^1_d
\quad \in \L^{X,(\Eu,E,l)}
\]
with $J^1_d$ independent of $Q_i'\ (i=1,\cdots,n)$, we have
\[
J^2 = \Eu(\mu^{-1} E) \cdot \sum_d \ \prod_{i=1}^n \left(Q_i' \cdot \left(-q^l \right)^{c_1(E)_i}\right)^{d_i} \cdot J^1_d \quad \in \L^{X,(\Eu^{-1},E^\vee,l+1)}.
\]
\end{manualtheorem}
\begin{manualtheorem}{B}[K-Theoretic Quantum Serre Duality] \label{MainB}
For any point 
\[
J^1 = \sum_d \ \prod_{i=1}^n \left(Q_i''\right)^{d_i} \cdot J^1_d \quad \in  \L^{X,(\Eu,E,l)}
\]
with $J^1_d$ independent of $Q_i''\ (i=1,\cdots,n)$, we have
\[
J^2 = \Eu(\mu E^\vee) \cdot \sum_d \ \prod_{i=1}^n \left( Q_i'' \cdot \left(-q^{l+1}\right)^{c_1(E)_i}\right)^{d_i} \cdot J^1_d \quad \in \L^{X,(\Eu^{-1},E^\vee,l+1)} 
\]
\end{manualtheorem}
The more precise statements, as well as the definition (and motivation) of the modified Novikov variables $Q_i'$ and $Q_i''$, will be given in Section \ref{StatementQuantumSerre}. We state the two versions above separately because they hold true under different settings. Yet, they both follow from one ``root'' theorem (Theorem \ref{TheoremQuantumSerre}).

\par The two theorems indicate that, no matter expressed in terms of $\mu$ or $\mu^{-1}$, $\L^{X,(\Eu,E,l)}$ and $\L^{X,(\Eu^{-1},E^\vee,l+1)}$ differ by a global scaling followed by a change on the Novikov variables. In this way, we obtain a K-theoretic version of the Coates-Givental result \cite{Coates-Givental}. The content of this duality is different from the one with the same name in \cite{Givental:perm11}.

\subsection{Idea of proof and the rational loop space formalism}
\par The proof involves the other key observation of this paper, more technical this time, that the level structures may in fact be interpreted as certain limit of the composition of two invertible twistings, i.e. of certain types $(c,E)$, whose effects on the Givental cones are already understood by the qARR theorem given in \cite{Givental:perm11}. By definition, an invertible twisting of the type $(c,E)$ modifies the virtual structure sheaf $\mathcal{O}^{\virt}_{g,m,d}$ of the moduli space $\M_{g,m}(X,d)$ by replacing it by its tensor product with $c(V) = e^{\sum_k c_k \Psi^k(V)}$, where $V = \ft_*\ev^* E$ (see Section \ref{EulerTwisting}) is a (virtual) vector bundle over $\M_{g,m}(X,d)$. Meanwhile, the level structure $(E,l)$ modifies $\mathcal{O}^{\virt}_{g,m,d}$ by $\det^{-l}(V)$ (see Section \ref{LevelStructures}), which may formally be re-written as $\det^{-l}(V) = $
\[
\left.\left((-1)^{\rank V} \cdot \frac{\Lambda^*_{-\mu}V^\vee}{\Lambda^*_{-\mu^{-1}}V}\right)^l\right|_{\mu=1} = (-1)^{l\rank V} \cdot \left[\exp{\left(- l\sum_{k>0} \mu^k \frac{\Psi^{-k}(V)}{k}\right)} \cdot \exp{\left(l\sum_{k>0} \mu^{-k} \frac{\Psi^{k}(V)}{k}\right)}\middle]\right|_{\mu = 1}.
\]
Here $\Lambda^*_{-\mu}$ refers to the $\mu$-weighted alternate sum of exterior powers and $\Psi^k$ the $k$-th Adams operation. With the appearance of $\mu$ (and $\mu^{-1}$), the two exponential factors on RHS are both invertible K-theoretic characteristic classes of $V$. Taking the limit at $\mu = 1$ recovers the level structure. As we will see, when $l=1$, the first factor comes exactly from twisting of the type $(\Eu,E)$; the second factor, of the type $(\Eu,E)^\vee$ according to our notation later in this paper, is closer when $E = TX$ to the one considered in Liu \cite{Liu} which produces balanced ``vertex functions''. In fact, the KqSD implies that the effect on the Givental cones of the $(\Eu,E)^\vee$-twisting is similar to that of the $(\Eu^{-1},E^\vee)$-twisting. 

\par In order to take proper limit of $\mu$ and thus realize the idea of composition above, however, we need a careful treatment of the convergence issues. Since the expression above involves both $\mu$ and $\mu^{-1}$, we encounter essentially a two-sided infinite sum which is ill-defined. 

\par Our solution involves an overhaul of the loop space formalism, where we extend the domain of definition of the big $\J$-functions, generating functions of genus-zero K-theoretic GW invariants whose images are exactly the Givental cones in consideration, to allow inputs being not only Laurent polynomials in $q$ but also certain rational functions in $q$ with poles away from roots of unity. It is the main goal of Section \ref{LoopSpacesGeneralized} (after necessary preparation in Section \ref{Preliminaries} regarding the classical theory) to understand such extension as well as several of its variations, called the rational loop spaces, and study their relations. Section \ref{LoopSpacesGeneralized} is also where the main technicality of this paper resides. As we will see, the extended domain of definition is somehow more intrinsic to quantum K-theory. We expect the strategy of studying such rational loop spaces and their reduction to have broader applications in genus-zero quantum K-theory than the consideration in this paper of level structures and quantum Serre duality.

\par Combining the new loop space formalism with the $D_q$-module structure on the Givental cones, we arrive at interpreting the effect of the two factors on RHS (of the formula of $\det^{-l}(V)$) above in a way without reference to power series expansion. Moreover, it acquires a well-defined limit at $\mu=1$. In this way, we provide a simpler proof to the result in \cite{R-Z} (see Theorem \ref{JfuncLevel3}) which does not require working directly with the adelic characterization \cite{Givental:perm9}. The details are carried out in Section \ref{LevelStructuresRevisited}. Such re-interpretation enables us to study the effect of the level structures, in genus-zero case, under the more classical scenario of (invertibly) twisted K-theoretic GW theory (in the sense of \cite{Tonita,Givental:perm11} as mentioned above). 

\par Eventually, we prove the KqSD in Section \ref{QuantumSerreDuality}, under the assumption that $E$ is primitive in $K(X)$ (that is to say, $E = \sum_j \pm f_j \in K(X)$ for line bundles $f_j$). The idea is still to express the level structures in terms of invertible twistings, so the proof resembles largely those in the preceding section, except that we will no longer need to take the limit at $\mu=1$ in the end, as the fiber-wise $\C^\times$-action on $E$ is indispensable here. The deduction will be done in several different variations of the generalized loop spaces defined in Section \ref{LoopSpacesGeneralized}, leading to different versions of the KqSD. When $E$ is convex, the appearance of big $\J$-functions suggest indeed that a non-equivariant limit at $\mu=1$ should exist, perhaps in the context of quantum $D$-modules \cite{IMM}, but we do not discuss such limit in the present paper. 

\subsection{Torus-equivariant theory and the assumption on primitivity}

\par In Section \ref{TorusEquivariantTheory}, we explore how the assumption of $E$ being primitive in the KqSD may be relaxed over (partial) flag varieties via consideration of the torus-equivariant theory. For such purpose, we set up the rational loop space formalism further under torus-equivariant settings, and derive a torus-equivariant version of the KqSD using the recursive characterization developed in Givental \cite{Givental:perm2}. We give the precise statements in Section \ref{TorusQuantumSerre}. 

\par When $X = \text{Flag}(v_1,\cdots,v_n;N)$, the torus-equivariant KqSD on its associated abelian quotient $Y$, where all bundles are primitive, reduces to the non-torus-equivariant KqSD for any given non-primitive bundle $E$ over $X$. Here we write out the statements only for the case where $E = V_j \ (1\leq j\leq n)$ is the $j$-th tautological bundle over $X$. They are corollaries of Theorem \ref{TheoremEquivariantQuantumSerre}.
\begin{manualcorollary}{A}[KqSD for flag varieties] \label{CorollaryA}
For any point 
\[
J^1 = \sum_{d = (d_{is}),d_{is}\geq 0} \ \prod_{i=1}^n (Q_i')^{\sum_{s=1}^{v_i}d_{is}} \cdot J^1_d \quad \in \L^{X,(\Eu,V_j,l)}
\]
with $J^1_d$ independent of $Q$ (or equivalently, $Q'$), we have
\[
J^2 = \Eu(\mu^{-1}V_j) \cdot \sum_{d = (d_{is}),d_{is}\geq 0} \ \prod_{i=1}^n \left( Q_i' \cdot \left(-q^{-l\cdot \delta_{i,j}}\right)\right)^{\sum_{s=1}^{v_i}d_{is}} \cdot J^1_d \quad \in \L^{X,(\Eu^{-1},V_j^\vee,l+1)}.
\]
\end{manualcorollary}
\begin{manualcorollary}{B}[KqSD for flag varieties] \label{CorollaryB}
For any point 
\[
J^1 = \sum_{d = (d_{is}),d_{is}\geq 0} \ \prod_{i=1}^n (Q_i'')^{\sum_{s=1}^{v_i}d_{is}} \cdot J^1_d \quad \in \L^{X,(\Eu,E,l)}
\]
with $J^1_d$ independent of $Q$ (or equivalently, $Q''$), we have
\[
J^2 = \Eu(\mu V_j^\vee) \cdot \sum_{d = (d_{is}),d_{is}\geq 0} \ \prod_{i=1}^n \left( Q_i'' \cdot \left(-q^{-(l+1)\cdot \delta_{i,j}}\right)\right)^{\sum_{s=1}^{v_i}d_{is}} \cdot J^1_d
\quad \in \L^{X,(\Eu^{-1},E^\vee,l+1)}.
\]
\end{manualcorollary}

\subsection{Acknowledgement}
I owe much thanks to my PhD advisor Alexander Givental for suggesting this problem and for providing many beautiful insights. This material is based upon work supported by the National Science Foundation under Grant DMS-1906326 and by ERC Consolidator Grant ROGW-864919. 

\section{Preliminaries} \label{Preliminaries}
\par Throughout the paper unless otherwise stated, we assume $X$ is a smooth projective variety. We fix line bundles $P_1,\cdots,P_n$ such that $\{c_1(P_i)\}_{i=1}^n$ form a basis of $H^2(X)$. Then, we may express any integral curve class $d\in H_2(X)$ in terms of a tuple of integers $(d_i)_{i=1}^n$ where $d_i = - \langle c_1(P_i), d\rangle$. For simplicity, we assume that the effective cone of curve classes is contained in $(\mathbb{Z}_{\geq 0})^n\subset \mathbb{Z}^n$. We also fix an additive basis $\{\phi_\alpha\}_{\alpha\in A}$ of $K(X)$. We denote by $\{\phi^\alpha\}_{\alpha\in A}$ its dual basis under the Poincaré pairing on $K(X)$:
\[
\langle \phi_\alpha,\phi^{\alpha'}\rangle := \chi(X; \phi_\alpha\otimes\phi^{\alpha'}) = \delta_{\alpha}^{\alpha'}.
\]

\subsection{Big $\mathcal{J}$-function} \label{bigJ}
\par We review the basics of permutation-invariant quantum K-theory below. For a more complete introduction, see for example Givental \cite{Givental:perm1, Givental:perm7}. 

\par Let $\M_{g,m}(X,d)$ be the moduli space of genus-$g$ stable maps to $X$ with homological degree $d$ and $m$ marked points. We denote by $\ev_i: \M_{g,m}(X,d) \rightarrow X\ (1\leq i\leq m)$ the evaluation map at the $i$-th marked point, sending the stable map $f:(C;p_1,\cdots,p_m)\rightarrow X$ to the value $f(p_i)\in X$. We denote by $L_i\ (1\leq i\leq m)$ the universal cotangent line bundle at the $i$-th marked point. Then, the K-theoretic permutation-invariant Gromov-Witten invariants (or ``correlators'') are defined as
\[
\langle aL_1^{k}, \cdots, aL_m^{k}\rangle^{S_m}_{g,m,d} := \chi^{S_m}\left(\widebar{\mathcal{M}}_{g,m}(X,d); \mathcal{O}^{\virt}_{g,m,d} \otimes \bigotimes_{i=1}^m\ev_i^*(a)L_i^{k}\right),
\]
for $a\in K(X)$ and $k\in\mathbb{Z}$. Here $\mathcal{O}^{\virt}_{g,m,d}$ is the virtual structure sheaf on $\M_{g,m}(X,d)$ introduced by Lee \cite{Lee}, a K-theoretic analogue of the virtual fundamental cycle from the usual Gromov-Witten theory. $\chi^{S_m}$ denotes the $S_m$-invariant part of the virtual holomorphic Euler characteristic, which is naturally an $S_m$-representation. In fact, $S_m$ acts on $\M_{g,m}(X,d)$ by permutation of marked points, and the bundle given above is invariant under pull-back along such action because 
\begin{itemize}
    \item $\ev_i\ (1\leq i\leq m)$ and $L_i\ (1\leq i\leq m)$ are both compatible to permutation of marked points;
    \item we have the same input $a$ and $k$ at all marked points.
\end{itemize} 
Equivalently, one may directly take the holomorphic Euler characteristic over $[\M_{g,m}(X,d)/S_m]$. 

\par When not all of $a$ and $k$ are the same, however, the virtual holomorphic Euler characteristic is no longer representation with respect to the entire $S_m$-action, but only action of a subgroup $H\subset S_m$. In such cases, given partition $m = m_1+\cdots+m_s$ and $H = S_{m_1}\times\cdots\times S_{m_s}\subset S_m$, we may define correlators of the form
\[
\langle a_1L_1^{k_1}, \cdots, a_1L_{m_1}^{k_1}, a_2L_{m_1+1}^{k_2}, \cdots, a_2L_{m_1+m_2}^{k_2},\cdots, a_{s}L_m^{k_s}\rangle_{g,m,d}^H
\]
for any $a_1,\cdots,a_s\in K(X)$ and $k_1,\cdots,k_s\in\mathbb{Z}$. Moreover, it is not hard to see that the same construction applies not only for monomial inputs, but for Laurent polynomials in $L_i$ as well. In other words, for Laurent polynomials $\t_1(q),\cdots,\t_s(q)$ with coefficients in $K(X)$, we may define correlators of the form
\[
\langle \t_1(L_1), \cdots, \t_1(L_{m_1}), \t_2(L_{m_1+1}), \cdots, \t_2(L_{m_1+m_2}),\cdots, \t_s(L_m)\rangle_{g,m,d}^H
\]

\par From now, we focus on the genus-zero theory. In Gromov-Witten-type theories, genus-zero invariants are recorded in the generating function commonly known as the big $\J$-function. By definition, the big $\J$-function of the permutation-invariant quantum K-theory of $X$ is
\[
\mathcal{J}^X(\t) = 1 - q + \t(q) + \sum_{d_i\geq 0,m\geq 0,\alpha\in A} \ \prod_{i=1}^n Q_i^{d_i} \cdot \phi^\alpha \langle \frac{\phi_\alpha}{1-qL_0}, \t(L_1), \cdots, \t(L_m)\rangle^{S_m}_{0,m+1,d},
\]
where $Q_i\ (1\leq i\leq n)$ are formal variables called the Novikov variables keeping track of the degrees of curves, and $\t = \t(q) = \sum_{k}\t_kq^k \in K(X)\otimes\Lambda[q,q^{-1}]$ is a Laurent polynomial. Here $\Lambda$ is a base coefficient ring containing the Novikov variables. We take the invariant part with respect to $H = S_m = S_1\times S_m\subset S_{m+1}$ for correlators in $\J^X$ with $m+1$ marked points, as only the last $m$ marked points have the same input $\t$.

\par Two remarks are in place regarding the definition above of the big $\J$-function. First, the expression $\phi_\alpha / (1-qL_0)$ is not a Laurent polynomial in $L_0$, so its meaning is still unclear. We postpone the details until Section \ref{LoopSpacesGeneralized}, and only remark here that the idea is to view it as certain Taylor expansion in $L_0$. The expansion will be different on different ``Kawasaki strata'' of the inertia stack $I[\M_{0,m+1}(X,d)/S_m]$, but will always have only finite terms, and will have coefficients being rational functions in $q$ with poles at roots of unity. 

\par Second, one may notice that we have only defined the correlators for input $\t = \t(q)$ being a Laurent monomial in $q$ with coefficient in $K(X)$, but in the big $\J$-function we need coefficients in $K(X)\otimes\Lambda$, which are $\Lambda$-linear combinations of elements in $K(X)$. Such extension to $\Lambda$-linear combinations is, however, not direct due to the $S_m$-action. Essentially, we need to decompose the $S_m$-invariant dimension as the average of $\tr_h$ over all $h\in S_m$, and eventually the result relies on the so-called $\lambda$-algebra structure on $\Lambda$. For motivation and details of such construction, the readers are referred to \cite{Givental:perm1}.

\par A $\lambda$-algebra is an algebra incorporated with a series of endomorphisms $\Psi^k\ (k\in\mathbb{Z}_{\geq 0})$ called the Adams operations. For the simplest case of $\Lambda = \C[[Q_1,\cdots,Q_n]]$, they are defined by $\Psi^k(Q_i) = Q_i^{k}\ (\forall 1\leq i\leq n)$. We note that for this $\Lambda$, there exists an ideal $\Lambda_+ = (Q_1,\cdots,Q_n) \subset \Lambda$, with respect to the adic topology of which $\Lambda$ is complete. The ideal is compatible with the Adams operation in the sense that $\Psi^k(\Lambda_+)\subset\Lambda_+$. Moreover, $K(X)$ naturally admits a structure of $\lambda$-algebra as well, through $\Psi^k(P) = P^k$ for any line bundle $P\in K(X)$. 

\par For convergence issues, we will consider only $\Lambda_+$-small inputs $\t = \t(q)\in K(X)\otimes \Lambda_+[q,q^{-1}]$ for our big $\J$-function, and in such case $\J^X(\t)\in K(X)\otimes\Lambda(q)$ is well-defined as a rational function in $q$ with coefficients in $K(X)[[Q_1,\cdots,Q_n]]$, or alternatively, a formal power series in $Q_1,\cdots,Q_n$ with coefficients being $K(X)$-valued rational functions in $q$.  We will come back to the topic of convergence in more detail.

\par Later in this paper, we will very often need to extend the pair $(\Lambda,\Lambda_+)$ above to a larger $\lambda$-algebra with chosen ideal, usually by adding in certain equivariant parameters. We always assume that properties above remain unchanged
\begin{itemize}
    \item the Adams operations on $\Lambda$ still satisfy $\Psi^k(Q_i) = Q_i^{k}\ (\forall 1\leq i\leq n)$;
    \item there is still an ideal $\Lambda_+\subset\Lambda$ with $\Psi^k(\Lambda_+)\subset\Lambda_+$ and $Q_i\in\Lambda_+\ (1\leq i\leq n)$, such that $\Lambda$ is complete with respect to the $\Lambda_+$-adic topology. 
\end{itemize}

\subsection{Loop space formalism and pseudo-finite-difference (PFD) operators} \label{LoopSpacePFD}
\par The geometry of $\J^X(\t)$ is better understood under the loop space formalism, which we introduce below. For simplicity of notation, we consider the minimal coefficient ring $\Lambda = \C[[Q_1,\cdots,Q_n]]$ here. Other cases may be treated in the same way.

\par Let $\mathcal{K} =  K(X)\otimes\Lambda(q)$ be the space of $K(X)\otimes\Lambda$-valued rational functions in $q$. It admits an symplectic (i.e. non-degenerate anti-symmetric) pairing over $\Lambda = \C[[Q_1,\cdots,Q_n]]$
\[
\Omega(\f,\g) = \Res_{q\neq 0,\infty}\langle \f(q^{-1}), \g(q)\rangle \frac{dq}{q},
\]
where $\langle\cdot,\cdot\rangle$ is the Poincaré pairing on $K(X)$ defined earlier. $\mathcal{K}$ admits a Lagrangian polarization $\mathcal{K} = \mathcal{K}_+ \oplus \mathcal{K}_-$ where the two Lagrangian subspaces are given by
\[
\mathcal{K}_+ = K(X) \otimes \Lambda[q,q^{-1}],\quad 
\mathcal{K}_- = \{\f\in\mathcal{K}|\f(0)\neq\infty, \f(\infty) = 0\}.
\]
More precisely, $\mathcal{K}_-$ contain reduced rational functions (i.e. the degree of numerator is strictly less than the degree of denominator) in $q$ with no pole at $q=0$. Under such polarization, the input $\t = \t(q)$ of the big $\J$-function is an element in $\K_+$, while the value $\J^X(\t)$ is an element in $\mathcal{K}$ whose projection to $\mathcal{K}_+$ is exactly $1 - q + \t(q)$. 

\par We denote by $\L^X\subset\mathcal{K}$ the image of $\J^X$. It can be proven \cite{Givental:perm7} that $\L^X$ resides in an overruled cone through the origin. In fact, by our assumption on $\t(q)\in\mathcal{K}_+$ being $\Lambda_+$-small, $\L^X$ is exactly the $\Lambda_+$-small germ of the overruled cone at $1-q$. Moreover, $\L^X$ admits a $D_q$-module structure, a property we will use throughout the paper. We formulate this property as the lemma below, in terms of the operators $P_iq^{Q_i\partial_{Q_i}}\ (1\leq i\leq n)$ which act on $\mathcal{K}$ by $P_iq^{Q_i\partial_{Q_i}} \cdot \prod_{j=1}^nQ_j^{d_j} = P_iq^{d_i}\cdot \prod_{j=1}^nQ_j^{d_j}$ and trivially on terms without Novikov variables.
\begin{lemma} \label{PFD}
$\L^X$ is invariant under operators of the following two forms
\[
O_1 = \exp{\left(D(Pq^{Q\partial_Q}, q)\right)},\quad
O_2 = \exp{\left(\sum_{k>0}\frac{\Psi^k(D(Pq^{kQ\partial_Q}, q))}{k(1-q^k)}\right)},
\]
where $D$ is any $\Lambda_+$-small Laurent polynomial expression.
\end{lemma}
Here we use the abbreviated notation
\[
Pq^{Q\partial_Q} = (P_1q^{Q_1\partial_{Q_1}}, \cdots, P_nq^{Q_n\partial_{Q_n}})
\]
and the Adams operations act by
\[
\Psi^k(Q_i) = Q_i^{k}, \quad \Psi^k(q) = q^k, \quad \Psi^k(P_iq^{kQ_i\partial_{Q_i}}) = P_i^kq^{kQ_i\partial_{Q_i}}.
\]
We call operators like $O_1$ and $O_2$ the pseudo-finite-difference (PFD) operators, and denote by $\mathcal{P}$ the group generated by them. $\L^X$ is invariant under the action of $\mathcal{P}$. A proof of the lemma may be found in \cite{Givental:perm8}.

\subsection{Invertible twisting} \label{EulerTwisting}
\par We review in this section the twisted K-theoretic Gromov-Witten invariants as considered in \cite{Tonita,Givental:perm11}. We call such twisting invertible as it is determined by an invertible K-theoretic characteristic class $c$, given in the exponential form
\[
c(V) = \exp{\left(\sum_k c_k \cdot \Psi^k(V) \right)}
\]
by constants $c_k$. Note that in K-theory, the Adams operation $\Psi^k$ may also be defined for $k\in\mathbb{Z}_{\leq 0}$, simply as $\Psi^{-k}$ of the dual. Fixing such $c$ and an element $E$ in $K(X)$, we may define over $\M_{g,m}(X,d)$ a new virtual structure sheaf
\[
\mathcal{O}^{\text{virt},(c,E)}_{g,m,d} := \mathcal{O}^{\text{virt}}_{g,m,d} \otimes c(E_{g,m,d}) = \mathcal{O}^{\text{virt}}_{g,m,d} \otimes \exp{\left(\sum_k \frac{c_k\cdot \Psi^k(E_{g,m,d})}{k}\right)},
\]
where $E_{g,m,d} = \ft_{m+1,*}\ev_{m+1}^* E$ is an element in the K-group of $\M_{g,m}(X,d)$ defined through the forgetting map $\ft_{m+1}$ and evaluation map $\ev_{m+1}$ appearing in the diagram below.
\begin{equation*}
\xymatrix{
	\overline{\mathcal{M}}_{g,m+1}(X,d)\ \ar^{\quad\quad \ev_{m+1}}[r]\ar^{\ft_{m+1}}[d] & \ X \\
	\overline{\mathcal{M}}_{g,m}(X,d)\ & 
}
\end{equation*}
With the modified virtual structure sheaves, the same construction as in Section \ref{bigJ} gives rise to the $(c,E)$-twisted correlators $\langle\cdot\rangle^{H,(c,E)}_{g,m,d}$ for any subgroup $H = S_{m_1}\times\cdots\times S_{m_s}\subset S_m$, and thus the $(c,E)$-twisted big $\J$-function, where $\{\phi^\alpha\}_{\alpha\in A}$ is now taken as the dual basis to $\{\phi_\alpha\}_{\alpha\in A}$ of $K(X)$ under the $(c,E)$-twisted Poincaré pairing
\[
\langle \phi_\alpha,\phi^{\alpha'}\rangle^{(c,E)} := \chi(X; \phi_\alpha\otimes\phi^{\alpha'}\otimes c(E)) = \delta_{\alpha}^{\alpha'}.
\]
We denote by $\J^{X,(c,E)}$ the twisted big $\J$-function, and by $\L^{X,(c,E)}$ its image. $\L^{X,(c,E)}$ is still an overruled cone in $\K$, and is invariant under the action of PFD operators. For a more complete treatment, see \cite{Givental:perm11}.

\par One special case is what we call the Euler-type twisting. It appears originally in the context of quantum Lefschetz theorems, which express the quantum K-theory on complete intersections in $X$ in terms of the theory on $X$. The Euler-type twisting is also indispensable for our statement of the KqSD.

\par To this end, we endow $E$ with a auxiliary fiber-wise $\C^\times$-action with character $\mu^{-1}$, and consider $c = \Eu$ the equivariant Euler class. For convenience of notation, we sometimes write $\mu^{-1}E$ instead of $E$ to emphasize the $\C^\times$-action. By definition, the modification applied to the virtual structure sheaf $\mathcal{O}^{\text{virt}}_{g,m,d}$, in the $(\Eu,E)$-twisted theory, is exactly $\Eu(\mu^{-1} E_{g,m,d})$, where for any element $V-W$ in the K-group of $\M_{g,m}(X,d)$ with $V$ and $W$ both vector bundles, the equivariant Euler class is given by 
\[
\Eu(V - W) = \frac{\sum_{i} (-1)^i \Lambda^i V^\vee}{\sum_{i} (-1)^i \Lambda^i W^\vee}.
\]
$\mu$ is introduced here solely for the invertibility and convergence purposes, as we now explain. We specialize to the genus-zero theory. When $E$ is nef, the fiber of $E_{0,m,d}$ over the stable map $f:(C;p_1,\cdots,p_m)\rightarrow X$ is simply $H^0(C;f^*E)$, so $\Eu(\mu^{-1}E_{0,m,d})$ is a well-defined element in $K(\M_{0,m}(X,d))[\mu]$. When $E$ is not nef, however, the fiber looks like $H^0(C;f^*E) - H^1(C;f^*E)$, and thus $\Eu(\mu^{-1}E_{0,m,d})$ lives no longer in $K(\M_{0,m}(X,d))[\mu]$ but rather in its fraction field, which creates problem in computation.

\par One of the most natural ways of resolving this problem is to consider Taylor expansion in $\mu$ so that $\Eu(\mu^{-1}E_{0,m,d})$ becomes a well-defined element in $K(\M_{0,m}(X,d))[[\mu]]$. In this way, for the $(\Eu,E)$-twisted theory, we consider 
\begin{itemize}
    \item the enlarged coefficient ring $\Lambda = \C[[Q_1,\cdots,Q_n,\mu]]$ with $\Psi^k(\mu) = \mu^{k}$,
    \item and the adic topology given by the ideal $\Lambda_+ = (Q_1,\cdots,Q_n,\mu)$.
\end{itemize}
In this case, it is not hard to see that $\J^{X,(\Eu,E)}$ is well-defined for any $\t(q)\in K(X)\otimes\Lambda_+[q,q^{-1}]$, and the value $\J^{X,(\Eu,E)}(\t)$ is well-defined as a formal power series in $Q_1,\cdots,Q_n$ and $\mu$ with coefficients being $K(X)$-valued rational functions in $q$, located in the $\Lambda_+$-small neighborhood of $1-q$ in $\mathcal{K}$. In other words, $\L^{X,(\Eu,E)}$ is a $\Lambda_+$-germ at $1-q$ for the enlarged pair $(\Lambda,\Lambda_+)$ above. Now, expanding $\Eu(\mu^{-1}E_{g,m,d})$ with respect to $\mu$, we have
\[
\Eu(\mu^{-1} \cdot E_{g,m,d}) = \exp{\left( - \sum_{k>0} \frac{\mu^k \cdot \Psi^{-k}(E_{g,m,d})}{k} \right)}.
\]

\par Through the quantum Adams-Riemann-Roch (qARR) formula developed in \cite{Givental:perm11}, one may express rather explicitly the $(c,E)$-twisted image cones $\L^{X,(c,E)}$ in terms of the untwisted cone $\L^{X}$. Applying directly the qARR formula, we obtain
\begin{proposition} \label{Prop1}
The $(\Eu,E)$-twisted big $\J$-function has image cone $\L^{X,(\Eu,E)} = \Box^{(\Eu,E)} \cdot \L^X$, where 
\[
\Box^{(\Eu,E)} = \exp{\left( - \sum_{k>0} \frac{\mu^{k} \cdot q^k \Psi^{-k}(E)}{k(1-q^k)} \right)}
\]
Both cones are regarded as in the loop space $\K$ with the enlarged coefficient ring $\Lambda = \C[[Q_1,\cdots,Q_n,\mu]]$. 
\end{proposition}

\par Assume now $E$ (and thus $E^\vee$) is generated by line bundles in $K(X)$. In other words, 
\[
E^\vee = \sum_{j=1}^{r} \pm\ f_j = \sum_{j=1}^{r} \pm\ f_j(P_1,\cdots,P_n) \quad\in K(X)
\]
for line bundles $f_j\ (1\leq j\leq \rank E)$ which are monomials in $P_1,\cdots,P_n$. We assume for simplicity of notation there are only positive terms, and thus $r = \rank E$. By Lemma \ref{PFD}, we have $\L^{X,(\Eu,E)} = D^{(\Eu,E)} \cdot \L^X$, where
\[
D^{(\Eu,E)} = \exp{\left( - \sum_{j=1}^{\rank E}\sum_{k>0} \frac{\mu^{k} \cdot q^k (f_j(P_1^k,\cdots,P_n^k) - f_j(P_1^kq^{kQ_1\partial_{Q_1}},\cdots,P_n^kq^{kQ_n\partial_{Q_n}}))}{k(1-q^k)} \right)}
\]
differs from $\Box^{(\Eu,E)}$ only by a PFD operator. Despite of its intimidating form, $D^{(\Eu,E)}$ acts on Novikov variables in an enjoyable way:
\begin{align*}
D^{(\Eu,E)} \cdot Q^d & = \exp{\left( - \sum_{j=1}^{\rank E}\sum_{k>0} \frac{\mu^{k} \cdot q^k (1 - q^{-k \langle c_1(f_j),d\rangle}) \cdot f_j(P_1^k,\cdots,P_n^k) }{k(1-q^k)} \right)} \cdot Q^d \\
& = \exp{\left( - \sum_{j=1}^{\rank E}\ \sum_{s=1}^{ - \langle c_1(f_j),d\rangle}\sum_{k>0} \frac{\mu^{k} \cdot q^{ks} \cdot f_j(P_1^k,\cdots,P_n^k) }{k} \right)} \cdot Q^d \\
& = \prod_{j=1}^{\rank E}\ \prod_{s=1}^{ - \langle c_1(f_j),d\rangle} \left( 1 - \mu \cdot q^{s} \cdot f_j \right) \cdot Q^d.
\end{align*}
If a certain summand $f_j$ comes with a negative sign in the decomposition of $E^\vee$, we simply need to invert the factors involving $f_j$. Therefore, we obtain the following proposition, equivalent to Proposition \ref{Prop1}. 
\begin{proposition} \label{Prop2}
For any point $J = \sum_{d} Q^d \cdot J_d \in \L^X$ with $J_d\in K(X)[[\mu]](q)$, we have
\[
\sum_d Q^d \cdot J_d \cdot \prod_{j=1}^{\rank E}\ \prod_{s=1}^{-\langle c_1(f_j),d\rangle} \left( 1 - \mu \cdot q^{s} \cdot f_j \right)^{\pm 1} \quad \in \L^{X,(\Eu,E)}
\]
\end{proposition}
Note that $\prod_{s=1}^{-\langle c_1(f_j),d\rangle}$ here refers more precisely to $\prod_{s=-\infty}^{-\langle c_1(f_j),d\rangle} / \prod_{s=-\infty}^{0}$, and thus resides actually in the denominator if $\langle c_1(f_j),d\rangle>0$. Therefore, unless $E$ may be decomposed into a sum of nef bundles in $K(X)$, terms involving $\mu$ emerge in the denominator, and they should always be regarded as its Taylor expansion at $\mu = 0$. Such interpretation will lead to a problem in our treatment of the KqSD, which is exactly the reason why the ``classical'' loop formalism that we have been using by far will not be sufficient. We will explain the details later.

\par For later use, we also consider modified virtual structure sheaves of the following form, which we call $(\Eu,E)^\vee$-twisted, closely related to the $(\Eu,E)$-twisted version above:
\[
\mathcal{O}^{\text{virt},(\Eu,E)^\vee}_{g,m,d} := \mathcal{O}^{\text{virt}}_{g,m,d} \otimes \Eu^{-1}((\mu^{-1}E_{g,m,d})^\vee).
\]
By applying the qARR formula to its $\mu^{-1}$-expansion
\[
\Eu^{-1}((\mu^{-1}E_{g,m,d})^\vee) = \exp{\left( \sum_{k>0} \frac{\mu^{-k} \cdot \Psi^{k}(E_{g,m,d})}{k} \right)}
\]
and pairing with suitable PFD operators as above, we obtain the following proposition.
\begin{proposition} \label{Prop3}
For any point $J = \sum_{d} Q^d \cdot J_d \in \L^X$ with $J_d\in K(X)[[\mu^{-1}]](q)$, we have
\[
\sum_d \ \frac{Q^d \cdot J_d} {\prod_{j=1}^{\rank E}\ \prod_{s=1}^{-\langle c_1(f_j),d\rangle} \left( 1 - \mu^{-1} \cdot q^{-s} \cdot f_j^{-1} \right)^{\pm 1}} \quad \in \L^{X,(\Eu,E)^\vee}.
\]
\end{proposition}
The equality holds in the loop space with the enlarged coefficient ring $\Lambda = \C[[Q_1,\cdots,Q_n,\mu^{-1}]]$ this time, and all $\mu^{-1}$-terms appearing in the denominator should be regarded as its Taylor expansion in $\mu^{-1}$.

\subsection{Level structures} \label{LevelStructures}
\par Level structures are introduced by Ruan and Zhang in \cite{R-Z}, where a determinant-type modification is applied to the virtual structure sheaves and the corresponding quantum K-theory with level structures is considered. Such twisting arises only in quantum K-theory but not in its cohomological prototype. As we will see later in this paper, level structures play an essential role in the K-theoretic quantum Serre duality. 

\par Given vector bundle $E$ over $X$ and integer $l$, we may consider the modified virtual structure sheaves
\[
\mathcal{O}^{\text{virt},(E,l)}_{g,m,d} := \mathcal{O}^{\text{virt}}_{g,m,d} \otimes \det^{-l}(E_{g,m,d}),
\]
where by definition $\det(V-W) = \det(V) \otimes \det(W)^{-1}$. The generating function $\J^{X,(E,l)}$ of correlators defined in the same way as in Section \ref{bigJ} but using virtual structure sheaves so modified, is called the big $\J$-function with level structure $(E,l)$. We denote its image cone in the loop space by $\L^{X,(E,l)}$. Note that similar to the case of Euler-type twisting of Section \ref{EulerTwisting}, the dual basis vectors $\{\phi^\alpha\}_{\alpha\in A}$ used in the definition of $\J^{(E,l)}$ are now taken with respect to the correspondingly modified Poincaré pairing on $K(X)$
\[
\langle \phi_\alpha,\phi^{\alpha'}\rangle^{(E,l)} := \chi(X; \phi_\alpha\otimes\phi^{\alpha'}\otimes \det^{-l}(E)) = \delta_{\alpha}^{\alpha'}.
\]

\par According to \cite{R-Z}, $\L^{X(E,l)}$ is related to $\L^X$ in the following way. We still assume the decomposition in $K(X)$ of $E^\vee$ into line bundles has no negative terms
\[
E^\vee = \sum_{j=1}^{\rank E} f_j = \sum_{j=1}^{\rank E} f_j(P_1,\cdots,P_n)
\]
and recognize that we invert the corresponding factors if any $f_j$ is endowed with a negative sign.
\begin{proposition} \label{Prop4}
For any point $J = \sum_{d} Q^d \cdot J_d \in \L^X$ with $J_d\in K(X)(q)$, we have
\[
\sum_d Q^d \cdot J_d \cdot \prod_{j=1}^{\rank E} \left[ f_j^{-\langle c_1(f_j),d\rangle} q^{\frac{\langle c_1(f_j),d\rangle(\langle c_1(f_j),d\rangle - 1)}{2}} \right]^{l} \in \L^{X,(E,l)}
\]
\end{proposition}
As an easy corollary, $\L^{X,(E,l)}$ and $\L^{X,(E^\vee,l)}$ differ only by the change of variables $Q_i\mapsto Q_i\cdot q^{c_1(E)_i} \ (1\leq i\leq n)$, where $c_1(E)_i\ (1\leq i\leq n)$ are the integer coefficients when we express $c_1(E)\in H^2(X)$ as a linear combination of the basis vectors $c_1(P_i)\ (1\leq i\leq n)$ of $H^2(X)$. 

\par We end this section with another (attempt of) interpretation of the level structures. Indeed, we do not need \emph{a priori} the auxiliary equivariant parameter $\mu$ for the consideration of level structures. However, the resemblance between the formulae in Proposition \ref{Prop2}, \ref{Prop3}, and \ref{Prop4} suggests certain intrinsic relation. Indeed, as we have mentioned in Introduction, for any element $V$ in the K-theory, 
\[
\det^{-1}(\mu^{-1} V) = (-1)^{\rank V} \cdot \frac{\Eu(\mu^{-1}V)}{\Eu(\mu V^\vee)}.
\]
Taking $V = E_{g,m,d}$, it is not hard to see that the twisting of the $(E,l)$-level structure on the virtual structure sheaves is, up to a sign, exactly the composition of the $(\Eu,E)$-twisting and the $(\Eu,E)^\vee$-twisting above, specialized at $\mu=1$. The sign here is, for genus-zero case,
\[
(-1)^{\rank E_{0,m,d}} = (-1)^{\rank E + \langle c_1(E),d \rangle}.
\]

\par In algebraic terms, to elucidate the relation even more, we may further re-write the expression appearing in Proposition \ref{Prop4} by
\[
f_j^{-\langle c_1(f_j),d\rangle} q^{\frac{\langle c_1(f_j),d\rangle(\langle c_1(f_j),d\rangle - 1)}{2}} = \left.(-1)^{-\langle c_1(f_j),d\rangle} \cdot \frac{\prod_{s=1}^{-\langle c_1(f_j),d\rangle} \left( 1 - \mu \cdot q^{s} \cdot f_j \right)}{\prod_{s=1}^{-\langle c_1(f_j),d\rangle} \left( 1 - \mu^{-1} \cdot q^{-s} \cdot f_j^{-1} \right)}\right|_{\mu=1}.
\]
The numerator and denominator on RHS are exactly what have appeared in Proposition \ref{Prop2} and \ref{Prop3}, created by the $(\Eu,E)$- and $(\Eu,E)^\vee$-twisting respectively; the sign here differs from earlier only by $(-1)^{\rank E}$, which may exactly be absorbed into the discrepancy of the correspondingly twisted Poincaré pairings on $K(X)$.

\par The deduction above is inviting as it seems to provide us with a strategy of realizing the level structures in terms of (the composition of) the more well-studied invertible twisting. However, we should note that the computation is so far only formal, as Proposition \ref{Prop2} requires expansion in $\mu$ while Proposition \ref{Prop3} in $\mu^{-1}$. In such way, the realization of the level structure as their composition will result in a two-sided infinite sum, which no longer has a proper meaning. 

\par In order to remedy this situation, we would need a formalism that does not demand expansion but rather treats the rational functions involving $\mu$ as themselves. We develop such a formalism in Section \ref{LoopSpacesGeneralized}, and carry out the intuitive re-interpretation of level structures above in rigorous terms in Section \ref{LevelStructuresRevisited}.

\section{Loop Spaces Generalized} \label{LoopSpacesGeneralized}
\par For the enlarged coefficient ring $\Lambda = \C[[Q_1,\cdots,Q_n,\lambda]]$ with Adams operation $\Psi^k(\lambda) = \lambda^{k}$ and ideal $\Lambda_{+} = (Q_1,\cdots,Q_n,\lambda)$, the (classical) loop space formalism introduced in Section \ref{Preliminaries} involves the space $\K$ consisting of rational functions in $q$ taking values in $K(X)[[Q_1,\cdots,Q_n,\lambda]]$. It admits a polarization given by the two Lagrangian subspaces $\K_{\pm}$, consisting of Laurent polynomials and reduced rational functions in $q$ respectively, and the $\Lambda_+$-small elements in $\K_+$ may serve as input of the big $\J$-function. The settings of Proposition \ref{Prop2} and \ref{Prop3} may be recovered by taking $\lambda$ either equal to $\mu$ or $\mu^{-1}$.

\par Nevertheless, as we have explained in previous sections, considering merely enlarged coefficient rings does not suffice for our purposes, due to the emergence of two-sided infinite sums. Therefore, we explore in this section a different way of extending the loop space formalism, by allowing directly rational functions of both $q$ and $\lambda$ in $\K_+$ (and thus as inputs to the big $\J$-function) instead of treating them as their Taylor expansion in $\lambda$. The central definition is given in Section \ref{CentralDefinition}. Yet, we start by considering the intermediate step where inputs are allowed to be more general rational functions in $q$ than merely Laurent polynomials, where much technicality arises already. We explore in Section \ref{PowerSeriesExpansion} the relation between the two spaces.

\par The construction in Section \ref{CentralDefinition} may be further generalized, either by increasing the number of parameters $\lambda$ or by posing more subtle constraints on the allowed denominators of the rational inputs. We give the most general form in Section \ref{variations}, and provide a dictionary of special cases that will be used later.

\par The symplectic form on the loop spaces may also be modified, but the construction need not to be changed as long as such form comes from a non-degenerate symmetric pairing on $K(X)$. For such reason, we proceed with the construction for now assuming the standard symplectic form, and postpone making the choice of symplectic form until Section \ref{variations}.

\par To distinguish the extended loop spaces defined in this Section from the classical ones where only Laurent polynomials are allowed as inputs, we call the new ones the rational loop spaces. (The author apologizes for his lack of creativity in giving names and is open to suggestions of improvements.)

\par Before moving on to the constructions, we first take a closer look at the big $\J$-function, or more precisely the part $\J^X(\t)-(1-q+\t(q))$, i.e. the part in the big $\J$-function coming from the correlators. Its $q$-dependence comes completely from the input $\phi_\alpha/(1-qL_0)$ at the $0$-th marked point, where $L_0$ is the $0$-th universal cotangent bundle. By Kawasaki-Riemann-Roch formula \cite{Kawasaki,Givental-Tonita},
\begin{align}
& \left\langle \frac{\phi_\alpha}{1-qL_0}, \t(L_1),\cdots,\t(L_m) \right\rangle^{S_m}_{0,m+1,d} \\
= & \chi\left([\M_{0,m+1}(X,d)/S_m]; \mathcal{O}^{\virt}_{0,m+1,d} \otimes \frac{\ev_0^* \phi_\alpha}{1-qL_0}\otimes \bigotimes_{i=1}^m \ev_i^*\t(L_i)\right) \\
= & \sum_{C} \ \int_{[C]^{\virt}} \ch\left(\tr_{g}\frac{1}{1-qL_0|_C}\right) \cdot \ch(\tr_{g} \ev_0^* \phi_\alpha) \cdot \ch\left(\tr_{g}\bigotimes_{i=1}^m \ev_i^*\t(L_i)|_C\right) \cdot \mathfrak{T}(\tr_{g} N_{C},T_C), \label{KawasakiRR}
\end{align}
where
\begin{itemize}
    \item the sum is over connected components $C$ of the inertia stack $I[\M_{0,m+1}(X,d)/S_m]$;
    \item each such component $C$ is by definition labeled by a conjugacy class of the group of local isotropy fixing each stable map in it, and $g$ is any representative of the conjugacy class;
    \item $\mathfrak{T}$ depends on the normal and the tangent bundle of $C$ and is explicitly known.
\end{itemize}
Since $g$ acts trivially on $C$, $\tr_{g} L_0|_C = \zeta\cdot \overline{L_0}$ for a constant $\zeta\in\C$ (the fiber-wise eigenvalue of $g$) and a non-equivariant line bundle $\overline{L_0}$. $\zeta$ has to be a root of unity as $g$ has finite order. We may thus expand
\[
\ch\left(\tr_{g}\frac{1}{1-qL_0|_C}\right) = \ch\left(\frac{1}{1-q\zeta \overline{L_0}}\right) = \ch\left(\frac{1}{1-q\zeta - q\zeta(\overline{L_0}-1)} \right) = \sum_{s=1}^\infty \frac{(q\zeta)^s}{(1-q\zeta)^{s+1}}\cdot \ch(\overline{L_0}-1)^s,
\]
and the expansion terminates after finite terms, as $\ch(\overline{L_0}-1)$ is nilpotent for dimension reasons. This is exactly how we understand the input at the $0$-th marked point.

\par Moreover, it is not hard to see that only poles at $q = \zeta^{-1}$ with $\zeta$ a root of unity arise (though possibly with multiplicities) in $\J^X(\t)-(1-q+\t(q))$, and the $q$-degree of the numerator is always less than that of the denominator. In other words, reduced rational functions in $q$ with poles solely at roots of unity should always be contained in $\K_-$.

\subsection{Rational functions in $q$ as inputs} \label{RationalInputs}
\par In this section, we move a small step forward from the usual loop space formalism by allowing rational functions in $q$ with no poles at roots of unity as inputs $\t = \t(q)$ of $\J^X(\t)$, which amounts to keeping the same $\K$ but taking a larger $\K_+$. We work over the coefficient ring $\Lambda = \C[Q,\lambda] := \C[[Q_1,\cdots,Q_n,\lambda]]$ and its ideal $\Lambda_+ = (Q,\lambda)$ to illustrate the idea, but the construction works for other $\Lambda$ and $\Lambda_+$ as well. To emphasize the dependence on $\lambda$, we denote the loop space by $\K[[\lambda]]$ and its Lagrangian subspaces by $\K_{\pm}[[\lambda]]$. 
\begin{definition}[Loop space with rational functions in $q$ allowed as input]
\begin{align*}
\K_+[[\lambda]] & := K(X)\otimes \Lambda \otimes \left\{ \frac{a(q)}{b(q)}\in\C(q)\ \middle|\ (a,b)=1, \ b \text{ has no zero at roots of unity} \right\}, \\
\K_-[[\lambda]] & := K(X)\otimes \Lambda \otimes \left\{ \frac{a(q)}{b(q)}\in\C(q)\ \middle|\ \deg a < \deg b, \ b \text{ has only zeros at roots of unity} \right\}, \\
\K[[\lambda]] \ & := K(X)\otimes \Lambda \otimes \ \C(q),
\end{align*}
where the tensor products are all over $\C$.
\end{definition}
We use the convention that $\deg c = 0$ for non-zero constant $c$ and $\deg 0 = -\infty$; $(\cdot,\cdot)$ means taking the greatest common divisor. Roughly speaking, $\K$ contains rational functions in $q$ taking values in power series of $Q = (Q_1,\cdots,Q_n)$ and $\lambda$; $\K_-$ contains those with only poles at roots of unity; $\K_+$ contains those without any poles at roots of unity.
\begin{proposition}\label{PolarizationLambdaExpanded}
$\K_{\pm}[[\lambda]]$ give a Lagrangian polarization of $\K[[\lambda]]$ under the $\Lambda$-bilinear symplectic pairing
\[
\Omega(\f,\g) = \Res_{q = \text{roots of unity}}\ \langle \f(q^{-1}), \g(q)\rangle \cdot \frac{dq}{q},
\]
for the base coefficient ring $\Lambda = \C[[Q,\lambda]]$. Here $\langle\cdot,\cdot\rangle$ is the K-theoretic Poincaré pairing on $K(X)$.
\end{proposition}
The pairing is anti-symmetric because $\langle\cdot,\cdot\rangle$ is symmetric but $dq/q$ admits a minus sign when we take change of variable $q\mapsto q^{-1}$. 

\par We start by proving the two subspaces are both Lagrangian. Any element in $\K[[\lambda]]$ can be regarded as an infinite vector of rational functions in $\C(q)$, through expansion in power series in $Q_1,\cdots,Q_n$ and $\lambda$ and then under the basis $\{\phi_\alpha\}$ of $K(X)$. Since $\Omega$ is $\Lambda$-bilinear and $\langle\cdot,\cdot\rangle$ on $K(X)$ is non-degenerate, it suffices to consider the $\C(q)$-entries and show that the two corresponding subspaces of $\C(q)$ are Lagrangian with respect to the pairing 
\[
\omega(f,g) = \Res_{q = \text{roots of unity}} f(q^{-1})\cdot g(q) \cdot \frac{dq}{q}.
\]
\begin{itemize}
    \item For $f,g$ both with poles only away from roots of unity: the poles of $f(q^{-1})$ are either at reciprocals of the poles of $f(q)$ or at $q=0$, and thus all poles of $f(q^{-1})\cdot g(q)\cdot q^{-1}$ are away from roots of unity. The residue computed by $\omega(f,g)$ thus vanishes naturally.
    \item For $f,g$ both reduced rational functions with poles only at roots of unity: the poles of $f(q^{-1})\cdot q^{-1}$ are still only at roots of unity for degree reasons, and as a consequence $h(q) = f(q^{-1})\cdot q^{-1}\cdot g(q)$ has no pole away from roots of unity. The residue computed by $\omega(f,g)$ thus vanishes by the residue theorem:
    \[
    \Res_{q = \text{roots of unity}} h(q) \cdot dq = - \Res_{q \neq \text{roots of unity}} h(q) \cdot dq.
    \]
\end{itemize}
That the pairing $\omega$ (and thus $\Omega$) is non-degenerate may be proved using the same idea. We omit the details.

\par That $\K[[\lambda]] = \K_{+}[[\lambda]]\oplus \K_{-}[[\lambda]]$ may be reduced to an argument on $\C(q)$ as well. It then follows from a standard deduction using the Euclidean division on $\C[q]$. We postpone the details until the proof of Proposition \ref{PolarizationCentral}, which will use exactly the same idea and will be under the even more complicated settings of Section \ref{CentralDefinition}. This completes our proof of Proposition \ref{PolarizationLambdaExpanded}.

\par We devote the rest of this section to defining $\J^X(\t)$ for any $\Lambda_+$-small input $\t = \t(q)\in\K_+[[\lambda]]$. For
\[
\t = \t(q) = \sum_{d\in\mathbb{Z}^n_{\geq 0}} \sum_{i\in\mathbb{Z}_{\geq 0}} \sum_{\alpha'\in A} Q^d \lambda^i \phi^{\alpha'} \cdot \frac{a_{d,i,\alpha'}(q)}{b_{d,i,\alpha'}(q)} \quad\in \K_{+}[[\lambda]],
\]
where $a_{d,i,\alpha'}(q), b_{d,i,\alpha'}(q)\in\C[q]$ and $b_{d,i,\alpha'}(q)$ has no pole at roots of unity for any $d,i,\alpha'$, we first define its contribution to the individual correlator of $\J^X(\t)$
\[
\left\langle \frac{\phi_\alpha}{1-qL_0}, \t(L_1),\cdots,\t(L_m) \right\rangle^{S_m}_{0,m+1,d}.
\]

\par The idea is once again to use the Kawasaki-Riemann-Roch formula, i.e. Formula (\ref{KawasakiRR}), and it suffices to understand the contribution of $\t(L_j)$ to $\tr_{g}(\otimes_i \ev_i^*\t(L_i)|_C)$ for each Kawasaki stratum $C$ of the inertia stack $I[\M_{0,m+1}(X,d)/S_m]$. As we have mentioned, each stratum $C$ comes with a fixed conjugacy class (in which we take a representative $g$) of the local isotropy group, and all points in the component are fixed by $g$. Since $S_m$ permutes the marked points, however, for a stable map $f:(\Sigma;p_0,\cdots,p_m)\rightarrow X$ on $C$, $g$ represents an automorphism of $\Sigma$ that commutes with $f$, fixes $p_0$, but possibly permutes $p_1,\cdots,p_m$. Therefore, unlike the input at the $p_0$ that we have discussed at the beginning of this section, $\ev_j^*\t(L_j)$ is not preserved by the action of $g$ but rather permuted, unless $p_j$ happens to be fixed by $g$, and thus deserves a different treatment.

\par Given $1\leq j\leq m$, we assume that $p_j$ is contained in a cycle of length $l$ of the permutation in $S_m$ induced by $g$, and that the cycle looks like $(j_1,j_2,\cdots,j_l)$ with $j_1=j$. That is to say, $g(p_{j_k}) = p_{j_{k+1}}$ for $1\leq k\leq l-1$ and $g(p_{j_l}) = p_{j_1}$. $\langle g\rangle$ acts trivially on the component $C$, but only the subgroup $\langle g^l\rangle$ is lifted to an action on $L_j$. We need the following lemma to better understand the contribution of $\t(L_j)$ (and of all other $\t(L_{j_k})$).
\begin{lemma} \label{AdamsFromPermutation}
Let $M$ be a smooth projective variety on which $\langle g \rangle$ acts trivially, and let $V$ be an $\langle g^l \rangle$-equivariant vector bundle over $M$. Then,
\[
\tr_g \left( V^{\otimes l} \right) = \Psi^l(\tr_{g^l} V).
\]
where the action on LHS is given by
\[
g\cdot (v_1\otimes v_2\otimes\cdots\otimes v_l) = g^l(v_{l})\otimes v_1\otimes\cdots\otimes v_{l-1}.
\]
\end{lemma}
The special case of Lemma \ref{AdamsFromPermutation} where $\langle g^l\rangle$ acts trivially on $V$ is proven in, for instance, \cite{Givental:perm3} or \cite{Givental-Tonita}. Let $r = \rank V$, then the idea of proof is to reduce to the case of $V$ being the universal $\rank=r$ vector bundle over $BGL(r)$, and thus further to the case of $V$ being the standard representation of $GL(r)$ over $M = pt$, where the $\langle g \rangle$-action on $V$ commutes with that of $GL(r)$. The idea works equally well in our case here. In fact, since $\langle g^l\rangle$ acts trivially on the base, $K_{\langle g^l\rangle}(X) = \text{Rep}(\langle g^l\rangle)\otimes K(X)$, so any $V$ in our case may be decomposed into building blocks on which $\langle g^l\rangle$ acts with the same constant scaling by root of unity $\zeta$. Now, we consider the case of one such building block of rank $r$, which we denote still by $V$ by abuse of notation. Let $T\subset GL(r)$ be a maximal torus, and $\{\e_{i}\}_{i=1}^r$ be an eigenbasis of $V$ with respect to $T$ where $h\in T$ acts with eigenvalues $\{x_i\}_{i=1}^r$, then $\{\e_{i_1}\otimes \cdots\otimes \e_{i_l}\}_{i_1,\cdots,i_l=1}^r$ form a basis of $V^{\otimes l}$. Let $\langle\cdot,\cdot\rangle_V$ be the inner product such that $\{\e_{i}\}_{i=1}^r$ is orthonormal, then
\begin{align*}
\tr_{(g,h)} \left( V^{\otimes l} \right) & = \sum_{i_1,\cdots,i_l = 1}^r \langle (g,h) \cdot \e_{i_1}\otimes\cdots\otimes\e_{i_l}\ , \ \e_{i_1}\otimes\cdots\otimes\e_{i_l} \rangle_{V^{\otimes l}} \\
& = \sum_{i_1,\cdots,i_l = 1}^r x_{i_1}\cdots x_{i_l} \cdot \langle g^l(\e_{i_l})\otimes\cdots\otimes\e_{i_{l-1}}\ , \ \e_{i_1}\otimes\cdots\otimes\e_{i_l} \rangle_{V^{\otimes l}} \\
& = \sum_{i=1}^r x_i^l \cdot \langle g^l(\e_i)\otimes\cdots\otimes\e_{i}\ , \ \e_{i}\otimes\cdots\otimes\e_{i} \rangle_{V^{\otimes l}} \\
& = \sum_{i=1}^r x_i^l \cdot \langle g^l(\e_i)\ , \ \e_{i}\rangle_{V} =  \sum_{i=1}^r \zeta\cdot x_i^l = \tr_{h} \Psi^l(\tr_{g^l} V).
\end{align*}
This completes our proof of Lemma \ref{AdamsFromPermutation}.

\par We drop $Q^d \lambda^i$ for now and take $V$ as $\ev^*_{j}\phi^{\alpha'}\cdot a_{d,i,\alpha'}(L_j)/b_{d,i,\alpha'}(L_j)$ for fixed $d,i,\alpha'$ and $M$ as (the étale chart of) $C$. Indeed, $g$ permutes $\{L_{j_k}\}_{k=1}^l$ cyclically, and thus after $l$-times $L_j = L_{j_1}$ goes back to itself. Moreover, pull-back along $g$ identifies $\{L_{j_k}\}_{k=1}^l$ as line bundles. Under such identification,
\[
\tr_g \bigotimes_{k=1}^l \left(\ev^*_{j_{k}}\phi^{\alpha'}\cdot\frac{a_{d,i,\alpha'}(L_{j_k})}{b_{d,i,\alpha'}(L_{j_k})}\right) = \tr_g \left(\ev^*_{j}\phi^{\alpha'}\frac{a_{d,i,\alpha'}(L_j)}{b_{d,i,\alpha'}(L_j)}\right)^{\otimes l},
\]
where the action on the RHS is exactly of the form in Lemma \ref{AdamsFromPermutation}. Therefore, by Lemma \ref{AdamsFromPermutation} we have
\[
\tr_g \bigotimes_{k=1}^l \left(\ev^*_{j}\phi^{\alpha'}\cdot\frac{a_{d,i,\alpha'}(L_{j_k})}{b_{d,i,\alpha'}(L_{j_k})}\right) = \Psi^l\left(\tr_{g^l} \ev^*_{j}\phi^{\alpha'}\cdot \frac{a_{d,i,\alpha'}(L_j)}{b_{d,i,\alpha'}(L_j)}\right) = \Psi^l\left(\ev^*_{j}\phi^{\alpha'}\cdot\frac{a_{d,i,\alpha'}(\zeta \overline{L})}{b_{d,i,\alpha'}(\zeta \overline{L})}\right),
\]
where we use $\tr_{g^l} L_j = \zeta \overline{L}$, with $\zeta\in\C$ the eigenvalue of the fiber-wise $g^l$-action on $L_j$, and $\overline{L}$ the underlying non-equivariant of $L_j$ (i.e. with $g^l$-action removed). It indicates that $\zeta$ is a root of unity and (the Chern character of) $(\overline{L} - 1)$ is nilpotent. By our assumption, $b_{d,i,\alpha'}(q)$ has no poles at roots of unity, so $b_{d,i,\alpha'}(\zeta \overline{L})$ has no zero at $\overline{L} = 1$ in our case, which in turn means that $a_{d,i,\alpha'}(\zeta \overline{L})/b_{d,i,\alpha'}(\zeta \overline{L})$ admits a well-defined Taylor expansion in the nilpotent element $\overline{L}-1$ which terminates after finite terms. The expansion gives a well-defined element in the K-theory of the base locus, and thus $\Psi^l$ acts naturally on it. 

\par In this way, we develop a proper understanding of the term above. Note that this is compatible with the usual interpretation under the known cases of $a_{d,i,\alpha'}(q)/b_{d,i,\alpha'}(q)$ being Laurent polynomials in $q$. It remains to show that such definition is independent of the choice of $a_{d,i,\alpha'}$ and $b_{d,i,\alpha'}$. In fact, the only ambiguity comes from multiplying $a_{d,i,\alpha'}$ and $b_{d,i,\alpha'}$ by the same polynomial $c\in\C[q]$ where $c$ is also assumed to have no zeros at roots of unity. It then suffices to note that the $\overline{L}-1$ expansion of $c(\zeta \overline{L})$ and $1/c(\zeta \overline{L})$ are indeed inverse to each other.

\par Combining this result with the dropped term $Q^d\lambda^i$ and summing over all $d,i,\alpha'$, we obtain our definition of the contribution of $\t$ on the chosen stratum $C\subset I[\M_{0,m+1}(X,d)/S_m]$ through the marked points $p_{j}=p_{j_1},\cdots,p_{j_l}$ (cyclically permuted by the local isotropy $g$) as follows.
\[
\tr_g \bigotimes_{k=1}^l \ev_{j_k}^*\t(L_{j_k}) := \Psi^l(\tr_{g^l}\ev_j^*\t(L_j)) := \sum_{d\in\mathbb{Z}^n_{\geq 0}} \sum_{i\in\mathbb{Z}_{\geq 0}} \sum_{\alpha'\in A} Q^{ld} \lambda^{li} \cdot \ev^*_{j}\Psi^l(\phi^{\alpha'}) \cdot \frac{a_{d,i,\alpha'}\left(\zeta\overline{L}^l\right)} {b_{d,i,\alpha'}\left(\zeta\overline{L}^l\right)},
\]
where $Q^d\lambda^i$ as a coefficient from $\Lambda$ contributes by definition through $\Psi^l(Q^d\lambda^i) = Q^{ld}\lambda^{li}$ on any cycle of length $l$, and the $b_{d,i,\alpha'}(\zeta\overline{L}^l)$ in the denominator is understood as its (finite-term) expansion in $\overline{L}-1$. 

\par The contribution of $\t$ on $C$ through other marked points, grouped by cycles of permutation $g$, may be defined in a similar way. The total contribution of $\t$ on $C$ is then the tensor product of contribution from all such groups. This then completes our definition of the correlator $\langle \phi_\alpha/(1-qL_0), \t(L_1),\cdots,\t(L_m) \rangle^{S_m}_{0,m+1,d}$ for general $\t\in\K_+[[\lambda]]$ through Formula (\ref{KawasakiRR}), as the integration formula over $I[\M_{0,m+1}(X,d)/S_m]$ depends on $\t$ only through the Chern characters of $\tr_g$'s that we considered above.

\par In the end, we define $\J^X(\t) \in \K[[\lambda]]$ by adding together the correlators. Since $\t$ is assumed to be $\Lambda_+$-small, $a_{0,0,\alpha'}(q) = 0$ for all $\alpha'\in A$. In this way, the coefficient of the $Q^D \lambda^I$-term of $\J^X(\t)$, for each fixed $D\in\mathbb{Z}^n_{\geq 0}$ and $I\in\mathbb{Z}_{\geq 0}$, comes at most from finitely many different correlator terms, and is thus well-defined as element in $K(X)(q)$ (that is to say, no infinite sum arises in our summation process).

\par To emphasize the expansion in $\lambda$, we denote the big $\J$-function of this loop space by $\J^X[[\lambda]]$. It is not hard to see that the regular loop space $\K$ in Section \ref{LoopSpacePFD}, where only Laurent polynomials in $q$ are allowed, is a subspace of $\K[[\lambda]]$. Moreover, $\J^X[[\lambda]]$ reduces to the regular big $\J$-function $\J^X$ when restricted to $\K$.

\par For the coefficient ring $\Lambda = \C[[Q_1,\cdots,Q_n,\lambda_1,\cdots,\lambda_r]]$, the same construction as above carries over. We denote by $\K[[\lambda_1,\cdots,\lambda_r]]$ and $\J^X[[\lambda_1,\cdots,\lambda_r]]$ the corresponding loop space and big $\J$-function. 

\par For instance, in Section \ref{CompositionEulerTwisting}, we will encounter the case of $r=2$ but with the two parameters $\lambda_1$ and $\lambda_2$ replaced by $a$ and $b$ instead. $\K[[a,b]]$ is defined over the coefficient ring $\Lambda[[a,b]] := \C[[Q_1,\cdots,Q_r,a,b]]$ with $\Lambda_+[[a,b]] := (Q_1,\cdots,Q_r,a,b)$, and rational functions in $q$ are allowed in $\K_+[[a,b]]$ as long as they do not contain poles at roots of unity.

\subsection{Rational functions in both $q$ and parameters as inputs} \label{CentralDefinition}
\par Recall that in Section \ref{LevelStructures}, in order to avoid two-side infinite sum, we would like to treat rational functions in $\mu$ (and $\mu^{-1}$) directly as themselves instead of as their power series expansions. For such purpose, we upgrade further our loop space $\K$ by forbidding infinite series of $\lambda$, but including rational functions in both $\lambda$ and $q$ instead. We denote by $\K^\lambda$ the resulting generalized loop space, and by $\K^\lambda_{\pm}$ its two subspaces. It is the central definition of this paper.

\par Let $\C[q,\lambda]$ be the ring of polynomials in $q$ and $\lambda$. Then the following three subsets of $\C[q,\lambda]$ are closed under multiplication:
\begin{align*}
S_+ & := \{g(q,\lambda)\ | \ g(\zeta,0)\neq 0, \text{ for any root of unity } \zeta\}, \\
S_- & := \{f_0(q)\cdot f_1(\lambda)\ |\ f_0(q) = \prod_{s=1}^N (1-q\zeta_s) \text{ where each } \zeta_s \text{ is a root of unity; } f_1\in \C[\lambda], f_1(0)\neq 0\}, \\
S\ & := S_+\times S_- = \{g(q,\lambda)\cdot f(q,\lambda) \ | \ g\in S_+,f\in S_-\}.
\end{align*}
\begin{definition}[Loop space with rational functions in both $q$ and $\lambda$ allowed as input] \label{DefinitionCentralLoopSpace}
\begin{align*}
\K^\lambda & := K(X)\otimes \C[[Q]] \otimes S^{-1}\C[q,\lambda], \\
\K^\lambda_+ & := K(X)\otimes \C[[Q]] \otimes S_+^{-1}\C[q,\lambda], \\
\K^\lambda_- & := K(X)\otimes \C[[Q]] \otimes \left\{ \frac{a(q,\lambda)}{b(q,\lambda)} \in S_-^{-1}\C[q,\lambda] \ \middle| \ a,b\in\C[q,\lambda], \deg_q(a) < \deg_q(b) \right\}.
\end{align*}
The tensor products are all over $\C$, and as before $\C[[Q]]:= \C[[Q_1,\cdots,Q_n]]$.
\end{definition}
\par The motivation for such construction is two-fold:
\begin{itemize}
    \item We would like to relate $\K^\lambda$ to the loop space $\K[[\lambda]]$ of Section \ref{RationalInputs} (or even the more classical ones as in Section \ref{LoopSpacePFD} and \ref{EulerTwisting}) where power series, though not Laurent series, in the extra parameter $\lambda$ are considered. It leads to our construction of $S\subset \C[q,\lambda]\setminus (\lambda)$, where $(\lambda)$ is the principal ideal, such that elements in $\K^\lambda$ have no poles at $\lambda=0$. We will see in Section \ref{PowerSeriesExpansion} that power series expansion in $\lambda$ induces an inclusion from $\K^\lambda$ to $\K[[\lambda]]$.
    \item By the same idea as the previous section, in order that both itself and its expansion in $\lambda$ may serve as legitimate inputs to the big $\J$-functions, an element in $\K_+^\lambda$ should avoid any factors in its denominator that would give rise to poles of $q$ at roots of unity. It leads to our construction of $S_+$. 
\end{itemize}
The coefficient ring is now $\Lambda = \C[\lambda]_{(\lambda)}[[Q]]$, in the sense that $\mathcal{K}^\lambda$ and $\mathcal{K}_{\pm}^\lambda$ are all free $\Lambda$-modules. Here the subscript $(\lambda)$ means localization at the principal ideal $(\lambda)\subset\C[\lambda]$. The Adams operation on $\Lambda$ is defined by $\Psi^k(\lambda) = \lambda^k$. We consider the adic topology given by the ideal $\Lambda_+ = (Q) := (Q_1,\cdots,Q_n) \subset \Lambda$. In general, one may also consider the case where $\Lambda_+ = (Q,\lambda)$, but it will not be necessary for our purposes.

\par We still have the symplectic pairing and the Lagrangian polarization property on $\K^\lambda$.
\begin{proposition}\label{PolarizationCentral}
$\K_{\pm}^\lambda$ give a Lagrangian polarization of $\K^\lambda$ under the $\Lambda$-bilinear symplectic pairing
\[
\Omega(\f,\g) = \Res_{q = \text{roots of unity}}\ \langle \f(q^{-1}), \g(q)\rangle \cdot \frac{dq}{q},
\]
for the coefficient ring $\Lambda = \C[\lambda]_{(\lambda)}[[Q]]$. $\langle\cdot,\cdot\rangle$ is as usual the Poincaré pairing on $K(X)$.
\end{proposition}
It is worth mentioning that $\Omega$ takes values only in $\Lambda$, although \emph{a priori} it may land in the larger ring $\C(\lambda)[[Q]]$. The reason is that by construction of $\K^\lambda$, any of its element when expanded as a Laurent series in $q-\zeta$ at a root of unity $\zeta$, has Laurent coefficients that never exceed $\Lambda\otimes K(X)$.

\par Now we prove $\K^\lambda = \K^\lambda_+\oplus \K^\lambda_-$ as $\Lambda$-modules in full detail, as is promised in the proof of Proposition \ref{PolarizationLambdaExpanded}. It is not hard to see that $\K^\lambda_+$ and $\K^\lambda_-$ has trivial intersection, so it suffices to prove $\K^\lambda = \K^\lambda_+ + \K^\lambda_-$. Moreover, since we may always expand any element in $\K^\lambda$ with respect to a basis $\{\phi_\alpha\}_{\alpha\in A}$ of $K(X)$ and then into power series in $Q_1,\cdots,Q_n$, it suffices to prove the decomposition for the $S^{-1}\C[q,\lambda]$ part. Indeed, any element $I$ in $S^{-1}\C[q,\lambda]$ may be written as
\[
I = \frac{a(q,\lambda)}{f_0(q) \cdot f_1(\lambda) \cdot g(q,\lambda)},
\]
where $a(q,\lambda)\in\C[q,\lambda]$ and $f_0(q), f_1(\lambda),g(q,\lambda)$ are as described in the definition of $S_+$ and $S_-$. We take a factor $(1-\zeta q)$ of $f_0(q)$, where $\zeta$ is a root of unity, and do induction on $\deg f_0(q)$. Applying the Euclidean division to $g(q,\lambda)$ by $(1-\zeta q)$, 
\[
g(q, \lambda) = (1-\zeta q)\cdot Q(q,\lambda) + g(\zeta^{-1},\lambda)
\]
for $Q(q,\lambda)\in\C[q,\lambda]$ and $g(\zeta^{-1},\lambda)\in\C[\lambda]\setminus (\lambda)$ (because $g(\zeta^{-1},0)\neq 0$ by assumption on $g$). Rearranging the terms and multiplying by $a(q,\lambda) / f_1(\lambda)g(\zeta^{-1},\lambda)$, we have
\[
\frac{a(q,\lambda)}{f_1(\lambda)g(\zeta^{-1},\lambda)}\cdot g(q, \lambda) - \frac{Q(q,\lambda)a(q,\lambda)}{f_1(\lambda)g(\zeta^{-1},\lambda)}\cdot (1-\zeta q) = \frac{a(q,\lambda)}{f_1(\lambda)}.
\]
Applying the Euclidean division by $(1-\zeta q)$ again to the numerator $a(q,\lambda)$ in the coefficient of $g(q, \lambda)$, we may merge the quotient part into the coefficient of $(1-\zeta q)$ in the formula above, so only the remainder, constant in $q$, stays in the coefficient of $g(q, \lambda)$:
\[
\frac{a(\zeta^{-1},\lambda)}{f_1(\lambda)g(\zeta^{-1},\lambda)}\cdot g(q, \lambda) + \frac{A(q,\lambda)}{f_1(\lambda)g(\zeta^{-1},\lambda)}\cdot (1-\zeta q) = \frac{a(q,\lambda)}{f_1(\lambda)}.
\]
Note that all denominators appearing here are in $\C[\lambda]\setminus (\lambda)$. Therefore, the above formula, divided further by $f_0(q)\cdot g(q,\lambda)$, gives a decomposition of $I$
\[
\frac{a(\zeta^{-1},\lambda)}{f_1(\lambda)g(\zeta^{-1},\lambda) \cdot f_0(q)} + \frac{A(q,\lambda)}{F_0(q) \cdot F_1(\lambda) \cdot g(q,\lambda)} = I.
\]
Here the first term on the LHS is in $\K^\lambda_-$; in the second term $F_0(q) := f_0(q)/(1-\zeta q)\in\C[q]$ has only zeros at roots of unity but has one less degree than $f_0(q)$, and $F_1(\lambda) := f_1(\lambda)g(\zeta^{-1},\lambda) \in\C[\lambda]\setminus (\lambda)$. We may then repeat what we have done for $I$ above for the second term on the LHS, until the degree of $f_0(q)$ reduces to zero, in which case the second term on LHS will eventually be an element in $\K^\lambda_+$. 

\par The proof of $\Omega$ being symplectic and the proof of $\K_{\pm}^\lambda$ both being Lagrangian are very similar to that of Proposition \ref{PolarizationLambdaExpanded}, so we omit the details here. This finishes our proof of Proposition \ref{PolarizationCentral}.

\par The big $\J$-function may be defined over $\K^\lambda$ in roughly the same way as in Section \ref{RationalInputs}. It suffices to provide a proper interpretation of (the Chern character of) $\t(L_j)$ on any stratum $C\subset I[\M_{0,m+1}(X,d)/S_m]$ with prescribed local isotropy $g$. Recall that if $p_j$ lies in a cycle of length $l$ of the permutation induced by $g$, then by Lemma \ref{AdamsFromPermutation}, for any $\Lambda_+$-small input 
\[
\t = \t(q) = \sum_{d\in\mathbb{Z}^n_{\geq 0}} \sum_{\alpha\in A} Q^d \phi_{\alpha} \cdot \frac{a_{d,\alpha}(q,\lambda)}{b_{d,\alpha}(q,\lambda)} \quad\in \K_{+}^\lambda
\]
with $a_{d,\alpha}, b_{d,\alpha} \in\C[q,\lambda]$ and $b_{d,\alpha}(\zeta,0)\neq 0$, its contribution boils down to $\Psi^l(\tr_{g^l}\t(L_j))$. We define
\begin{align} \label{DefineTrgl}
\Psi^l(\tr_{g^l}\ev_j^*\t(L_j)) := \sum_{d\in\mathbb{Z}^n_{\geq 0}} \sum_{\alpha\in A} Q^{ld} \cdot \ev^*_{j}\Psi^l(\phi_{\alpha}) \cdot \frac{a_{d,\alpha}\left(\zeta\overline{L}^l, \lambda^l\right)} {b_{d,\alpha}\left(\zeta\overline{L}^l, \lambda^l\right)},
\end{align}
where $\zeta$ and $\overline{L}$ are as before, and $1 / b_{d,\alpha}(\zeta \overline{L}^l, \lambda^l)$ is understood as its finite-term Taylor expansion in the nilpotent element $\overline{L} - 1$. Due to our assumption on $b_{d,\alpha}$, such Taylor expansion in $\overline{L} - 1$ has coefficients all lying in $\C[\lambda]\setminus (\lambda)$. We take $\lambda^l$ instead of $\lambda$ in the definition of $\Psi^l(\tr_{g^l}\t(L_j))$ so that it is compatible with the Adams operation on the coefficient ring. In this way, taking holomorphic Euler characteristics, we see that any correlator 
\[
Q^d \phi^\alpha \left\langle \frac{\phi_\alpha}{1-qL_0}, \t(L_1),\cdots,\t(L_m) \right\rangle^{S_m}_{0,m+1,d}
\]
is indeed an element in $\K^\lambda_-$, with the poles of $q$ at roots of unity coming from the leading input $\phi_\alpha/(1-qL_0)$, and the poles at $\lambda$ away from zero coming from the (permuted) inputs $\t(L_j)\ (1\leq j\leq m)$.

\par Similar to before, the assumption of $\t$ being $\Lambda_+$-small indicates that $a_{0,\alpha}(q,\lambda) = 0$ for all $\alpha\in A$, i.e. there is no $Q^0$-term in $\t$. Hence, the coefficient of the $Q^D$-term of $\J^X(\t)$, for each fixed $D\in\mathbb{Z}^n_{\geq 0}$, comes at most from finitely many different correlator terms and is thus well-defined taking values $\C[\lambda]_{(\lambda)}\otimes K(X)$. It is also clear that the projection of $\J^X(\t)$ to $\K^\lambda_+$ is $(1-q)+\t(q)$ while all the remaining terms live in $\K^\lambda_-$. 

\par To avoid conflict of notation, we denote the big $\J$-function over $\K^\lambda$ by $\J^{X,\lambda} = \J^{X,\lambda}(\t)$.

\subsection{Power series expansion} \label{PowerSeriesExpansion}
\par In this section, we prove that the power series expansion in $\lambda$ defines an inclusion from $\K^\lambda$ to $\K[[\lambda]]$. We also point out how the big $\J$-functions (see Section \ref{bigJ} and \ref{LoopSpacePFD}) on the two spaces are related.

\par We start with the coefficient rings. It is not hard to see that the power series expansion in $\lambda$ gives an injective $\C$-algebra homomorphism $\iota: \Lambda^\lambda := \C[\lambda]_{(\lambda)}[[Q]]\rightarrow \Lambda[[\lambda]] := \C[[Q,\lambda]]$, which preserves the Adams operations and satisfies $\iota(\Lambda^\lambda_+)\subset \Lambda[[\lambda]]_+$. Recall that $\Lambda^\lambda_+ = (Q)\subset \Lambda^\lambda$ and $\Lambda[[\lambda]]_+ = (Q,\lambda)\subset \Lambda[[\lambda]]$. Here $Q$ always serves as the abbreviation for $Q_1,\cdots,Q_n$.

\par Then on the level of loop spaces, we have the following proposition.
\begin{proposition} \label{LoopSpaceHom}
The power series expansion in $\lambda$ gives an injective homomorphism of symplectic $\Lambda^\lambda$-modules 
\[
\iota:  \K^\lambda \rightarrow \K[[\lambda]],
\]
where the $\Lambda^\lambda$-module structure on the codomain, naturally a $\Lambda[[\lambda]]$-module, is induced from the coefficient ring homomorphism above. Under such homomorphism, $\iota(\K^\lambda_{\pm}) = \K[[\lambda]]_{\pm}$. Moreover, the restriction of $\iota$ on $\K^\lambda_+$, which we denote by $\iota_+$, satisfies the following properties.
\begin{itemize}
    \item[(a)] $\iota_+(1-q) = 1-q$.
    \item[(b)] $\iota_+$ commutes with the Adams operations on the coefficients. That is to say,
    \[
    \iota_+(\Psi^l(\t(q))) = \Psi^l(\iota_+(\t(q))), \quad \forall \t = \t(q)\in \K^\lambda_+, l\in\mathbb{Z}_{>0}.
    \]
    Here $\Psi^l(\lambda) = \lambda^l$, $\Psi^l(Q) = Q^l$, and $\Psi^l$ acts on $K(X)$ through the standard action.
    \item[(c)] $\iota_+$ commutes with the changes of variables $q\mapsto \zeta q^l$ for any root of unity $\zeta$ and $l\in\mathbb{Z}_{>0}$. That is to say,
    \[
    \iota_+(\t(\zeta q^l)) = \iota_+(\t)(\zeta q^l), \quad \forall \t = \t(q)\in \K^\lambda_+.
    \]
\end{itemize}
\end{proposition}
The statements regarding $\iota$ may be checked directly. We only remark that by construction of $S_+$ in Definition \ref{CentralDefinition}, elements in  $\K^\lambda_+$ have no poles at $(q,\lambda) = (\zeta,0)$ for any root of unity $\zeta$, which ensures that no poles of $q$ at roots of unity arise under $\lambda$-power series expansion of $\K^\lambda_+$. As for the the properties regarding $\iota_+$, (a) and (c) are obvious, as by definition $\iota_+$ acts only non-trivially on $\lambda$, and (b) follows from the fact that the Taylor expansion in $\lambda$ commutes with the change of variable $\lambda\mapsto\lambda^l$ for any $l\in\mathbb{Z}_{>0}$. 
\begin{proposition}[Identification of big $\J$-functions] \label{BigJCommutes}
\[
\iota\left( \J^{X,\lambda}(\t) \right) = \J^{X}[[\lambda]](\iota_+(\t)), \quad \forall \t = \t(q) \in \Lambda^\lambda_+ \cdot \K^\lambda_+.
\]
\end{proposition}
It follows from our definition of the big $\J$-functions. Indeed,
\begin{align*}
\mathcal{J}^{X,\lambda}(\t) & = 1 - q + \t(q) + \sum_{d,m,\alpha} \ Q^{d} \cdot \phi^\alpha \langle \frac{\phi_\alpha}{1-qL_0}, \t(L_1), \cdots, \t(L_n)\rangle^{S_m}_{0,m+1,d}, \\
\mathcal{J}^{X}[[\lambda]](\iota_+(\t)) & = 1 - q + \iota_+(\t)(q) + \sum_{d,m,\alpha} \ Q^{d} \cdot \phi^\alpha \langle \frac{\phi_\alpha}{1-qL_0}, \iota_+(\t)(L_1), \cdots, \iota_+(\t)(L_n)\rangle^{S_m}_{0,m+1,d},
\end{align*}
so it suffices to look at the individual correlators. By our argument in Section \ref{RationalInputs} and \ref{CentralDefinition} using Kawasaki-Riemann-Roch formula, $\t(L_j)$ and $\iota_+(\t)(L_j)$ contribute to the correlators only through, respectively,
\[
\Psi^l(\tr_g \ev_j^*\t(L_j)) = \Psi^l(\ev_j^*\t)(\zeta q^l)|_{q = \overline{L}}, \quad\text{and}\quad \Psi^l(\tr_g \ev_j^*\iota_+(\t)(L_j)) = \Psi^l(\ev_j^*\iota_+(\t))(\zeta q^l)|_{q = \overline{L}}.
\]
Here, the non-equivariant bundle $\overline{L}$, the positive integer $l$, and the root of unity $\zeta$ depend only on the Kawasaki stratum $C$ and are thus the same in both cases; $\Psi^l$ is the $l$-th Adams operation (on the coefficients); $\overline{L}-1$ is nilpotent, and the RHS's should both be understood as the finite-term power series expansion in $\overline{L}-1$, with the former having coefficients in $\Lambda^\lambda$ and the latter having coefficients in $\Lambda[[\lambda]]$. The two RHS's are identified under $\iota$ by (b) and (c) of Proposition \ref{BigJCommutes}, and thus the corresponding correlators as $\iota$ commutes with taking Chern characters and integration. Finally, since $\iota(\Lambda^\lambda_+)\subset \Lambda[[\lambda]]_+$, all legitimate inputs $\t$ of $\J^{X,\lambda}$ are sent to legitimate inputs $\iota_+(\t)$ of $\J^{X}[[\lambda]]$. This completes our proof.
\begin{corollary} \label{LCommutes}
Let $\L^{X,\lambda}\subset \K^{X,\lambda}$ and $\L^{X}[[\lambda]]\subset \K^{X}[[\lambda]]$ be the images of $\J^{X,\lambda}(\t)$ and $\J^{X}[[\lambda]](\t)$ respectively. Then,
\[
\L^{X,\lambda} = \iota^{-1}\left(\L^{X}[[\lambda]]\right) \ \bigcap \ ((1-q) + \Lambda^\lambda_+\cdot\K^\lambda).
\]
In other words, if a point $I\in\K^\lambda$, $\Lambda^\lambda_+$-close to $(1-q)$, satisfies $\iota(I)\in\L^{X}[[\lambda]]$, then $I\in \L^{X,\lambda}$.
\end{corollary}
In particular, taking $\lambda = \mu$ (or $\mu^{-1}$), we are finally able to regard the expression appearing in Proposition \ref{Prop2} (or \ref{Prop3}) as an honest rational function in $\lambda$ and $q$ rather than its Taylor expansion in $\lambda$, and the corollary tells us that it resides indeed on $\L^{X,\lambda}$. Re-interpretation of this kind will be very useful later.

\par For simplicity of notation, starting from now, we omit the decorations and write simply $\L^X\subset \K^{X,\lambda}$ and $\L^X\subset \K^{X}[[\lambda]]$ to distinguish between the image cones of big $\J$-functions defined in different loop spaces. We follow the same notational convention in the rest of the paper, as no confusion will be caused in this way.

\subsection{Variations} \label{variations}
\par The rational loop space $\K^\lambda$ in Section \ref{CentralDefinition} may be further generalized in several ways. We introduce the variations below, and list at the end a few special cases which will be used later in this paper.

\par First, we may generalize the construction by increasing the number of parameters. More precisely, we consider the coefficient ring $\Lambda = T^{-1}\C[\lambda_1,\cdots,\lambda_r] \otimes \C[[Q_1,\cdots,Q_n]]$, where
\[
T = \left\{\prod_{i=1}^r f_i(\lambda_i)\ \middle|\ f_i(0)\neq 0, 1\leq i\leq r\right\} = \prod_{i=1}^r \C[\lambda_i]\setminus (\lambda_i),
\]
and the ideal $\Lambda_+ = (Q_1,\cdots,Q_r)$. Below are three multiplicatively closed subsets of $\C[q,\lambda_1,\cdots,\lambda_r]$
\begin{align*}
S_+^{\lambda_1,\cdots,\lambda_r} & := \left\{\prod_{i=1}^r g_i(q,\lambda_i)\ \middle| \ g_i(\zeta,0)\neq 0, \text{ for any root of unity } \zeta \text{ and any } i\right\}, \\
S_-^{\lambda_1,\cdots,\lambda_r} & := \left\{f_0(q)\cdot \prod_{i=1}^r f_i(\lambda_i)\ \middle|\ f_0(q) = \prod_{s=1}^N (1-q\zeta_s) \text{ where each } \zeta_s \text{ is a root of unity; } \prod_{i=1}^r f_i(\lambda_i) \in T\right\}, \\
S^{\lambda_1,\cdots,\lambda_r} \ & := S_+^{\lambda_1,\cdots,\lambda_r} \times S_-^{\lambda_1,\cdots,\lambda_r}.
\end{align*}
We define 
\begin{align*}
\K^{\lambda_1,\cdots,\lambda_r} & = K(X) \otimes \C[[Q_1,\cdots,Q_r]] \otimes (S^{\lambda_1,\cdots,\lambda_r})^{-1}\C[q,\lambda_1,\cdots,\lambda_r], \\
\K^{\lambda_1,\cdots,\lambda_r}_{+} & = K(X) \otimes \C[[Q_1,\cdots,Q_r]] \otimes (S^{\lambda_1,\cdots,\lambda_r}_{+})^{-1} \C[q,\lambda_1,\cdots,\lambda_r], \\
\K^{\lambda_1,\cdots,\lambda_r}_{-} & = K(X) \otimes \C[[Q_1,\cdots,Q_r]] \otimes \left\{\frac{a(q)}{b(q)} \in (S^{\lambda_1,\cdots,\lambda_r}_{-})^{-1} \C[q,\lambda_1,\cdots,\lambda_r] \ \middle| \ \deg a < \deg b\right\}.
\end{align*}
Note that if we regard $T$ as a multiplicative subset of $\C[q,\lambda_1,\cdots,\lambda_r]$, it is a subset of both $S_+$ and $S_-$. In this way, $\K^{\lambda_1,\cdots,\lambda_r}$ and the two subspaces $\K^{\lambda_1,\cdots,\lambda_r}_{\pm}$ are all free $\Lambda$-modules. The definition above indicates:
\begin{itemize}
    \item When we expand $\Psi^l(\tr_{g^l}\t(L_j))$ at $\overline{L}-1$ as in Formula (\ref{DefineTrgl}) for $\t = \t(q) \in \K^{\lambda_1,\cdots,\lambda_r}_+$, the expansion coefficients lie in $\Lambda = T^{-1}\C[\lambda_1,\cdots,\lambda_r] \otimes \C[[Q_1,\cdots,Q_n]]$, as the denominators come entirely from $S_+^{\lambda_1,\cdots,\lambda_r}$. Consequently, the correlators appearing in $\J(\t)$ for such $\t$ give indeed elements in $\K^{\lambda_1,\cdots,\lambda_r}_-$, which in turn implies that the big $\J$-function is well-defined for $\Lambda_+$-small elements in $\K^{\lambda_1,\cdots,\lambda_r}_+$.
    \item Elements in $\K^{\lambda_1,\cdots,\lambda_r}$ all admit well-defined Taylor expansion in $\lambda_1,\cdots,\lambda_r$. Moreover, such expansion gives rise to an inclusion $\K^{\lambda_1,\cdots,\lambda_r} \rightarrow \K[[\lambda_1,\cdots,\lambda_r]]$, under which the two Lagrangian polarizations are compatible (expansion coefficients of elements in $\K^{\lambda_1,\cdots,\lambda_r}_+$ have no poles of $q$ at roots of unity). As a result, Proposition \ref{BigJCommutes} and Corollary \ref{LCommutes} both hold for this general version. We omit the statements.
\end{itemize}

\par Second, one should realize that the only purpose of the choice of $T$ above is to ensure that the Taylor expansion in $\lambda_1,\cdots,\lambda_r$ gives an inclusion $\K^{\lambda_1,\cdots,\lambda_r} \rightarrow \K[[\lambda_1,\cdots,\lambda_r]]$. When there is no requirement for such expansion, one may instead consider any subset 
\[
T\subset \C[\lambda_1,\cdots,\lambda_r],
\]
closed under multiplication, invariant under the Adams operations, and not containing zero. A loop space formalism over the coefficient ring $\Lambda = T^{-1}\C[\lambda_1,\cdots,\lambda_r]\otimes \C[[Q_1,\cdots,Q_n]]$ parallel to the one above may be developed as below. $\Lambda$ is complete under the adic topology of its ideal $\Lambda_+ = (Q_1,\cdots,Q_n)$.

\par Consider any three multiplicatively closed subsets $S^T_+,S^T_-,S^T$ of $\C[q,\lambda_1,\cdots,\lambda_r]$ satisfying
\begin{itemize}
    \item $S^T_{\min}:=\{q^k | k\in\mathbb{Z}_{\geq 0}\}\times T \ \subset S^T_+ \ \subset S^T_{\max}:=\{ g(q,\lambda_1,\cdots,\lambda_r)|g(\zeta,\lambda_1,\cdots,\lambda_r) \in T, \forall \text{ root of unity }\zeta \},$
    \item $S^T_- = \left\{f_0(q)\ \middle|\ f_0(q) = \prod_{s=1}^N (1-q\zeta_s) \text{ where each } \zeta_s \text{ is a root of unity } \right\} \times T,$
    \item $S^T = S^T_+ \times S^T_-$,
\end{itemize}
Note that while $S^T_-$ is determined entirely by $T$, the above construction allows for some flexibility on $S^T_+$ even for a fixed $T$. We define
\begin{align*}
\K & = K(X)\otimes \C[[Q_1,\cdots,Q_n]] \otimes (S^T)^{-1}\C[q, \lambda_1,\cdots,\lambda_r], \\
\K_{+} & = K(X)\otimes \C[[Q_1,\cdots,Q_n]] \otimes (S^T_{+})^{-1}\C[q, \lambda_1,\cdots,\lambda_r], \\
\K_{-} & = K(X)\otimes \C[[Q_1,\cdots,Q_n]] \otimes \left\{\frac{a(q)}{b(q)} \in (S^T_{-})^{-1}\C[q, \lambda_1,\cdots,\lambda_r] \ \middle| \ \deg a < \deg b \right\}.
\end{align*}
This variated loop space $\K$ is a free $\Lambda$-module, and $\K_{\pm}$ are free sub-modules. The big $\J$-function may be defined on $\K$ in the same way as in Section \ref{CentralDefinition}, taking $\Lambda_+$-small elements in $\K_+$ as input and giving correlators in $\K_-$ as output. In fact, we impose the constraints on $S^T_+$ as above exactly because
\begin{itemize}
    \item [(a)] we want $S^T_{\min} \subset S^T_+$ so that Laurent polynomials in $q$ with coefficients in $\Lambda$ are allowed in $\K_+$;
    \item[(b)] we need $S^T_+ \subset S^T_{\max}$ so that the expansion of $\Psi^l(\tr_{g^l}\t(L_j))$ at $\overline{L}-1$ as in Formula (\ref{DefineTrgl}) for $\t = \t(q) \in \K_+$ has coefficients lying in $\Lambda$, and thus the big $\J$-function is well-defined on $\Lambda_+\cdot\K_+$.
\end{itemize}
In particular, in our construction of $\K^{\lambda_1,\cdots,\lambda_r}$ above, $S^T_+$ was taken as $S^T_{\max}$ for the corresponding $T$. 

\par Moreover, $\K_{\pm}$ gives a Lagrangian polarization of $\K$ under the symplectic pairing 
\[
\Omega(\f,\g) = \Res_{q = \text{roots of unity}}\ \langle \f(q^{-1}), \g(q)\rangle \cdot \frac{dq}{q}
\]
for any non-degenerate symmetric pairing $\langle \cdot , \cdot\rangle$ on $K(X)$. For simplicity, we take it always as the Poincaré pairing on $K(X)$ unless stated otherwise (in certain cases we need to apply twisting to it). We omit the justification of the polarization being Lagrangian, as it largely resembles the one for $\K^\lambda$. We only remark here that $\K_+\cap\K_- = 0$ because $S^T_+\cap S^T_- = T$, and that $\Omega$ is non-degenerate because we have taken $S^T_+\supset S^T_{\min}$, which implies that Laurent polynomials in $q$ with coefficients in $\Lambda$ are included in $\K_+$.

\par So far we have described the most general settings under which a rational loop space in the spirit of Section \ref{CentralDefinition} may be defined. In summary, the entire construction is determined by
\begin{itemize}
    \item[(1)] a base coefficient ring $\Lambda$,
    \item[(2)] a set $S^T_+$ of ``allowed denominators of input'', and
    \item[(3)] a pairing $\langle \cdot , \cdot \rangle$ on $K(X)$.
\end{itemize}
For instance, when $\Lambda = T^{-1}\C[\lambda_1,\cdots,\lambda_r]$ with $T = \prod_{i} \C[\lambda_i]\setminus (\lambda_i)$ and $S^T_+ = S^T_{\max}$, we recover $\K^{\lambda_1,\cdots,\lambda_r}$ above. In this paper, we only need certain special cases which we list below. The rest of this section may be used as a dictionary notation to appear later.

\par We start with the trivial cases. When $r=0$, i.e. there are no parameters at all, we may consider the multiplicatively closed subset $T = \C\setminus 0 \subset \C$ and thus the base ring $\Lambda = T^{-1}\C \otimes \C[[Q]] = \C[[Q]]$. The choice of $S^T_+ = S^T_{\min}$ will give as back the classical loop space $\K$, where $\K_+$ is set of Laurent polynomials in $q$ taking values in $K(X)\otimes\C[[Q]]$. The choice of $S^T_+ = S^T_{\max}$ will give us a slightly larger loop space $\widetilde{\K}$, where $\widetilde{\K}_+$ contains now all rational functions in $q$ with poles away from the roots of unity, still taking values in $K(X)\otimes\C[[Q]]$. $\widetilde{\K}$ is exactly what we have considered in Section \ref{RationalInputs}, but with the extra parameter $\lambda$ removed. 

\par We will need the rational loop space $\K^{a,b}$ in Section \ref{CompositionEulerTwisting}, which is the special case of $\K^{\lambda_1,\cdots,\lambda_r}$ (introduced at the beginning of this section) when $r=2$, with the two parameters re-named $a$ and $b$ instead. We re-iterate that the Taylor expansion in $a$ and $b$ induces an inclusion of loop spaces $\iota^{a,b}: \K^{a,b}\rightarrow \K[[a,b]]$, where $\K[[a,b]]$ is constructed as in Section \ref{RationalInputs} but with base coefficient ring $\Lambda[[a,b]] = \C[[Q]]\otimes \C[[a,b]]$ and adic topology defined by the ideal $\Lambda_+ = (a,b,Q)$. $\iota^{a,b}$ share the same properties with $\iota$ studied in Section \ref{PowerSeriesExpansion}. In Section \ref{CompositionEulerTwisting}, the symplectic form on both $\K^{a,b}$ and $\K[[a,b]]$ is defined with respect to the $R(a,b)$-twisted Poincaré pairing on $K(X)$, namely by taking $\mathcal{F} = (-1)^{l\cdot \rank E} \cdot \Eu(a^{-1}E)^l / \Eu(b^{-1}E^\vee)^l$ in the following formula
\begin{align}
\langle V,W \rangle = \chi\left(X;V\otimes W\otimes \mathcal{F}\right) \label{TwistedPoincaréPairing}
\end{align}

\par We define $\K(\mu)$ as the rational loop space determined by $\Lambda = \C(\mu)[[Q]] = T^{-1}\C[\mu] \otimes \C[[Q]]$ with $T = \C[\mu]\setminus 0$, and $S^T_+ = S^T_{\max}$. In particular, $\K(\mu)$ contains all rational functions in $q$ and $\mu$, and $\K(\mu)_+$ those without factors $(1-q\zeta)$ in the denominator for any root of unity $\zeta$. $\K(\mu)$ will appear in Section \ref{ReductionCoefficients} and then in Section \ref{QuantumSerreDuality}. The symplectic form on $\K(\mu)$ in Section \ref{ReductionCoefficients} is induced from the Poincaré pairing on $K(X)$ with level structure $(\mu^{-1}E,l)$, and in Section \ref{QuantumSerreDuality} from the Poincaré pairing on $K(X)$ with either $(\Eu,E,l)$- or $(\Eu^{-1},E^\vee,l+1)$-twisting. More precisely, the three pairings are given respectively by taking $\mathcal{F} = \det^{-l}(\mu^{-1}E)$, $\mathcal{F} = \Eu(\mu^{-1}E)\otimes \det^{-l}(\mu^{-1}E)$ and $\mathcal{F} = \Eu^{-1}(\mu E^\vee) \otimes \det^{-(l+1)}(\mu E^\vee)$ in Formula (\ref{TwistedPoincaréPairing}).

\par We define $\K[\mu,\mu^{-1}]$ as the rational loop space determined by $\Lambda = \C[\mu,\mu^{-1}][[Q]] = T^{-1}\C[\mu]\otimes \C[[Q]]$ with $T = (\C\setminus 0)\times \{\mu^k | k\in\mathbb{Z}_{>0}\}$, and $S^T_+ = S^T_{\min}$. In particular, 
\begin{align*}
\K_+[\mu,\mu^{-1}] & := K(X)\otimes \C[[Q]] \otimes \C[q,q^{-1},\mu,\mu^{-1}], \\
\K_-[\mu,\mu^{-1}] & := K(X)\otimes \C[[Q]] \otimes \left\{ \frac{f(q,\mu)}{\mu^{\gamma}\cdot \prod_s (1-\zeta_s q)}\ \middle| \ f\in \C[q,\mu], \gamma\in\mathbb{Z}_{>0}, \zeta_s \text{ is root of unity}\right\}.
\end{align*}
$\K[\mu,\mu^{-1}]$ is a subspace of $\K(\mu)$, and will appear also in Section \ref{ReductionCoefficients} with symplectic form defined by the Poincaré pairing on $K(X)$ with level structure $(\mu^{-1}E,l)$ (see above). 

\par We define $\K(\mu)^{a,b}$ as the rational loop space determined by the coefficient ring $\Lambda = T^{-1}\C[\mu,a,b]\otimes \C[[Q]]$ with $T = (\C(\mu)\setminus 0) \times (\C[a]\setminus (a)) \times (\C[b]\setminus (b))$, and the set of allowed denominators of input 
\[
S^T_+ = \left\{ f_1(q,\mu)\cdot f_2(q,a)\cdot f_3(q,b) \ | \ f_1(\zeta,\mu)\cdot f_2(\zeta,a)\cdot f_3(\zeta,b) \in T, \forall \text{ root of unity } \zeta \right\}.
\]
It will appear in Section \ref{ProofQuantumSerre}. Intuitively, it is a combination of $\K(\mu)$ and $\K^{a,b}$ that we have already seen. Through Taylor expansion in $a$ and $b$, $\K(\mu)^{a,b}$ admits an inclusion into $\K(\mu)[[a,b]]$, the space defined in the same way as $\K(\mu)$ but with the coefficient ring further tensored with $\C[[a,b]]$ and the adic topology defined by the ideal $(Q,a,b)$. Such inclusion resembles the one studied in Section \ref{PowerSeriesExpansion}, and it is not hard to see that they enjoy similar properties.

\par Finally, given any torus $G$, we define a rational loop space $\widetilde{\K}^{pt}$ for point target space for use in Section \ref{TorusEquivariantLoopSpaceFormalism}. Let $\chi_1,\cdots,\chi_r$ be a basis of characters of $G$, then $\widetilde{\K}^{pt}$ is the rational loop space over determined by the coefficient ring $\Lambda = T^{-1}\C[\chi_1,\cdots,\chi_r] \otimes \C[[Q]] = \C(\chi_1,\cdots,\chi_r) \otimes \C[[Q]]$ and $S^T_+ = S^T_{\max}$. Here $\C(\chi_1,\cdots,\chi_r) = \Rep(G)_0$ is the field of fractions of $\Rep(G) = \C[\chi_1^{\pm 1},\cdots,\chi_r^{\pm 1}]$. Intuitively, $\widetilde{\K}^{pt}$ may be regarded as the ``$G$-equivariant'' loop space of $pt$, but is somehow trivial since only trivial action exists on $pt$. For more general target spaces $X$ equipped with $G$-action, the construction of the $G$-equivariant loop space requires a slightly different formalism, which we save until Section \ref{TorusEquivariantLoopSpaceFormalism}.

\par To end this section, we remark that when $\Lambda_+$ is larger than $(Q_1,\cdots,Q_n)$, $\Lambda$ will have to be completed correspondingly with respect to the adic topology and thus include certain power series in the parameters. In this way, we obtain a construction which intuitively combines that of Section \ref{RationalInputs} and that of the current section. $\K(\mu)[[a,b]]$ that we have seen above serves as an example. To avoid over-complicated notation, however, we will not systematically discuss such loop spaces in this paper. In this paper, such spaces will only arise as co-domains of Taylor expansion of certain parameters, which we have already understood.

\section{Level Structures Revisited} \label{LevelStructuresRevisited}
\par In this section, we realize into rigorous terms the intuitive idea (see Section \ref{LevelStructures}) of re-interpreting level structures via invertible twistings, with the help of rational loop spaces like $\K^\lambda$. In particular, combining the new loop space formalism with quantum Lefschetz theorems, we provide a new proof to Proposition \ref{Prop4}, which was originally proved in \cite{R-Z} using adelic characterization. During our proof, generalized rational loop spaces of Section \ref{variations} will appear as a result of certain ``reduction of coefficients'', i.e. base change of coefficient rings. We study how the big $\J$-functions are preserved under such reduction of coefficients, and prove a lemma in regard of this issue which applies to rather general settings. 

\par Throughout the section, we use the notation $\K$ with no decoration for the classical loop space, defined in Section \ref{LoopSpacePFD} with only Laurent polynomials in $q$ allowed as inputs in $\K_+$. It is over the minimal coefficient ring $\Lambda = \C[[Q_1,\cdots,Q_n]]$, and endowed with the adic topology by $\Lambda_+ = (Q) = (Q_1,\cdots,Q_n)$.

\subsection{Composition of invertible twistings} \label{CompositionEulerTwisting}

\par Recall that for vector bundle $E$ over $X$, the $(E,l)$-level structure modifies the virtual structure sheaf on $\M_{0,m+1}(X,d)$ by $\det^{-l}(E_{0,m+1,d})$. Consider the modification
\[
R(a,b)_{0,m+1,d} := (-1)^{l\cdot\rank E_{0,m+1,d}} \cdot \frac{\Eu(a^{-1}E_{0,m+1,d})^l}{\Eu(b^{-1}E_{0,m+1,d}^\vee)^l},
\]
where $E_{0,m+1,d}^\vee$ denotes the dual of $E_{0,m+1,d}$ and is different from $(E^\vee)_{0,m+1,d}$. Up to a sign, it comes from the composition of two invertible twistings, of the type $(\Eu,a^{-1}E)$ and $(\Eu,bE)^\vee$ respectively (see Section \ref{EulerTwisting}), and reduces to the $(E,l)$-level structure as $a=b=1$. Our idea is to recover $\L^{X,(E,l)}\subset \K$, the image cone with $(E,l)$-level structure, by studying reduction of coefficients $a,b$ on $\L^{X,R(a,b)}$. For now, we may understand $a^{-1}$ and $b$ as two independent characters of a $(\C^\times)^2$ acting fiber-wise on $E$. The process of taking $a = \mu$ and $b = \mu^{-1}$ below may be understood as taking a one-dimensional subtorus, and eventually the process of taking $\mu=1$ as taking the non-equivariant limit. 

\par We denote by $\J^{X,R(a,b)}$ the big $\J$-function defined using the modified virtual structure sheaves as above, and $\L^{X,R(a,b)}$ its image cone. We start by interpreting the twisting $R(a,b)$ through Taylor expansion in $a$ and $b$ as it allows us to apply the qARR formula. In this way, $R(a,b)_{0,m+1,d}\in K(\M_{0,m+1}(X,d))[[a,b]]$, and thus $\L^{X,R(a,b)}$ lives naturally in $\K[[a,b]]$ (see Section \ref{RationalInputs}) equipped with a symplectic form defined by the following twisted pairing on $K(X)$
\[
\langle V,W\rangle^{R(a,b)} = \chi\left(X;V\otimes W\otimes (-1)^{l\cdot \rank E} \cdot \frac{\Eu(a^{-1}E)^l}{\Eu(b^{-1}E^\vee)^l}\right).
\]
It reduces back to the Poincaré pairing with $(E,l)$-level structure when $a=b=1$. 

\par With loss of generality, we assume as usual that there are no negative terms in the decomposition of $E^\vee$
\[
E^\vee = \sum_{j=1}^{\rank E} f_j = \sum_{j=1}^{\rank E} \ f_j(P_1,\cdots,P_n) \quad\in K(X)
\]
where $f_j$ are line bundles. By the qARR formulae (Proposition \ref{Prop2} and \ref{Prop3}), we have
\begin{proposition} \label{JfunctionRab}
For any point $J = \sum_d Q^d\cdot J_d\in\mathcal{L}^X\subset \K$ with $J_d\in K(X)(q)$, we have
\[
J^{X,R(a,b)} := \sum_d Q^d \cdot J_d \cdot \prod_{j=1}^{\rank E} (-1)^{-\langle c_1(f_j),d\rangle} \cdot \frac{\prod_{s=1}^{-l\cdot \langle c_1(f_j),d\rangle} \left( 1 - a \cdot q^{s} \cdot f_j \right)^l}{\prod_{s=1}^{-\langle c_1(f_j),d\rangle} \left( 1 - b \cdot q^{-s} \cdot f_j^{-1} \right)^l} \quad\in \L^{X,R(a,b)}\subset \K[[a,b]],
\]
where $J^{X,R(a,b)}$ here is regarded as a power series in $a$ and $b$,
\end{proposition}
The sign $\prod_{j} (-1)^{-\langle c_1(f_j),d\rangle}$ comes from canceling $(-1)^{l\cdot\rank E_{0,m+1,d}}$ of $R(a,b)_{0,m+1,d}$ with $(-1)^{l\cdot \rank E}$ of the twisted Poincaré pairing. The factors involving $f_j$ should be inverted if $f_j$ admits a negative sign in $E^\vee$.

\par In order to obtain genuine rational functions and study the limit of $a=b=1$, we consider restriction of $\L^{X,R(a,b)}$ to the ``smaller'' loop space $\K^{a,b}$ (see Section \ref{variations}). Recall that $\K^{a,b}$ has base coefficient ring
\[
\Lambda^{a,b} = (T^{a,b})^{-1}\C[a,b]\otimes \C[[Q_1,\cdots,Q_n]], \quad \text{ where }\quad T^{a,b} = (\C[a]\setminus (a)) \times (\C[b]\setminus (b)),
\]
with adic topology given by the ideal $\Lambda^{a,b}_+ = (Q)$. For $\t\in\Lambda^{a,b}_+\cdot \K_+^{a,b}$, although we may understand the contribution of $\t(L_j)$ to any correlator in the same way as in Section \ref{variations}, in order to define $\J^{X,R(a,b)}(\t) \in \K^{a,b}$, we still need to obtain a proper understanding of the twisting sheaf $R(a,b)_{0,m+1,d}$ under the settings of $\K^{a,b}$. Since the correlators may be computed by the Kawasaki-Riemann-Roch formula as the sum of integrals over the Kawasaki strata of $I[\M_{0,m+1}(X,d)/S_m]$, it suffices to realize $\tr_g R(a,b)_{0,m+1,d}|_C$ as an element in $\Lambda^{a,b}\otimes K(C)$, for any given (étale chart of) Kawasaki stratum $C$ and its associated local isotropy $g$ acting trivially on $C$. We start with the following observations:
\begin{itemize}
    \item For any $\langle g \rangle$-equivariant vector bundle $V$ over $C$, $\tr_g \Eu(a^{-1}V)$ is well-defined and invertible in $\Lambda^{a,b}\otimes K(C)$. In fact, since $g$ acts trivially on $C$, $V$ admits a decomposition
    \[
    \tr_g V = \zeta_1 \cdot V_1 + \cdots + \zeta_r \cdot V_r,
    \]
    where $V_i \ (i=1,\cdots,r)$ are non-equivariant bundles in $K(C)$, and $\zeta_i \ (i=1,\cdots,r)$ are fiber-wise eigenvalues of $g$ which are roots of unity. In order to obtain an inverse of 
    \[
    \tr_g \Eu(a^{-1}V) = \prod_{i=1}^r \Eu(a^{-1} \cdot \zeta_i \cdot V_i),
    \]
    it suffices to look at the individual factors. Let $v_i$ be the rank of $V_i$, then
    \[
    \frac{1}{\Eu(a^{-1} \cdot \zeta_i \cdot V_i)} = \frac{1}{\sum_{j=0}^{v_i} (-a\zeta_i^{-1})^j \cdot \bigwedge^{j}V_i} = \frac{1}{(1-a\zeta_i)^{v_i} + \sum_{j=0}^{v_i} (-a\zeta_i^{-1})^j \cdot \left(\bigwedge^{j}V_i - \rank\bigwedge^{j}V_i\right)}
    \]
    is well-defined as long as $(1-a\zeta_i^{-1})$ is localized, as the second term in the denominator of RHS is nilpotent. In $\Lambda^{a,b}$, $(1-a\zeta)$ is inverted for any root of unity $\zeta$.
    \item Similarly, for any $\langle g \rangle$-equivariant vector bundle $W$ over $C$, $\tr_g \Eu(b^{-1}W)$ is well-defined and invertible in $\Lambda^{a,b}\otimes K(C)$.
\end{itemize}
Now, by construction
\[
R(a,b)_{0,m+1,d} = (-1)^{l\cdot\rank E_{0,m+1,d}} \cdot \frac{\Eu(a^{-1}(l\cdot R^0\ft_*\ev^*E)) \cdot \Eu(b^{-1}(l\cdot (R^1\ft_*\ev^*E)^\vee))}{\Eu(a^{-1}(l\cdot R^1\ft_*\ev^*E)) \cdot \Eu(b^{-1}(l\cdot (R^0\ft_*\ev^*E)^\vee))}.
\]
Its denominator is the product of $\Eu(a^{-1}V)$ and $\Eu(b^{-1}W)$ with
\[
V = l\cdot R^1\ft_*\ev^*E, \quad \text{and} \quad W = l\cdot (R^0\ft_*\ev^*E)^\vee
\]
being $\langle g \rangle$-equivariant vector bundles when restricted to $C$, and thus the trace of its restriction to $C$ is invertible in $\Lambda^{a,b}\otimes K(C)$. In this way, we obtain an interpretation of $\tr_g R(a,b)_{0,m+1,d}|_C$. On the other hand, it is hard to realize $R(a,b)_{0,m+1,d}$ directly as an element in $\Lambda^{a,b}\otimes K([\M_{0,m+1}(X,d)/S_m])$, 
\begin{definition} \label{DefinitionRab}
We define the $R(a,b)$-twisted big $\J$-function $\J^{X,R(a,b)} = \J^{X,R(a,b)}(\t)$ on $\K^{a,b}$ as above. We define $\L^{X,R(a,b)}\subset\K^{a,b}$ as its image cone.
\end{definition}

\par Taylor expansion in $a$ and $b$ induces an inclusion $\iota: \K^{a,b}\rightarrow\K[[a,b]]$ which shares similar properties with the one considered in Section \ref{PowerSeriesExpansion}. We denote by $\iota_+: \K^{a,b}_+\rightarrow\K_+[[a,b]]$ its restriction to the Lagrangian subspace of inputs. The big $\J$-functions on the two loop spaces may be identified under $\iota$ and $\iota_+$.
\begin{proposition}[Identification of big $\J$-functions] \label{TwistedBigJCommutes}
\[
\iota \left( \J^{X,R(a,b)}(\t) \right) = \J^{X,R(a,b)}(\iota_+(\t)), \quad \forall \t = \t(q) \in \Lambda^{a,b}_+ \cdot \K^{a,b}_+.
\]
In particular, $\L^{X,R(a,b)}\subset \K^{a,b}$ is exactly the intersection of (the preimage under $\iota$ of) $\L^{X,R(a,b)}\subset \K[[a,b]]$ with the $\Lambda^{a,b}_+$-neighborhood of $(1-q)\in \K^{a,b}$.
\end{proposition}
The proof of Proposition \ref{BigJCommutes} and Corollary \ref{LCommutes} may be directly transplanted here, with $\lambda$ replaced by $a$ and $b$. The only extra ingredient here is that the two big $\J$-functions are both twisted,
but it is not a problem as the twisting sheaves $R(a,b)_{0,m+1,d}$ are also identified under the Taylor expansion in $a$ and $b$.

\par Therefore, in order to prove that the expression $J^{X,R(a,b)}$ appearing in Proposition \ref{JfunctionRab}, now regarded as an honest rational function instead of its expansion in $a$ and $b$, lies on the smaller cone $\L^{X,R(a,b)}\subset \K^{a,b}$, it suffices to prove that it lies in the $\Lambda^{a,b}_+$-neighborhood at $(1-q)\in\K^{a,b}$. To prove that $J^{X,R(a,b)}$ is indeed in $\K^{a,b}$, we simply observe that all factors appearing in its denominator are either purely in $(a,q)$ or purely in $(b,q)$ (but never in $(a,b,q)$ in the same factor), and that the factors never have zeros when $a=b=0$ and $q = \zeta$ a root of unity. That $J^{X,R(a,b)}$ is $\Lambda^{a,b}_+$-close to $(1-q)$ is, in turn, due to our assumption in Proposition \ref{JfunctionRab} of $J$ being in (the $\Lambda_+$-germ at $(1-q)$ of) the classical loop space. Any such expression must satisfy $J_0 = J|_{Q_1=\cdots=Q_n=0} = 1-q$, so $J^{X,R(a,b)}_0 = J^X_0 = 1-q$ by construction in Proposition \ref{JfunctionRab}. All other terms in $J^{X,R(a,b)}$ contain non-trivial monomials in the Novikov variables, so $J^{X,R(a,b)}_0$ is indeed close to $(1-q)$ under the adic topology given by $\Lambda^{a,b}_+=(Q)$. In this way, we have proved 
\begin{proposition} \label{JfunctionRabRestricted}
For any point $J = \sum_d Q^d\cdot J_d\in\mathcal{L}^X\subset \K$ with $J_d\in K(X)(q)$,
\[
J^{X,R(a,b)} := \sum_d Q^d \cdot J_d \cdot \prod_{j=1}^{\rank E} (-1)^{-\langle c_1(f_j),d\rangle} \cdot \frac{\prod_{s=1}^{-\langle c_1(f_j),d\rangle} \left( 1 - a \cdot q^{s} \cdot f_j \right)^l}{\prod_{s=1}^{-\langle c_1(f_j),d\rangle} \left( 1 - b \cdot q^{-s} \cdot f_j^{-1} \right)^l} \quad\in \L^{X,R(a,b)}\subset \K^{a,b}.
\]
\end{proposition}
The difference from Proposition \ref{JfunctionRab} is that $J^{X,R(a,b)}$ here is regarded as a genuine rational function.

\subsection{Reduction of coefficients} \label{ReductionCoefficients}
\par The twisting of the type $R(a,b)$ studied in Section \ref{CompositionEulerTwisting} reduces to the $(E,l)$-level structure as $a=b=1$. We will treat the reduction of coefficients in two steps:
\begin{itemize}
    \item first, we take $a\mapsto\mu, b\mapsto\mu^{-1}$;
    \item second, we take $\mu \mapsto 1$. 
\end{itemize}
Vaguely speaking, the power series expansion considered in Section \ref{PowerSeriesExpansion} and \ref{CompositionEulerTwisting} may also be regarded as a reduction of coefficients between rational loop spaces, where we take the identity map $\lambda \mapsto \lambda$ on the parameters but alter the adic topology. In this section, we fix the adic topology (generated by the Novikov variables), but allow non-trivial (and not necessarily injective in particular) maps in between the parameters. 

\par We start by proving a general lemma regarding reduction of coefficients. Let $(\K^i,\Omega^i) = (\K^{\Lambda^i,S^i_+}, \Omega^i) \ (i=1,2)$ be two loop spaces as defined in Section \ref{variations}. In other words, denote
\begin{itemize}
    \item by $\Lambda^i = (T^i)^{-1}\C[\lambda^i]\otimes \C[[Q]]$ their respective coefficient ring (complete with respect to the prescribed ideal $\Lambda^i_+ = (Q)$), where $\lambda^i = (\lambda^i_1,\cdots,\lambda^i_{r_i})$ is a vector of parameters (of the polynomial ring), $T^i$ a multiplicatively closed subset, and $Q = (Q_1,\cdots,Q_n)$,
    \item by $\K^i = \K^i_+\oplus\K^i_-$ their respective Lagrangian polarization, where
    \[\K^i_{\pm} = K(X)\otimes \C[[Q]] \otimes (S^i_{\pm})^{-1} \C[\lambda^i,q],
    \]
    \item and by $\Omega^i$ their respective symplectic pairing.
\end{itemize}
Consider any $\C$-algebra homomorphism $\pi: \C[\lambda^1]\rightarrow\C[\lambda^2]$ satisfying
\begin{itemize}
    \item $\pi(T^1)\subset T^2$, and
    \item $\pi\otimes\id(S^1_+)\subset S^2_+$, for $\pi\otimes\id: \C[\lambda^1]\otimes\C[q]\rightarrow\C[\lambda^2]\otimes\C[q]$.
\end{itemize}
By extending trivially to the parts not involving $\lambda^1$ and $\lambda^2$, it gives rise to a $\C[[Q]]$-algebra homomorphism 
\[
\pi: \Lambda^1 = (T^1)^{-1}\C[\lambda^1]\otimes\C[Q] \rightarrow \Lambda^2 = (T^2)^{-1}\C[\lambda^2]\otimes\C[Q]
\]
which automatically preserves the Adams operations, and moreover a $\Lambda^1$-module homomorphism 
\[
\pi: \K^1 = K(X)\otimes \C[[Q]] \otimes (S^1_+\times S^1_-)^{-1}\C[q,\lambda^1]  \longrightarrow \K^2 = K(X)\otimes \C[[Q]] \otimes (S^2_+\times S^2_-)^{-1}\C[q,\lambda^2]
\]
where the $\Lambda^1$-module structure on the codomain is induced as usual from its natural $\Lambda^2$-module structure. Note that $S^i_-\ (i=1,2)$ are determined completely by $T^i\ (i=1,2)$, so our assumption that $\pi(T^1)\subset T^2$ implies already $\pi\otimes\id(S^1_-)\subset S^2_-$. It motivates our definition below.
\begin{definition} \label{DefReductionCoefficients}
A $\C$-algebra homomorphism $\pi: \C[\lambda^1] \rightarrow \C[\lambda^2]$ is called a reduction of coefficients from $\K^1$ to $\K^2$, if it satisfies
\begin{itemize}
    \item $\pi(T^1)\subset T^2$, 
    \item $\pi\otimes\id(S^1_+)\subset S^2_+$, for $\pi\otimes\id: \C[\lambda^1]\otimes\C[q]\rightarrow\C[\lambda^2]\otimes\C[q]$, and
    \item $\pi^*\Omega^2 = \Omega^1$, for the induced map $\pi: \K^1\rightarrow\K^2$
\end{itemize}
\end{definition}
The induced symplectic $\Lambda^1$-module homomorphism $\pi: \K^1\rightarrow \K^2$ has the following immediate properties:
\begin{itemize}
    \item $\pi(\K^1_{\pm}) \subset \K^2_{\pm}$, i.e. $\pi$ respects the Lagrangian polarization;
    \item $\pi(\Lambda^1_+) = \pi((Q)) \subset (Q) = \Lambda^2_+$, i.e. $\pi$ respects the adic topology.
\end{itemize}
Note that $\pi$ may not be injective. We denote by $\pi_{\pm}: \K^1_{\pm}\rightarrow \K^2_{\pm}$ the restriction of $\pi$ to the Lagrangian subspaces. The lemma below regarding $\pi_+$ resonates with Proposition \ref{LoopSpaceHom} regarding the Taylor expansion $\iota_+$.
\begin{lemma} \label{RedCoeffPropertiesLoopSpace}
\begin{itemize}
    \item[(a)] $\pi_+(1-q) = 1-q$.
    \item[(b)] $\pi_+$ commutes with the Adams operations on the coefficients. That is to say,
    \[
    \pi_+(\Psi^l(\t(q))) = \Psi^l(\pi_+(\t(q))), \quad \forall \t = \t(q)\in \K^1_+, l\in\mathbb{Z}_{>0}.
    \]
    \item[(c)] $\pi_+$ commutes with the changes of variables $q\mapsto \zeta q^l$ for any root of unity $\zeta$ and $l\in\mathbb{Z}_{>0}$. That is to say,
    \[
    \pi_+(\t(\zeta q^l)) = \pi_+(\t)(\zeta q^l), \quad \forall \t = \t(q)\in \K^\lambda_+.
    \]
\end{itemize}
\end{lemma}
As we have seen in the proof of Proposition \ref{BigJCommutes}, the big $\J$-functions are completely determined by the three ingredients above. Therefore, we have
\begin{lemma} \label{ReductionCoefficientsJfunc}
Let $\J^{X,i}(\t)\ (i=1,2)$ be the big $\J$-function of $X$ defined in the loop space $\K^i$, i.e. for inputs $\t\in \Lambda_+^i\cdot \K^i_+$. Then,
\[
\J^{X,2} \circ \pi_+ = \pi \circ \J^{X,1}.
\]
In particular, for the image cones $\L^{X,i}\subset \K^i$ of $\J^{X,i}(\t)$, we have $\pi(\L^{X,1})\subset \L^{X,2}$.
\end{lemma}
Note that the lemma is still correct when the big $\J$-functions $\J^{X,i}(\t)$ are twisted, as long as the twisting sheaves are well-defined both before and after the reduction of coefficients, and are identified under $\pi$. More precisely, for any Kawasaki stratum $C$ of $I[\M_{0,m+1}(X,d)/S_m]$, the homomorphism $\pi:\Lambda^1\rightarrow\Lambda^2$ between coefficient rings induces
\[
\pi_C = \pi\otimes\id_{K(C)}: \Lambda^1\otimes K(C) \rightarrow \Lambda^2\otimes K(C)
\]
and the following lemma holds regarding twisted big $\J$-functions.
\begin{lemma} \label{ReductionCoefficientsTwistedJfunc}
Let $\J^{X,R^i}(\t)\ (i=1,2)$ be, respectively, the big $\J$-function of $X$ over $\K^i$ with twisting $R^i$. If the reduction of coefficients $\pi$ further satisfies 
\[
\pi_C (\tr_g R^1_{0,m+1,d}|_C) = \tr_g R^2_{0,m+1,d}|_C
\]
for any Kawasaki stratum $C$ of $I[\M_{0,m+1}(X,d)/S_m]$ with prescribed local isotropy $g$. Then,
\[
\J^{X,R^2} \circ \pi_+ = \pi \circ \J^{X,R^1}.
\]
In particular, for the image cones $\L^{X,R^i}\subset \K^i$ of $\J^{X,R^i}(\t)$, we have $\pi(\L^{X,R^1})\subset \L^{X,R^2}$. Moreover, if $\pi$ is injective, $\L^{X,R^1}$ is exactly the intersection of $\pi^{-1}(\L^{X,R^2})$ with the $\Lambda^1_+$-neighborhood of $(1-q)\in \K^1$.
\end{lemma}
Lemma \ref{ReductionCoefficientsJfunc} and \ref{ReductionCoefficientsTwistedJfunc} are parallel to Proposition \ref{BigJCommutes} and \ref{TwistedBigJCommutes} on the level of big $\J$-functions. Note that a twisted big $\J$-function, defined with respect to the modification $R_{0,m+1,d}$ of virtual structure sheaf on $\M_{0,m+1}(X,d)$, is well-defined in a rational loop space with coefficient ring $\Lambda$ if and only if $\tr_g R_{0,m+1,d}|_C \in \Lambda\otimes K(C)$ over each Kawasaki stratum $C$. 

\par With the preparation above, we may now carry out the two-step reduction of coefficients described at the beginning of this section. For Step I,
\begin{itemize}
    \item we take $\K^1 = \K^{a,b}$ as in Section \ref{CompositionEulerTwisting}, for which in particular the coefficient ring is $\Lambda^1 = \Lambda^{a,b} = (T^1)^{-1}\C[a,b] \otimes \C[[Q]]$ with $T^1 = (\C[a]\setminus(a))\times (\C[b]\setminus(b))$, and endow it with the $R(a,b)$-twisted symplectic form $\Omega^1 = \Omega^{R(a,b)}$;
    \item we take $\K^2 = \K(\mu)$ (see the dictionary in Section \ref{variations}), for which in particular the coefficient ring is $\Lambda^2 = \C(\mu) \otimes \C[[Q]] = (T^2)^{-1}\C[\mu] \otimes \C[[Q]]$ with $T^2 = \C[\mu]\setminus 0$, and endow it with the symplectic form $\Omega^2 = \Omega^{(\mu^{-1}E,l)}$ with level structure $(\mu^{-1}E,l)$.
\end{itemize}
One may verify through direct computation that the $\C$-algebra homomorphism 
\begin{align*}
\pi:\quad & \C[a,b] \quad \longrightarrow \quad \C[\mu]\\
& \quad a,b \quad \ \longmapsto \quad \mu,\mu^{-1}
\end{align*}
satisfies all conditions of Definition \ref{DefReductionCoefficients} and is thus a legitimate reduction of coefficients. Moreover, it satisfies the extra assumption in Lemma \ref{ReductionCoefficientsTwistedJfunc} as well. Indeed, for any Kawasaki stratum $C$ of $I[\M_{0,m+1}(X,d)/S_m]$ with prescribed isotropy $g$, $\tr_g R(a,b)_{0,m+1,d}|_C$ reduces still to $ \tr_g \det^{-l}(\mu^{-1} E_{0,m+1,d})|_C$ as $a=\mu,b=\mu^{-1}$, because the process of ``expanding $R(a,b)_{0,m+1,d}$ in nilpotent part'' preceding Definition \ref{DefinitionRab} happens entirely in $\Rep(\langle g\rangle)\otimes K(C)$ and thus commutes naturally with any action on the coefficient ring, including in particular the specialization $a=\mu,b=\mu^{-1}$. Therefore, applying Lemma \ref{ReductionCoefficientsTwistedJfunc}, we have
\begin{proposition} \label{JfuncLevel1}
For any point $J = \sum_d Q^d\cdot J_d\in\mathcal{L}^X\subset \K$ with $J_d\in K(X)(q)$,
\[
J^{X,(\mu^{-1}E,l)} := \sum_d Q^d \cdot J_d \cdot \prod_{j=1}^{\rank E} \left[ (\mu f_j)^{-\langle c_1(f_j),d\rangle} q^{\frac{\langle c_1(f_j),d\rangle(\langle c_1(f_j),d\rangle - 1)}{2}} \right]^{l} \quad\in \L^{X,(\mu^{-1}E,l)}\subset \K(\mu).
\]
\end{proposition}
$J^{X,(\mu^{-1}E,l)}$ comes directly from specializing $a=\mu$ and $b=\mu^{-1}$ in the expression in Proposition \ref{JfunctionRabRestricted}. 

\par The result of Step I may be refined as follows. We note that $J^{X,(\mu^{-1}E,l)}$ as a point on the big $\J$-function $\J^{X,(\mu^{-1}E,l)}(\t)$ with $(\mu^{-1}E, l)$-level structure has ``Laurent polynomial'' input $\t = \t(q) \in \C[\mu,\mu^{-1}][q,q^{-1}]\otimes K(X)[[Q_1,\cdots,Q_n]]$. In fact, since we assumed $J$ to be a point in the classical loop space, it has $q$-Laurent polynomial input (even independent from $\mu$) and thus factors in the denominator of each $J_d$ can only be either $q$ or $(1-\zeta q)$ for root of unity $\zeta$. Moreover, the modification terms multiplied to $J_d$ in the expression of $J^{X,(\mu^{-1}E,l)}$ above involve only Laurent monomials in $q$ and $\mu$ as well, and thus introduce at most new factors being powers of $\mu$ into the denominator. In this way, $J^{X,(\mu^{-1}E,l)}_d$ lives naturally in the loop space (see again the dictionary in Section \ref{variations})
\[
\K[\mu,\mu^{-1}] := K(X)\otimes \C[[Q]] \otimes \left\{ \frac{f(q,\mu)}{q^{\gamma_1}\mu^{\gamma_2}\cdot \prod_s (1-\zeta_s q)}\ \middle| \ f\in \C[q,\mu], \gamma_1,\gamma_2\in\mathbb{Z}_{>0}, \zeta_s \text{ is root of unity}\right\}
\]
which is a subspace of $\K(\mu)$ and allows only ($\Lambda_+$-small) inputs in
\[
\K_+[\mu,\mu^{-1}] = K(X)\otimes \C[[Q]] \otimes \C[q,q^{-1},\mu,\mu^{-1}].
\]
$\K[\mu,\mu^{-1}]$ is a free module over the coefficient ring $\Lambda[\mu,\mu^{-1}] := \C[\mu]_\mu \otimes \C[[Q]] = \C[\mu,\mu^{-1}][[Q]]$ with adic topology given by the ideal $\Lambda_+ = (Q)$. It is not hard to see the projection of $J^{X,(\mu^{-1}E,l)}$ to $\K_+[\mu,\mu^{-1}]$ is indeed $\Lambda_+$-close to $(1-q)$.

\par The level structure $(\mu^{-1}E,l)$ is obviously still well-defined under these restricted settings. Denote still by $\J^{X,(\mu^{-1}E,l)}$ such narrower big $\J$-function with level structure $(\mu^{-1}E,l)$ and by $\L^{X,(\mu^{-1}E,l)} \subset \K[\mu,\mu^{-1}]$ its image. Then the above argument implies that $J^{X,(\mu^{-1}E,l)}$ is still a point on the narrower image cone.
\begin{proposition} \label{JfuncLevel2}
For any point $J = \sum_d Q^d\cdot J_d\in\mathcal{L}^X\subset \K$ with $J_d\in K(X)(q)$,
\[
J^{X,(\mu^{-1}E,l)} := \sum_d Q^d \cdot J_d \cdot \prod_{j=1}^{\rank E} \left[ (\mu f_j)^{-\langle c_1(f_j),d\rangle} q^{\frac{\langle c_1(f_j),d\rangle(\langle c_1(f_j),d\rangle - 1)}{2}} \right]^{l} \quad\in \L^{X,(\mu^{-1}E,l)}\subset \K[\mu,\mu^{-1}].
\]
\end{proposition}
The proposition above may also be regarded as direct implication of Lemma \ref{ReductionCoefficientsTwistedJfunc} for the reduction of coefficients $\pi: \K[\mu,\mu^{-1}]\rightarrow \K(\mu)$ induced from $\pi = \id: \C[\mu]\rightarrow \C[\mu]$ in the sense of Definition \ref{DefReductionCoefficients}.

\par Now we proceed to Step II. For this step, we simply take 
\begin{itemize}
    \item $\K^1 = \K[\mu,\mu^{-1}]$ at which we have arrived above, with the coefficient ring $\Lambda^1 = (T^1)^{-1}\C[\mu]\otimes\C[[Q]] = \C[\mu,\mu^{-1}]\otimes\C[[Q]]$ where $T^1 = \{c \cdot \mu^k \ | \ c\in\C^\times, k\in\mathbb{Z}_{\geq 0}\}$ and the symplectic form with $(\mu^{-1}E,l)$-level structure, and
    \item $\K^2 = \K$ the classical loop space defined in Section \ref{LoopSpacePFD}, with coefficient ring $\Lambda^2 = (T^2)^{-1}\C\otimes\C[[Q]] = \C[[Q]]$ where $T^2 = \C^\times$ but the symplectic form with $(E,l)$-level structure.
\end{itemize}
The $\C$-algebra homomorphism
\begin{align*}
\pi:\quad \C[\mu] & \longrightarrow \ \C \\
\mu \ & \longmapsto \ 1
\end{align*}
satisfies all conditions of Definition \ref{DefReductionCoefficients} and thus gives rise to a well-defined reduction of coefficients $\pi: \K^1\rightarrow\K^2$. Moreover, under the specialization $\mu=1$, the modification $\det^{-l}(\mu^{-1}E_{0,m+1,d})$ on the virtual structure sheaf reduces to $\det^{-l}(E_{0,m+1,d})$, and thus $\pi$ satisfies the extra condition in Lemma \ref{ReductionCoefficientsTwistedJfunc} by an argument similar to the one preceding Proposition \ref{JfuncLevel2}. Therefore, by Lemma \ref{ReductionCoefficientsTwistedJfunc}, the big $\J$-function with level structure $(\mu^{-1}E,l)$ in $\K^1$ reduces to the big $\J$-function with level structure $(E,l)$ in $\K^2 =\K$.
\begin{theorem} \label{JfuncLevel3}
For any point $J = \sum_d Q^d\cdot J_d\in\mathcal{L}^X\subset \K$ with $J_d\in K(X)(q)$,
\[
J^{X,(E,l)} := \sum_d Q^d \cdot J_d \cdot \prod_{j=1}^{\rank E} \left[ f_j^{-\langle c_1(f_j),d\rangle} q^{\frac{\langle c_1(f_j),d\rangle(\langle c_1(f_j),d\rangle - 1)}{2}} \right]^{l} \quad\in \L^{X,(E,l)}\subset \K.
\]
\end{theorem}
Here $f_j$ are the K-theoretic Chern roots of $E^\vee$. This recovers exactly the result in \cite{R-Z}.

\section{K-Theoretic Quantum Serre Duality} \label{QuantumSerreDuality}
\par In this section, we prove the genus-zero K-theoretic quantum Serre duality (KqSD).

\par We recall some definitions first. For simplicity of notation, we assume as usual
\[
E^\vee = \sum_{j=1}^{\rank E} f_j = \sum_{j=1}^{\rank E} \ f_j(P_1,\cdots,P_n) \in\K(X)
\]
for line bundles $f_j$. The case where negative terms do arise in the decomposition may be treated in exactly the same way, and the result will differ only by inverting the corresponding terms in the theorems below. 

\par Let $\J^{X,(\Eu,E,l)}$ be the $(\Eu,E,l)$-twisted big $\J$-function of $X$ defined with the modified virtual structure sheaves
\[
\mathcal{O}^{\virt,(\Eu,E,l)}_{g,m,d} = \mathcal{O}^{\virt}_{g,m,d} \otimes \Eu(\mu^{-1} E_{g,m,d}) \otimes \det^{-l}(\mu^{-1} E_{g,m,d})
\]
In other words, $\J^{X,(\Eu,E,l)}$ is the $(\Eu,\mu^{-1}E)$-twisted big $\J$-function with level structure $(\mu^{-1}E,l)$. Similarly, let $\J^{X,(\Eu^{-1},E^\vee,l+1)}$ be the $(\Eu^{-1},E^\vee,l+1)$-twisted big $\J$-function of $X$ defined with the modified virtual structure sheaves
\[
\mathcal{O}^{\virt,(\Eu^{-1},E^\vee,l+1)}_{g,m,d} = \mathcal{O}^{\virt}_{g,m,d} \otimes \Eu^{-1}(\mu (E^\vee)_{g,m,d}) \otimes \det^{-(l+1)}(\mu (E^\vee)_{g,m,d})
\]
$\J^{X,(\Eu^{-1},E^\vee,l+1)}$ is the $(\Eu^{-1},\mu E^\vee)$-twisted big $\J$-function with level structure $(\mu E^\vee,l+1)$. We denote by $\L^{X,(\Eu,E,l)}$ and $\L^{X,(\Eu^{-1},E^\vee,l+1)}$ the two image cones. Both cones live naturally in the loop space $\K(\mu)$ (see Section \ref{variations}) where arbitrary rational functions in $q$ and $\mu$ are allowed, but the symplectic structures on $\K(\mu)$ associated to the two theories are different.

\par We make two remarks regarding the definitions. First, the twisting of the type $(\Eu^{-1},\mu E^\vee)$ is indeed invertible, in the sense that its effect on the big $\J$-function is governed by the qARR formula after being expanded into the exponential form. However, it is different from the $(\Eu,\mu^{-1} E)^\vee$-twisting of Section \ref{EulerTwisting}. Second, we will not able to eventually remove the fiber-wise $\C^\times$-action with equivariant parameter $\mu$ like we did in Section \ref{LevelStructuresRevisited}, as otherwise the two Euler-type twistings arising in the KqSD will no longer be invertible. 

\subsection{Statement of the theorem} \label{StatementQuantumSerre}
\par We state (various versions of) the duality as follows. Note that the $(\Eu,E,l)$- and $(\Eu^{-1},E^\vee,l+1)$-twisted theories are defined on the same loop space but with different symplectic forms.
\begin{theorem}[K-Theoretic Quantum Serre Duality] \label{TheoremQuantumSerre}
For any point 
\[
J^1 = \sum_d \prod_{i=1}^n Q_i^{d_i} \cdot J^1_d \quad \in \L^{X,(\Eu,E,l)} \subset \K(\mu)
\]
with $J^1_d\in K(X)(q,\mu)$, we have $J^2 \in \L^{X,(\Eu^{-1},E^\vee,l+1)} \subset \K(\mu)$, where
\begin{align}
J^2 & = \sum_d \ \prod_{i=1}^n \left( Q_i \cdot \left(-q^{l}\right)^{c_1(E)_i}\right)^{d_i} \cdot J^1_d \cdot \prod_{j=1}^{\rank E}\frac{1-\mu\cdot f_j}{1-\mu \cdot f_j \cdot q^{-\langle c_1(f_j), d\rangle}} \label{J2mu}\\
& = \sum_d \ \prod_{i=1}^n \left(Q_i \cdot \left(-q^{(l+1)}\right)^{c_1(E)_i}\right)^{d_i} \cdot J^1_d \cdot \prod_{j=1}^{\rank E}\frac{1-\mu^{-1} \cdot f_j^{-1}}{1-\mu^{-1} \cdot f_j^{-1} \cdot q^{\langle c_1(f_j), d\rangle}}. \label{J2muinverse}
\end{align}
Here $c_1(E)_i$ refers to the $i$-th component of $c_1(E)$ under the basis $c_1(P_1),\cdots c_1(P_n)$ of $H_2(X)$. 
\end{theorem}
The expression of $J^2$ in the above theorem suggests that the non-equivariant limit at $\mu = 1$ should exist. Yet, we will not discuss implications of such limit in the present paper.

\par The result may be refined by looking into the powers of $\mu$. Indeed, in the expression of $J^2$ in Theorem \ref{TheoremQuantumSerre}, the factors multiplied to $J^1_d$ do not have poles at $\mu=0$ or $\mu = \infty$, so presumably the theorem is true when restricted to smaller rational loop spaces like $\K^\mu\subset \K(\mu)$ or $\K^{\mu^{-1}}\subset \K(\mu)$, where poles at at $\mu=0$ or $\mu = \infty$ are forbidden respectively. (Here $\K^\mu$ and $\K^{\mu^{-1}}$ are defined as replacing $\lambda$ respectively by $\mu$ and $\mu^{-1}$ in the loop space $\K^\lambda$ of Section \ref{CentralDefinition}.) However, the truth is that neither poles at $\mu=0$ nor at $\mu=\infty$ can be removed completely, because such poles, coming from the twisted correlators, are caused not entirely by the inputs $\t(q)\in\K_+(\mu)$, but rather more intrinsically from the modified virtual structure sheaves. For instance, in the $(\Eu,E,l)$-twisted theory, the modification on $\mathcal{O}^{\virt}_{0,m+1,d}$ takes the form
\[
\Eu(\mu^{-1}E_{0,m+1,d}) \otimes \det^{-l}(\mu^{-1}E_{0,m+1,d}).
\]
While the first term produces no poles at $\mu=0$ or $\mu = \infty$ to the correlator $\langle\cdot\rangle_{0,m+1,d}$, the second term gives
\[
\mu^{l\cdot\rank E_{0,m+1,d}} = \mu^{l \cdot \langle c_1(E),d \rangle + l\cdot \rank E}.
\]
The change in Poincaré pairing endows the entire dual basis $\{\phi^\alpha\}$ with an extra factor of $\Eu^{-1}(\mu^{-1}E) \otimes \det^{l}(\mu^{-1} E)$, and in particular the constant power $\mu^{l\cdot \rank E}$. In total, in a term of degree-$d$ in $\J^{X,(\Eu,E,l)}$ (i.e. marked with $Q^d$), the power of $\mu$ that arises due to the twisting is exactly $\mu^{l \cdot \langle c_1(E),d \rangle}$, which may be absorbed into the Novikov variables by 
\[
Q_i\mapsto Q_i' = Q_i \cdot \mu^{l\cdot c_1(E)_i}.
\]
The exact same powers of $\mu$ arise in the $(\Eu^{-1},E^\vee,l+1)$-twisted theory as well. In other words, while the untwisted big $\J$-function exists in a $\Lambda_+$-neighborhood of $(1-q)$ in $\K^\mu$, where $\Lambda_+ = (Q) \subset \Lambda = \C[\mu]_{(\mu)}\otimes\C[[Q]]$, the $(\Eu,E,l)$- and $(\Eu^{-1},E^\vee,l+1)$-twisted big $\J$-functions more naturally exist in a different subspace of $\K(\mu)$ completed under a different adic topology. More precisely, let $\K^{\mu,Q'}$ be the rational loop space resulting from the replacement $Q_i\mapsto Q_i' \ (i=1,\cdots,n)$. It is a free module over the coefficient ring $\Lambda^{\mu,Q'} = \C[\mu]_{(\mu)}\otimes\C[[Q']] = \C[\mu]_{(\mu)}\otimes\C[[Q'\cdot \mu^{l\cdot c_1(E)}]]$, complete under the adic topology of $\Lambda^{\mu,Q'}_+ = (Q')\subset\Lambda^{\mu,Q'}$, where $Q' = (Q_1',\cdots,Q_n')$. $\J^{X,(\Eu,E,l)}(\t)$ and $\J^{X,(\Eu^{-1},E^\vee,l+1)}(\t)$ are both defined for $\t\in\Lambda^{\mu,Q'}_+\cdot\K_+^{\mu,Q'}$, and we denote by $\L^{X,(\Eu,E,l)}\subset\K^{\mu,Q'}$ and $\L^{X,(\Eu^{-1},E^\vee,l+1)}\subset\K^{\mu,Q'}$ their image cones.
\begin{corollary} \label{Cor2}
For any 
\[
J^1 = \sum_d \prod_{i=1}^n (Q_i')^{d_i} \cdot J^1_d = \sum_d \prod_{i=1}^n Q_i^{d_i} \cdot J^1_d \cdot \mu^{l\langle c_1(E),d \rangle} \quad \in \L^{X,(\Eu,E,l)} \subset \K^{\mu,Q'} \subset \K(\mu),
\]
$J^2$ as is given by in Formula (\ref{J2mu}) but with each $Q_i$ replaced by $Q_i'$ resides in $\L^{X,(\Eu^{-1},E^\vee,l+1)}\subset\K^{\mu,Q'}$.
\end{corollary}

\par Similarly, if we extract a factor of $\mu^{(l+1) \cdot \langle c_1(E),d \rangle}$ from the correlators, the remaining expression will be more naturally written in terms of $\mu^{-1}$, still with poles neither at $\mu^{-1} = 0$ nor at $\mu^{-1} = \infty$. We take 
\[
Q_i \mapsto Q_i'' = Q_i \mu^{(l+1)\cdot c_1(E)_i},
\]
and consider the rational loop space $\K^{\mu^{-1},Q''}$ defined by replacing $Q_i\mapsto Q_i'' \ (i=1,\cdots,n)$ in $\K^{\mu^{-1}}$. $\K^{\mu^{-1},Q''}$ is a free module over the coefficient ring $\Lambda^{\mu^{-1},Q''} = \C[\mu^{-1}]_{(\mu^{-1})}[[Q'']]$, complete under the adic topology defined by $\Lambda^{\mu^{-1},Q''}_+ = (Q'')\subset\Lambda^{\mu^{-1},Q''}$, where $Q'' = (Q_1'',\cdots,Q_n'')$. $\J^{X,(\Eu,E,l)}(\t)$ and $\J^{X,(\Eu^{-1},E^\vee,l+1)}(\t)$ are both defined for $\t\in\Lambda^{\mu^{-1},Q''}_+ \cdot \K_+^{\mu^{-1},Q''}$, and we denote by $\L^{X,(\Eu,E,l)}\subset\K^{\mu^{-1},Q''}$ and $\L^{X,(\Eu^{-1},E^\vee,l+1)}\subset\K^{\mu^{-1},Q''}$ their image cones.
\begin{corollary} \label{Cor3}
For any 
\[
J^1 = \sum_d \prod_{i=1}^n (Q_i'')^{d_i} \cdot J^1_d = \sum_d \prod_{i=1}^n Q_i^{d_i} \cdot J^1_d \cdot \mu^{(l+1)\langle c_1(E),d \rangle} \quad \in \L^{X,(\Eu,E,l)} \subset \K^{\mu^{-1},Q''} \subset \K(\mu),
\]
$J^2$ as is given by in Formula (\ref{J2muinverse}) but with each $Q_i$ replaced by $Q_i''$ resides in $\L^{X,(\Eu^{-1},E^\vee,l+1)}\subset\K^{\mu^{-1},Q''}$.
\end{corollary}

\par In another direction, if we allow power series expansion in either $\mu$ or $\mu^{-1}$, the KqSD admits the following simpler forms like before. We denote by $\K[[\mu]]^{Q'}$ (resp. $\K[[\mu^{-1}]]^{Q''}$) the loop space as defined in Section \ref{RationalInputs}, but with $\lambda$ replaced by $\mu$ (resp. $\mu^{-1}$) and with $Q_i$ replaced by $Q_i'$ (resp. $Q_i''$) for $1\leq i\leq n$. Note that $\K[[\mu]]^{Q'}$ (resp. $\K[[\mu^{-1}]]^{Q''}$) is over the coefficient ring $\Lambda[[\mu]]^{Q'} = \C[[\mu,Q']] = \C[[\mu,Q_1',\cdots,Q_n']]$ (resp. $\Lambda[[\mu^{-1}]]^{Q''} = \C[[\mu^{-1},Q'']] = \C[[\mu^{-1},Q_1'',\cdots,Q_n'']]$) and is equipped with the adic topology defined by the ideal $\Lambda[[\mu]]^{Q'}_+ = (\mu,Q')$ (resp. $\Lambda[[\mu^{-1}]]^{Q''}_+ = (\mu^{-1},Q'')$). Both $\J^{X,(\Eu,E,l)}(\t)$ and $\J^{X,(\Eu^{-1},E^\vee,l+1)}(\t)$ are defined on both of the spaces, in a small neighborhood of $(1-q)$. Hence, Theorem \ref{TheoremQuantumSerre} admits two variations as below, paraphrasing respectively Theorem \ref{MainA} and \ref{MainB} in Introduction in a more rigorous way.
\begin{manualtheorem}{2A}[K-Theoretic Quantum Serre Duality] \label{1A}
For any point 
\[
J^1 = \sum_d \prod_{i=1}^n (Q_i')^{d_i} \cdot J^1_d \quad \in \L^{X,(\Eu,E,l)} \subset \K[[\mu]]^{Q'}
\]
with $J^1_d\in K(X)[[\mu]](q)$, we have
\[
J^2 = \Eu(\mu^{-1} E) \cdot \sum_d \ \prod_{i=1}^n \left( Q_i' \cdot \left(-q^{l}\right)^{c_1(E)_i}\right)^{d_i} \cdot J^1_d
\quad \in \L^{X,(\Eu^{-1},E^\vee,l+1)} \subset \K[[\mu]]^{Q'}.
\]
\end{manualtheorem}
\begin{manualtheorem}{2B}[K-Theoretic Quantum Serre Duality] \label{1B}
For any point 
\[
J^1 = \sum_d \prod_{i=1}^n (Q_i'')^{d_i} \cdot J^1_d \quad \in \L^{X,(\Eu,E,l)} \subset \K[[\mu^{-1}]]^{Q''}
\]
with $J^1_d\in K(X)[[\mu^{-1}]](q)$, we have
\[
J^2 = \Eu(\mu E^\vee) \cdot \sum_d \ \prod_{i=1}^n \left( Q_i'' \cdot \left(-q^{(l+1)}\right)^{c_1(E)_i}\right)^{d_i} \cdot J^1_d
\quad \in \L^{X,(\Eu^{-1},E^\vee,l+1)} \subset \K[[\mu^{-1}]]^{Q''}.
\]
\end{manualtheorem}
In other words, $\L^{X,(\Eu,E,l)}$ and $\L^{X,(\Eu^{-1},E^\vee,l+1)}$ differ by a constant scaling followed by a modification of the Novikov variables by signed powers of $q$. Note that $\Eu(\mu^{-1} E) \equiv 1 \ (\mod \ \mu)$ in $K(X)[[\mu]]$ and $\Eu(\mu E^\vee) \equiv 1 \ (\mod \ \mu^{-1})$ in $K(X)[[\mu^{-1}]]$. In this sense, the scalings are close to identity under the correponding adic topology.

\subsection{Proof of the theorem} \label{ProofQuantumSerre}
\par This section is mainly devoted to proving Theorem \ref{TheoremQuantumSerre} where we work with the loop space $\K(\mu)$. During the proof, we make small digressions to explain how it may be trimmed to justify Corollary \ref{Cor2} and \ref{Cor3}. Essentially, the problem is only due to the powers of $\mu$ created by the level structures, and may be settled by taking the change of variables from $Q_i$ to $Q_i'$ or $Q_i''$ as described earlier. In the end, we explain how Theorem \ref{1A} and \ref{1B} may be proved along the way by slightly altering the PFD operators that we employ.

\par Following the same idea of Section \ref{LevelStructuresRevisited} to realize the level structures in terms of invertible twisting. The virtual structure sheaf on $\M_{0,m+1}(X,d)$ modified by the $(\Eu,E,l)$-twisting, composed from the twisting of type $(\Eu,\mu^{-1} E)$ and the $(\mu^{-1} E,l)$-level structure, is given by $\mathcal{O}^{\virt,(\Eu,E,l)}_{g,m,d} = $
\[
\mathcal{O}^{\virt}_{g,m,d} \otimes (-1)^{l\cdot \rank E_{g,m,d}} \cdot \exp{\left((l+1)\sum_{k<0} \frac{a^{-k}\Psi^k(E_{g,m,d})}{k} + l \sum_{k>0} \frac{b^k\Psi^k(E_{g,m,d})}{k} \middle)\right|_{a=\mu, b=\mu^{-1}}}.
\]
Similarly, the virtual structure sheaf modified by the $(\Eu^{-1},E^\vee,l+1)$-twisting, composed from the twisting of type $(\Eu^{-1},\mu E^\vee)$ and the $(\mu E^\vee,l+1)$-level structure, is given by $\mathcal{O}^{\virt,(\Eu^{-1},E^\vee,l+1)}_{g,m,d} = $
\[
\mathcal{O}^{\virt}_{g,m,d} \otimes (-1)^{(l+1)\cdot \rank (E^\vee)_{g,m,d}} \cdot \exp{\left(l\sum_{k<0} \frac{b^{-k}\Psi^k((E^\vee)_{g,m,d})}{k} + (l+1)\sum_{k>0} \frac{a^k\Psi^k((E^\vee)_{g,m,d})}{k} \middle)\right|_{a=\mu, b=\mu^{-1}}}.
\]
The discrepancy between the two signs is
\[
\mathfrak{s} = (-1)^{l\cdot \rank E_{g,m,d}} \cdot (-1)^{(l+1)\cdot \rank (E^\vee)_{g,m,d}} = (-1)^{\langle c_1(E),d\rangle + \rank E}.
\]
We postpone the consideration of $\mathfrak{s}$ until later by removing the signs $(-1)^{l\cdot \rank E_{g,m,d}}$ and $(-1)^{(l+1)\cdot \rank (E^\vee)_{g,m,d}}$ and focusing only on the two sign-free exponential-type twistings for now. We denote by $\L^{X,1}$ and $\L^{X,2}$ the two image cones defined via the modified virtual structure sheaves before specialization to $a = \mu$ and $b=\mu^{-1}$, and by $\L^{X,|\Eu,E,l|}$ and $\L^{X,|\Eu^{-1},E^\vee,l+1|}$ the ones defined via the virtual structure sheaves after such specialization. Applying the qARR formula \cite{Givental:perm11}, we obtain
\[
\L^{X,1} = \Box_1^{-1} \L^X \quad \text{and} \quad \L^{X,2} = \Box_2^{-1} \L^X,
\]
where
\begin{align*}
\Box_1 & = \exp{\left((l+1)\sum_{k>0} \frac{a^{k}\Psi^{-k}(E)\cdot q^k}{k(1-q^k)} + l \sum_{k>0} \frac{b^k\Psi^k(E)}{k(1-q^k)} \right)},\\
\Box_2 & = \exp{\left(l\sum_{k>0} \frac{b^{k}\Psi^{k}(E)\cdot q^k}{k(1-q^k)} + (l+1) \sum_{k>0} \frac{a^k\Psi^{-k}(E)}{k(1-q^k)} \right)}.
\end{align*}
The qARR formula is still correct for the generalized loop space $\K(\mu)$ because our way of defining the big $\J$-functions on $\K(\mu)$ is compatible with the adelic characterization, and thus the proof in \cite{Givental:perm11} of the qARR formula. At the current stage, all three cones appearing above are regarded as living in the loop space $\K(\mu)[[a,b]]$, where power series in $a$ and $b$ are allowed. 

\par Modifying $\Box_1$ and $\Box_2$ respectively by the PFD operators 
\begin{align*}
F_1 & = \exp{\left( - (l+1) \sum_{k>0}\sum_{j=1}^{\rank E} \frac{a^{k}q^k f_j(P^kq^{kQ\partial_Q})}{k(1-q^k)} - l \sum_{k>0}\sum_{j=1}^{\rank E} \frac{b^k f_j(P^{-k}q^{-kQ\partial_Q})}{k(1-q^k)} \right)},\\
F_2 & = \exp{\left( - l \sum_{k>0}\sum_{j=1}^{\rank E} \frac{b^{k}q^k f_j(P^{-k}q^{-kQ\partial_Q})}{k(1-q^k)} - (l+1) \sum_{k>0}\sum_{j=1}^{\rank E} \frac{a^k f_j(P^kq^{kQ\partial_Q})}{k(1-q^k)} \right)},
\end{align*}
where $P^kq^{kQ\partial_Q} = (P_1^kq^{k Q_1\partial_{Q_1}},\cdots,P_n^kq^{k Q_n\partial_{Q_n}})$ and $P^{-k}q^{-kQ\partial_Q} = (P_1^{-k}q^{-kQ_1\partial_{Q_1}},\cdots,P_n^{-k}q^{-k Q_n\partial_{Q_n}})$,
we obtain new operators $D_1 := F_1\circ\Box_1$ and $D_2 := F_2\circ\Box_2$ whose effects on the Novikov variables look nicer
\begin{align*}
D_1\cdot\prod_{i=1}^n Q_i^{d_i} & = \prod_{j=1}^{\rank E}\prod_{s=1}^{-\langle c_1(f_j),d\rangle} \frac{(1-b\cdot f_{j}^{-1}\cdot q^{-s})^l}{(1-a\cdot f_j\cdot q^s)^{l+1}} \cdot \prod_{i=1}^n Q_i^{d_i},\\
D_2\cdot\prod_{i=1}^n Q_i^{d_i} & = \prod_{j=1}^{\rank E}\prod_{s=0}^{-\langle c_1(f_j),d\rangle - 1} \frac{(1-b\cdot f_{j}^{-1}\cdot q^{-s})^l}{(1-a\cdot f_j\cdot q^s)^{l+1}} \cdot \prod_{i=1}^n Q_i^{d_i}.
\end{align*}
Here the computation largely resembles that of Section \ref{EulerTwisting}, so we omit the detail. By Lemma \ref{PFD}, $D_1^{-1}$ and $D_2^{-1}$ still send $\L^X$ to $\L^{X,1}$ and to $\L^{X,2}$ respectively.

\par Now we are ready to restrict to a smaller loop space which we denote by $\K(\mu)^{a,b}$ (see the dictionary of Section \ref{variations} for definition). The embedding $\iota: \K(\mu)^{a,b}\rightarrow \K(\mu)[[a,b]]$ enjoys the same properties as its prototypes $\K^{\lambda}\rightarrow \K[[\lambda]]$ in Section \ref{PowerSeriesExpansion} and $\K^{a,b}\rightarrow \K[[a,b]]$ in Section \ref{CompositionEulerTwisting}. In particular, if $J^0 = \sum_d Q^d \cdot J^0_d \in \L^X\subset\K(\mu)$, we have
\begin{align*}
D_1^{-1} \cdot J^0 & = \sum_d Q^d \cdot J^0_d \cdot \prod_{j=1}^{\rank E}\prod_{s=1}^{-\langle c_1(f_j),d\rangle} \frac{(1-b\cdot f_{j}^{-1}\cdot q^{-s})^{-l}}{(1-a\cdot f_j\cdot q^s)^{-(l+1)}} \quad \in \L^{X,1} \subset \K(\mu)^{a,b},\\
D_2^{-1} \cdot J^0 & = \sum_d Q^d \cdot J^0_d \cdot \prod_{j=1}^{\rank E}\prod_{s=0}^{-\langle c_1(f_j),d\rangle - 1} \frac{(1-b\cdot f_{j}^{-1}\cdot q^{-s})^{-l}}{(1-a\cdot f_j\cdot q^s)^{-(l+1)}} \quad \in \L^{X,2} \subset \K(\mu)^{a,b}.
\end{align*}
\par By construction, the specialization $a = \mu$ and $b=\mu^{-1}$ reduces the modified virtual structure sheaves defining $\L^{X,1}$ and $\L^{X,2}$ respectively to those defining $\L^{X,|\Eu,E,l|}$ and $\L^{X,|\Eu^{-1},E^\vee,l+1|}$. Therefore, by Lemma \ref{ReductionCoefficientsTwistedJfunc}, the reduction of coefficients induced by
\begin{align*}
\pi: \quad \C[\mu,a,b] & \longrightarrow \quad \C[\mu] \\
 \mu, a, b \ & \longmapsto \ \mu, \mu, \mu^{-1}
\end{align*}
(easily checked to satisfy the premises of Definition \ref{DefReductionCoefficients}) sends $\L^{X,1}$ to $\L^{X,|\Eu,E,l|}$ and $\L^{X,2}$ to $\L^{X,|\Eu^{-1},E^\vee,l+1|}$. As a result,
\begin{align*}
I^1 & := (\overline{D_1})^{-1} \cdot J^0 = \sum_d Q^d \cdot J^0_d \cdot \prod_{j=1}^{\rank E}\prod_{s=1}^{-\langle c_1(f_j),d\rangle} (-\mu f_j q^s)^l \cdot (1-\mu f_j q^s)   \quad \in \L^{X,|\Eu,E,l|}\subset \K(\mu),\\
I^2 & := (\overline{D_2})^{-1} \cdot J^0 = \sum_d Q^d \cdot J^0_d \cdot \prod_{j=1}^{\rank E}\prod_{s=0}^{-\langle c_1(f_j),d\rangle - 1} (-\mu f_j q^s)^l \cdot (1-\mu f_j q^s)    \quad \in \L^{X,|\Eu^{-1},E^\vee,l+1|}\subset \K(\mu).
\end{align*}
Here $\overline{D_1}$ and $\overline{D_2}$ are the specialization of $D_1$ and $D_2$ under $a = \mu, b =\mu^{-1}$. It is not hard to see that if $J^0 \equiv I^1 \equiv I^2 \ (\mod\ Q)$. Hence, if $J^0$ is $\Lambda_+(\mu)$-close to $1-q$, so are $I^1$ and $I^2$, where $\Lambda_+(\mu)$ is the ideal generated by the Novikov variables in the coefficient ring $\Lambda(\mu)$ of $\K(\mu)$.

\par Here we make a digression toward Corollary \ref{Cor2} and \ref{Cor3}. In both $I^1$ and $I^2$, for a given degree $d$, the factor multiplied to $Q^d\cdot J_0^d$ may be split into a power of $\mu$ and a part without poles at $\mu=0$ or $\mu=\infty$, and the power of $\mu$ is, in both cases,
\[
\prod_{j=1}^{\rank E} \mu^{ - l\cdot \langle c_1(f_j),d\rangle} = \mu^{l\cdot \langle c_1(E),d\rangle} = \prod_{i=1}^n \mu^{c_1(E)_i\cdot d_i}.
\]
Whether it is a positive or negative power of $\mu$ that appears depends on the choice of $d$. It may well happen that the exponent is positive in a proper subcone of the Kahler cone, and negative in the complement of closure. For this reason, the argument above cannot be restricted directly to subspaces like $\K^\mu$ and $\K^{\mu^{-1}}$ (for Corollary \ref{Cor2} and \ref{Cor3}), or more generally $\K[[\mu]]$ and $\K[[\mu^{-1}]]$ (for Theorem \ref{1A} and \ref{1B}), where either positive or negative powers of $\mu$ are abandoned. We solve this problem by modifying the Novikov variables and considering the shifted loop spaces as described earlier, as the exponents of $\mu$ above are linear in $d$.

\par For the case of $\mu$, we take $Q_i' = Q_i \cdot \mu^{l\cdot c_1(E)_i}$ so that
\[
\prod_{i=1}^d (Q_i')^{d_i} = \mu^{l\cdot \langle c_1(E),d\rangle} \prod_{i=1}^d Q_i^{d_i}.
\]
Consequently, when written as a power series in $Q_1',\cdots,Q_n'$,
\begin{align*}
I^1 & = \sum_d \left(\prod_{i=1}^{n} (Q_i')^{d_i}\right) \cdot J^0_d \cdot \prod_{j=1}^{\rank E}\prod_{s=1}^{-\langle c_1(f_j),d\rangle} (- f_j q^s)^l \cdot (1-\mu f_j q^s),\\
I^2 & = \sum_d \left(\prod_{i=1}^{n} (Q_i')^{d_i}\right) \cdot J^0_d \cdot \prod_{j=1}^{\rank E}\prod_{s=0}^{-\langle c_1(f_j),d\rangle - 1} (- f_j q^s)^l \cdot (1-\mu f_j q^s),
\end{align*}
so they are both in $\K^{\mu,Q'}$ and $\Lambda^{\mu,Q'}_+$-close to $(1-q)$ (or more generally, $\Lambda[[\mu]]^{Q'}_+$-close to $(1-q)$ in $\K[[\mu]]^{Q'}$).

\par Similarly, for the case of $\mu^{-1}$, we take $Q_i'' := Q_i \cdot \mu^{(l+1)\cdot c_1(E)_i}$ so that 
\[
\prod_{i=1}^d (Q_i'')^{d_i} = \mu^{(l+1)\cdot \langle c_1(E),d\rangle} \prod_{i=1}^d Q_i^{d_i},
\]
and the remaining factors are purely in $\mu^{-1}$ with no poles at $\mu^{-1} = 0$. Consequently,
\begin{align*}
I^1 & = \sum_d \left(\prod_{i=1}^{n} (Q_i'')^{d_i}\right) \cdot J^0_d \cdot \prod_{j=1}^{\rank E}\prod_{s=1}^{-\langle c_1(f_j),d\rangle} (- f_j q^s)^l \cdot (\mu^{-1} - f_j q^s),\\
I^2 & = \sum_d \left(\prod_{i=1}^{n} (Q_i'')^{d_i}\right) \cdot J^0_d \cdot \prod_{j=1}^{\rank E}\prod_{s=0}^{-\langle c_1(f_j),d\rangle - 1} (- f_j q^s)^l \cdot (\mu^{-1} - f_j q^s),
\end{align*}
so they are both $\Lambda^{\mu^{-1},Q''}_+$-close to $(1-q)$ in $\K^{\mu^{-1},Q''}$ (or $\Lambda[[\mu^{-1}]]^{Q''}_+$-close to $(1-q)$ in $\K[[\mu^{-1}]]^{Q''}$).

\par Now we go back the main proof. Since $\overline{D_1}$ and $\overline{D_2}$ are both invertible on $\K(\mu)$, the point-wise maps $J^0\mapsto I^1$ and $J^0\mapsto I^2$ are both bijective, regarded respectively as from $\L^X$ to $\L^{X,|\Eu,E,l|}$ and $\L^{X,|\Eu^{-1},E^\vee,l+1|}$. Therefore, $(\overline{D_2})^{-1}\circ \overline{D_1}$, combined with the sign discrepancy from level structures of the two twistings that we have dropped, gives then a bijective map from $\L^{X,(\Eu,E,l)}$ to $\L^{X,(\Eu^{-1},E^\vee,l+1)}$. 

\par The sign discrepancy has two sources: one from the modification on the virtual structure sheaves which gives $\mathfrak{s}$ as we have computed at the beginning, and the other from the modification on the Poincaré pairings which gives $(-1)^{l\rank E - (l+1)\rank E^\vee} = (-1)^{\rank E}$. Merging the two sources, we obtain the total sign change on the degree-$d$ correlators
\[
\mathfrak{s} \cdot (-1)^{\rank E} = (-1)^{\langle c_1(E),d\rangle}.
\]
Therefore, for $J^1 = \sum_d Q^d\cdot J^1_d \in \L^{X,(\Eu,E,l)}\subset \K(\mu)$, we have
\begin{align*}
J^2 & = (\overline{D_2})^{-1}\circ \overline{D_1} \cdot \left(\sum_d Q^d\cdot (-1)^{\langle c_1(E),d\rangle} \cdot J^1_d \right) \\
& = \sum_d \prod_{i=1}^n \left((-1)^{c_1(E)_i}Q_i\right)^{d_i} \prod_{k=1}^{\rank E} \frac{(1-\mu^{-1}f_j^{-1}q^{\langle c_1(f_j),d\rangle})^l (1-\mu f_j)^{l+1}}{(1-\mu^{-1}f_j^{-1})^l (1-\mu f_j q^{-\langle c_1(f_j),d\rangle})^{l+1}} \cdot J^1_d \ \in \L^{X,(\Eu^{-1},E^\vee,l+1)}\subset \K(\mu).
\end{align*}
By direct computation, it reduces to Formulae (\ref{J2mu}) and (\ref{J2muinverse}) (equivalent in $\K(\mu)$) of Theorem \ref{TheoremQuantumSerre}. This completes our proof of Theorem \ref{TheoremQuantumSerre}. Had we composed directly $\Box_2^{-1}$ and $\Box_1$ without modification by PFD operators, we would have eventually arrived at $J^2 = J^1 \cdot \Eu^{l+1}(\mu^{-1}E) / \Eu^{l}(\mu E^{-1})$, which would have looked simpler. Nevertheless, it is impossible to obtain in such case a map from the germ at $(1-q)$ to itself.

\par As for Corollary \ref{Cor2}, we note that when restricted to $\K^{\mu}$, the invertible map $\overline{D_1}$ gives a one-to-one correspondence between $\L^X\subset\K^{\mu}$ and $\L^{X,|\Eu,E,l|}\subset\K^{\mu,Q'}$ by sending $J^0$ to $I^1$, and $\overline{D_2}$ between $\L^X\subset\K^{\mu}$ and $\L^{X,|\Eu^{-1},E^\vee,l+1|}\subset\K^{\mu,Q'}$ by sending $J^0$ to $I^2$. Therefore, incorporating the signs as before, we obtain a correspondence between points on $\L^{X,(\Eu,E,l)} \subset \K^{\mu,Q'}$ and $\L^{X,(\Eu^{-1},E^\vee,l+1)} \subset \K^{\mu,Q'}$ by sending 
\[
J^1 = \sum_d \prod_{i=1}^n (Q_i')^{d_i} \cdot J^1_d = \sum_d \prod_{i=1}^n Q_i^{d_i} \cdot J^1_d \cdot \mu^{l\langle c_1(E),d \rangle}
\]
to 
\begin{align*}
J^2 & = (\overline{D_2})^{-1}\circ \overline{D_1} \cdot \left(\sum_d Q^d\cdot (-1)^{\langle c_1(E),d\rangle} \cdot J^1_d \cdot \mu^{l\langle c_1(E),d \rangle} \right) \\
& = \sum_d \ \prod_{i=1}^n \left((-1)^{c_1(E)_i}Q_i\right)^{d_i} \prod_{k=1}^{\rank E} \frac{(1-\mu^{-1}f_j^{-1}q^{\langle c_1(f_j),d\rangle})^l (1-\mu f_j)^{l+1}}{(1-\mu^{-1}f_j^{-1})^l (1-\mu f_j q^{-\langle c_1(f_j),d\rangle})^{l+1}} \cdot J^1_d \cdot \mu^{l\langle c_1(E),d \rangle} \\
& = \sum_d \ \prod_{i=1}^n \left( Q_i' \cdot \left(-q^{l}\right)^{c_1(E)_i}\right)^{d_i} \cdot J^1_d \cdot \prod_{j=1}^{\rank E}\frac{1-\mu\cdot f_j}{1-\mu \cdot f_j \cdot q^{-\langle c_1(f_j), d\rangle}},
\end{align*}
which proves Corollary \ref{Cor2}. Note that Formula (\ref{J2muinverse}) is not well-defined in $\K^{\mu,Q'}$ so we stick with Formula (\ref{J2mu}) here. The same argument over $\K^{\mu^{-1},Q''}$ gives Corollary \ref{Cor3}, where we stick with Formula (\ref{J2muinverse}) instead of (\ref{J2mu}).

\par Recall that
\[
F_A = \exp{\left( \sum_{k>0} \frac{a^k f_j(P^{k}q^{kQ\partial_{Q}})}{k} \right)} \quad\text{ and }\quad F_B = \exp{\left( \sum_{k>0} \frac{b^k f_j(P^{-k}q^{-kQ\partial_{Q}})}{k} \right)}
\]
are both Type-1 PFD operators of Lemma \ref{PFD} and thus preserve the image cones as well. For Theorem \ref{1A}, the entire proof above needs not to be changed except that we replace the PFD operator $F_2$, which we used to modify $\Box_1$ above, by $F_2\circ F_A$. Eventually, $J^2$ gains an extra factor of $(1-\mu f_j q^{-\langle c_1(f_j),d\rangle})$, which reduces its formula to the one in Theorem \ref{1A}. Note that in this case $J^2 \equiv \Eu(\mu^{-1}E)\cdot J^1 \ (\mod \ Q')$ and $\Eu(\mu^{-1}E) \equiv 1 \ (\mod \ \mu)$, which means if $J^1$ is $\Lambda[[\mu]]^{Q'}_+$-close to $(1-q)$, so is $J^2$. For Theorem \ref{1B}, we replace $F_2$ by $F_2\circ F_B$, and eventually $J^2$ gains an extra factor of $(1-\mu^{-1} f_j^{-1} q^{\langle c_1(f_j),d\rangle})$, which reduces its formula to the one in Theorem \ref{1B}. If $J^1$ is $\Lambda[[\mu^{-1}]]^{Q''}_+$-close to $(1-q)$, so is $J^2$.

\section{Torus-Equivariant Theory} \label{TorusEquivariantTheory}

\par In this section, we assume further that the smooth projective variety $X$ admits a torus action with isolated fixed points and isolated 1-dim orbits connecting these fixed points which are isomorphic to $\C P^1$'s. Such $X$ is called a GKM variety in literature. In order to avoid notational conflict with the localized multiplicatively closed subset in the base coefficient ring which are always denoted by $T$ (see, for instance, Section \ref{variations}), we denote by $G$ the torus acting on $X$.

\par We prove a torus-equivariant version of the KqSD for such $X$. As an application, using the method of torus fixed point localization and abelian/non-abelian correspondence developed in \cite{Yan:flag}, we generalize in Section \ref{NonAbelQuantumSerre} the KqSD to the ``non-abelian'' case where $X$ is a flag variety, with the primitivity assumption on the vector bundle $E$ removed.

\subsection{Loop space formalism for torus-equivariant theory} \label{TorusEquivariantLoopSpaceFormalism}
\par The $G$-equivariant big $\J$-function may be defined in the same way as in Section \ref{bigJ}, except that 
\begin{itemize}
    \item we replace $\{\phi_\alpha\}_{\alpha\in A}$ and $\{\phi^\alpha\}_{\alpha\in A}$ now by bases of $K_G(X)$, the $G$-equivariant K-theory of $X$, dual with respect to the $G$-equivariant Poincaré pairing 
    \[
    \langle \phi_\alpha, \phi^{\alpha'} \rangle_G = \chi_G (X; \phi_\alpha \otimes \phi^{\alpha'}) \quad \in \Rep (G);
    \]
    \item we define the correlators as the $G$-equivariant $S_n$-invariant holomorphic Euler characteristics $\chi_G^{S_n}$ over the moduli spaces, taking values in $\Rep(G)$ as well.
\end{itemize}

\par We denote by $\J^X_G = \J^X_G(\t)$ the $G$-equivariant big $\J$-function of $X$. It takes values in $K_G(X)$ by definition. It is well-known that when $X$ has isolated fixed points, $K_G(X)$ is a free $\K_G(pt) = \Rep(G)$-module generated by the fixed point classes $\{\phi_\alpha\}_{\alpha\in\Fix_G(X)}$, up to a localization of $\Rep(G)$. For simplicity of notation, we denote such localization still by $\Rep(G)$. In order to include $\J^X_G$ into our picture, we need to extend our rational loop space formalism in the previous sections by replacing $K(X)$ with $K_G(X)$, and the minimal coefficient ring $\Lambda = \C[[Q_1,\cdots,Q_n]]$ with $\Lambda^G = \C[[Q_1,\cdots,Q_n]]\otimes \Rep(G)$. Note that the representation rings are naturally endowed with Adams operations, so is $\Lambda^G$. 

\par However, such extension is not enough. In fact, one crucial feature of the loop spaces for torus-equivariant quantum K-theory is that a wider variety of poles of $q$ may arise from the correlators. Indeed, one may use torus fixed point localization to compute the correlator
\[
\left\langle \frac{\phi_\alpha}{1-qL_0}, \t(L_1),\cdots,\t(L_m) \right\rangle^{S_m}_{0,m+1,d},
\]
and the observation in \cite{Givental:perm2} is that two types of poles of $q$ arise from such computation, depending on what happens on the irreducible component of the domain curve where the non-permuted $0$-th marked point $p_0$ resides. For a given stable map, let $\alpha \in X$ be its evaluation at $p_0$, then $\alpha$ must be a $G$-fixed point on $X$. 
\begin{itemize}
    \item If the stable map has degree zero on the irreducible component containing $p_0$, i.e. restricts to the constant map to $\alpha$, the resulting pole of $q$ is at a root of unity.
    \item If the stable map has non-trivial degree, which means the image of this irreducible component is a 1-dim orbit of $G$ branching out from $\alpha$, the resulting pole of $q$ is at $\lambda^{1/m}$, where $\lambda$ is a tangent $G$-character and $m$ is a positive integer. 
\end{itemize}
Similarly, the Kawasaki-Riemann-Roch contribution from the input $\t(L_j)$ on a certain stratum of the inertia stack of the moduli space may now take the form $\Psi^l(\t(c\cdot \overline{L}^l))$, where $l$ is a positive integer, $\overline{L}$ is the topological part of the universal cotangent bundle with $\overline{L}-1$ nilpotent, and $c$ is either a root of unity or a root of tangent character as above. It indicates that $\J^X_G$ is well-defined for input function $\t = \t(q)$ without poles at roots of unity or roots of tangent characters.

\par Therefore, the loop space formalism may be generalized to the $G$-equivariant settings as follows. Take
\[
R_X = \{ \text{ roots of unity } \} \cup \{ \lambda^{1/m}\ | \ \lambda \text{ is a tangent character at a $G$-fixed point of } X, \ m\in\mathbb{Z}_{>0} \}.
\]
\begin{definition} \label{GequivRationalLoopSpace} We define
\begin{align*}
    \widetilde{\K}_+^G & = \C[[Q]] \otimes_{\C} K_G(X) \otimes_{\Rep(G)} (\widetilde{S}^G_+)^{-1} (\Rep(G)\otimes \C[q]), \\
    \widetilde{\K}_-^G & = \C[[Q]] \otimes_{\C} K_G(X) \otimes_{\Rep(G)} \left\{ \frac{a(q)}{b(q)} \in (\widetilde{S}^G_-)^{-1}(\Rep(G)\otimes \C[q]) \ \middle| \ \deg(a) < \deg(b) \right\}, \\
    \widetilde{\K}^G & = \C[[Q]] \otimes_{\C} K_G(X) \otimes_{\Rep(G)} (\widetilde{S}^G)^{-1} (\Rep(G)\otimes \C[q]),
\end{align*}
where 
\[
\widetilde{S}^G_+ = \{ g(q) | g(c)\neq 0, \forall c\in R_X \}, \quad 
\widetilde{S}^G_- = \left\{ f(q) = c \prod_{s=1}^N (1-q c_s) | c\neq 0\in \Rep(G), c_s \in R_X \right\},
\]
\[
\widetilde{S}^G = \widetilde{S}^G_+ \times \widetilde{S}^G_- = (\Rep(G)\otimes \C[q]) \setminus 0.
\]
\end{definition}
Let $\Rep(G)_0$ the field of fractions of $\Rep(G)$. Then all three spaces appearing above are free modules over the coefficient ring $\widetilde{\Lambda}^G = \C[[Q_1,\cdots,Q_n]] \otimes \Rep(G)_0$. It is not hard to check that $\widetilde{\K}^G = \widetilde{\K}_+^G \oplus \widetilde{\K}_-^G$, and that the two subspaces $\widetilde{\K}^G_{\pm}$ are both Lagrangian under the symplectic form
\[
\Omega^G(\f,\g) = \Res_{q \in R_X}\ \langle \f(q^{-1}), \g(q)\rangle_G \cdot \frac{dq}{q}, \quad \f,\g \in \widetilde{\K}^G.
\]
The $G$-equivariant big $\J$-function $\J^X_G(\t)$ is well-defined for $\t\in\widetilde{\Lambda}_+^G \cdot \widetilde{\K}_+^G$ where $\widetilde{\Lambda}_+^G = (Q)\subset \widetilde{\Lambda}^G$. We denote its image by $\widetilde{\L}^X_G$. It is a $\widetilde{\Lambda}^G_+$-germ at $1-q$. 

\par By our assumption that $G$-fixed points on $X$ are isolated, $K_G(X)$ admits a $\Rep(G)$-basis given by the dual fixed point classes $\{\phi^\alpha\}_{\alpha\in A = \Fix_G(X)}$. Therefore, there exists decomposition
\[
\widetilde{\K}^G = \bigoplus_{\alpha\in\Fix_G(X)} \widetilde{\K}^{pt}, \quad \text{ where } \quad \widetilde{\K}^{pt} = \C[[Q]] \otimes_{\C} (\widetilde{S}^G)^{-1} (\Rep(G)\otimes \C[q]).
\]
In other words, the $G$-equivariant loop space over $X$ is a direct sum of (non-equivariant) rational loop spaces $\widetilde{\K}^{pt}$ (see Section \ref{variations}) over point target spaces (the $G$-fixed points), but with the enlarged coefficient ring $\widetilde{\Lambda}^G$. One should note, however, that $\widetilde{\K}^G_{\pm}$ do not decompose accordingly into $\widetilde{\K}^{pt}_{\pm}$ as $R_X\neq R_{pt}$. For instance, elements like $1/(1-\lambda q^m)$, where $m$ is an integer and $\lambda$ is a tangent character in $T_\alpha X$, appear in $\widetilde{\K}^G_{-}$ but localize to $\widetilde{\K}^{pt}_{+}$. Such feature is crucial in the recursive characterization which we will use later.

\par Under such decomposition, the component of $\J^X_G(\t)$ at the fixed point $\alpha$ is exactly
\[
\J^X_G(\t)|_\alpha = 1 - q + \langle \t(q), \phi_\alpha\rangle_G + \sum_{d,n} Q^d \left\langle \frac{\phi_\alpha}{1-qL_0}, \t(L_1),\cdots,\t(L_m) \right\rangle^{S_m}_{0,m+1,d},
\]
with correlator contribution coming solely from the locus in the moduli space where the $0$-th marked point is sent to $\alpha\in X$. In this way, through fixed point localization, it is not hard to realize that the correlators in $\J^X_G(\t)|_\alpha$ create only poles of $q$ at 
\[
R_\alpha = \{ \text{ roots of unity } \} \cup \{ \lambda^{1/m}\ | \ \lambda \text{ is a tangent character at } \alpha, \ m\in\mathbb{Z}_{>0} \},
\]
instead of the entire $R_X$. Therefore, one may even further extend the domain of definition of $\J^X_G(\t)$: instead of asking $\t = \sum_\alpha \langle \t(q), \phi_\alpha\rangle_G \phi^\alpha$ to never have poles at $R_X$, we need only to require that $\langle \t(q), \phi_\alpha\rangle_G$ have no poles at $R_\alpha$. We do not need such extension in this paper.

\par Even without such extension, not all points on $\widetilde{\L}^X_G$ descend to points on $\L^X$, the image cone of the non-equivariant big $\J$-function, when we take the non-equivariant limit by setting all $G$-equivariant parameters to $1$, and it is mainly due to two reasons. First, not all points in $\widetilde{\L}^X_G$ admit non-equivariant limits. For instance, we do not forbid elements of the form $1/(\lambda_1-\lambda_2)$ in $\widetilde{\K}^G$, where $\lambda_1$ and $\lambda_2$ are any two (one-dimensional) characters of $G$. The denominator vanishes as we take $\lambda_1=\lambda_2=1$. Second, some elements in $\widetilde{\K}^G_+$ may descend not to $\K_+$ but to $\K_-$ under non-equivariant limits. For instance, the element $1/(1-\lambda q^2)$ lives in $\widetilde{\K}^G_+$ when $\lambda\in\Rep(G)$ is not a tangent $G$-character in $X$, but its non-equivariant limit $1/(1-q^2)$ lives in $\K_-$ as it has only poles at roots of unity. 

\par For such reasons, we restrict the above formalism of $G$-equivariant loop space as follows. Let
\[
\widebar{\pi}: \Rep(G) \longrightarrow \C, \quad
\pi = \widebar{\pi}\otimes\id: \Rep(G)\otimes \C[q] \longrightarrow \C[q]
\]
be the maps of taking non-equivariant limits, i.e. sending any 1-dim character of $G$ to $1$. Then,
\begin{itemize}
    \item $T^G := \widebar{\pi}^{-1}(\C \setminus 0)$ is a multiplicatively closed subset of $\Rep(G)$;
    \item $S^G_+ := \pi^{-1}(S_+)$ is a multiplicatively closed subset of $\Rep(G)\otimes \C[q]$, where $S_+ = \{ a(q)/b(q) \ | \ b(\zeta)\neq 0, \forall \text{ root of unity } \zeta\}$ as we have used in Section \ref{RationalInputs}, and it is not hard to see $S^G_+\subset \widetilde{S}^G_+$.
\end{itemize}
\begin{definition}[$G$-equivariant rational loop space] \label{GenuineEquivRationalLoopSpace} 
\begin{align*}
    \K_+^G & = \C[[Q]] \otimes_{\C} K_G(X) \otimes_{\Rep(G)} (S^G_+)^{-1} (\Rep(G)\otimes \C[q]), \\
    \K_-^G & = \C[[Q]] \otimes_{\C} K_G(X) \otimes_{\Rep(G)} \left\{ \frac{a(q)}{b(q)} \in (S^G_-)^{-1}(\Rep(G)\otimes \C[q]) \ \middle| \ \deg(a) < \deg(b) \right\}, \\
    \K^G & = \C[[Q]] \otimes_{\C} K_G(X) \otimes_{\Rep(G)} (S^G)^{-1} (\Rep(G)\otimes \C[q]),
\end{align*}
where $S^G = S^G_+ \times S^G_-$, $S^G_+$ is as above, and 
\[
S^G_- = \left\{ f(q) = c \prod_{s=1}^N (1-q c_s) \ | \ c \in T^G \subset \Rep(G), c_s \in R_X \right\}.
\]
\end{definition}
$\K^G$ as well as the two subspaces $\K_{\pm}^G$ are now free modules over the coefficient ring $\Lambda^G = (T^G)^{-1}\Rep(G)\otimes \C[[Q]]$. We take the adic topology with respect to $\Lambda^G_+ = (Q)$ as usual. Moreover, $\K^G\subset\widetilde{\K}^G$ and $\K^G_{\pm}\subset\widetilde{\K}^G_{\pm}$. The symplectic form $\Omega^G$ on $\widetilde{\K}^G$ is still non-degenerate when restricted to $\K^G$.

\par Another advantage of considering $\K^G$ is that variations appearing in Section \ref{variations} are more readily generalized to $G$-equivariant settings: not much needs to be changed in Definition \ref{GenuineEquivRationalLoopSpace}. Let $\K^{T,S^T_+}$ be the rational loop space in Section \ref{variations} over the coefficient ring $\Lambda = T^{-1}\C[\lambda_1,\cdots,\lambda_r]\otimes\C[[Q]]$ determined by the multiplicatively closed subset $S^T_+ \subset \C[q, \lambda_1,\cdots,\lambda_r]$. We denote still by
\[
\widebar{\pi}: \Rep(G)\otimes \C[\lambda_1,\cdots,\lambda_r]\rightarrow \C[\lambda_1,\cdots,\lambda_r], \quad
\pi: \Rep(G)\otimes \C[\lambda_1,\cdots,\lambda_r,q] \longrightarrow \C[\lambda_1,\cdots,\lambda_r, q]
\]
the maps of taking non-equivariant limits. Then, $\K^{G,T,S^T_+}$ below is the $G$-equivariant counterpart of $\K^{T,S^T_+}$.
\begin{definition}[Variation of the $G$-equivariant rational loop space] \label{VariationDefinitionEquivLoopSpace}
\begin{align*}
    \K^{G,T,S^T_+}_+ & = \C[[Q]] \otimes_{\C} K_G(X) \otimes_{\Rep(G)} (S^{G,T}_+)^{-1} (\Rep(G)\otimes \C[q,\lambda_1,\cdots,\lambda_r]), \\
    \K^{G,T,S^T_+}_- & = \C[[Q]] \otimes_{\C} K_G(X) \otimes_{\Rep(G)} \left\{ \frac{a(q)}{b(q)} \in (S^{G,T}_-)^{-1}(\Rep(G)\otimes \C[q]) \ \middle| \ \deg(a) < \deg(b) \right\}, \\
    \K^{G,T,S^T_+} & = \C[[Q]] \otimes_{\C} K_G(X) \otimes_{\Rep(G)} (S^{G,T})^{-1} (\Rep(G)\otimes \C[q]),
\end{align*}
where $S^{G,T} = S^{G,T}_+ \times S^{G,T}_-$, with
\[
S^{G,T}_+ = \pi^{-1}(S^T_+), \quad
S^{G,T}_- = \left\{ f(q) = c \prod_{s=1}^N (1-q c_s) \ | \ c \in T^G := \widebar{\pi}^{-1}(T), c_s \in R_X \right\}.
\]
\end{definition}
$\K^{G,T,S^T_+}$ as well as the two subspaces $\K^{G,T,S^T_+}_{\pm}$ are free modules over the coefficient ring 
\[
\Lambda^{G,T} = (T^G)^{-1} (\Rep(G)\otimes\C[\lambda_1,\cdots,\lambda_r]) \otimes \C[[Q]].
\]
We denote still by $\J^X_G(\t)$ the $G$-equivariant big $\J$-function, defined for $\t \in \Lambda^{G,T}_+\cdot\K^{G,T,S^T_+}_+$ where $\Lambda^{G,T}_+ = (Q)$, and by $\L^X_G \subset \K^{G,T,S^T_+}$ its image cone.

\par We further denote by $\widetilde{\pi}_X: K_G(X)\rightarrow K(X)$ the map of forgetting the $G$-action, then 
\[
\Pi := \id \otimes \widetilde{\pi}_X \otimes \pi: \ \K^{G,T,S^T_+} \ \rightarrow \ \K^{T,S^T_+}
\]
is well-defined by construction of $S^{G,T}_+$ and is exactly the map of taking non-equivariant limits. Moreover, it is not hard to see $\Pi(\K^{G,T,S^T_+}_{\pm}) \subset \K^{T,S^T_+}_{\pm}$.
\begin{lemma} [Non-equivariant limit] \label{LemmaNonEquivLimit}
\[
\Pi(\L^X_G) \subset \L^X.
\]
\end{lemma}
Indeed, since $\Pi(\K^{G,T,S^T_+}_{\pm}) \subset \K^{T,S^T_+}_{\pm}$, the projection of $\Pi(\J^X_G(\t))$ to $\K^{T,S^T_+}_+$ is exactly $(1-q)+\Pi(\t)$. It suffices then to show that $\J^X(\Pi(\t)) = \Pi(\J^X_G(\t))$, which follows from the commutative diagram (by taking $\mathfrak{X} = [\M_{0,m+1}(X,d)/S_m]$) below as the correlators are ($G$-equivariant) holomorphic Euler characteristics. 
\begin{equation*}
\xymatrix{
K_G(\mathfrak{X}) \ \ar^{\chi^{G}\quad\quad}[r] \ar^{\widetilde{\pi}_{\mathfrak{X}}}[d] & K_G(pt) = \Rep(G) \ar^{\widetilde{\pi}_{pt}}[d] \\
K(\mathfrak{X}) \ \ar^{\chi\quad\quad\quad}[r] & \quad K(pt) = \C \quad\quad
}
\end{equation*}

\par Twisted $G$-equivariant big $\J$-functions may be defined on $\K^{G,T,S^T_+}$ as well, and the same lemma regarding non-equivariant limits holds. With such preparation, we are now ready to consider the quantum Serre duality in torus-equivariant quantum K-theory.

\subsection{Torus-equivariant K-theoretic quantum Serre duality} \label{TorusQuantumSerre}
\par Consider the $G$-equivariant primitive vector bundle
\[
E^\vee = \sum_{j=1}^{\rank E} f_j = \sum_{j=1}^{\rank E} \ f_j(P_1,\cdots,P_n) \quad\in K_G(X)
\] 
where $f_j$ are $G$-equivariant line bundles. We state the main theorem as follows.
\begin{theorem}[$G$-equivariant KqSD] \label{TheoremEquivariantQuantumSerre}
For any point 
\[
J^1 = \sum_d \prod_{i=1}^n Q_i^{d_i} \cdot J^1_d \quad \in \L^{X,(\Eu,E,l)} \subset \K^G(\mu)
\]
with $J^1_d$ independent of $Q_i\ (i=1,2,\cdots,n)$, we have $J^2 \in \L^{X,(\Eu^{-1},E^\vee,l+1)} \subset \K^G(\mu)$ where
\begin{align*}
J^2 & = \sum_d \ \prod_{i=1}^n \left( Q_i \cdot \left(-q^{l}\right)^{c_1(E)_i}\right)^{d_i} \cdot J^1_d \cdot \prod_{j=1}^{\rank E} \frac{1-\mu\cdot f_j}{1-\mu \cdot f_j \cdot q^{-\langle c_1(f_j), d\rangle}}\\
& = \sum_d \ \prod_{i=1}^n \left(Q_i \cdot \left(-q^{(l+1)}\right)^{c_1(E)_i}\right)^{d_i} \cdot J^1_d \cdot \prod_{j=1}^{\rank E}\frac{1-\mu^{-1} \cdot f_j^{-1}}{1-\mu^{-1} \cdot f_j^{-1} \cdot q^{\langle c_1(f_j), d\rangle}}.
\end{align*}
Here $\K^G(\mu)$ is the $G$-equivariant counterpart of $\K(\mu)$ (see Section \ref{variations}) given through Definition \ref{VariationDefinitionEquivLoopSpace}, and $c_1(E)_i$ refers to the $i$-th component of $c_1(E)$ under the basis $c_1(P_1),\cdots c_1(P_n)$ of $H_2(X)$. 
\end{theorem}
As usual, if any summand $f_j$ of $E^\vee$ admits a $(-)$ sign, the corresponding factors in $J^2$ should be inverted.

\par We remark that both $(1-\mu \cdot f_j \cdot q^{-\langle c_1(f_j), d\rangle})$ and $(1-\mu^{-1} \cdot f_j^{-1} \cdot q^{\langle c_1(f_j), d\rangle})$ appearing in the theorem are invertible in $\K^G(\mu)$. For instance, for the former we have
\[
\frac{1}{1-\mu \cdot f_j \cdot q^{-\langle c_1(f_j), d\rangle}} = \sum_{\alpha\in\Fix_G(X)} \frac{1}{1-\mu \cdot f_j|_\alpha \cdot q^{-\langle c_1(f_j), d\rangle}} \cdot \phi^\alpha,
\]
and the denominators on RHS, now in $\Rep(G)\otimes \C[\mu,q]$, are indeed invertible in $\K^G(\mu)$ as their non-equivariant limits take the form $(1 - \mu \cdot q^{-\langle c_1(f_j), d\rangle})$ which do not have any poles at roots of unity. 

\par Our proof of the theorem involves a criterion developed in \cite{Givental:perm2} of determining if a given rational function in $q$ represents a point on $\L^X_G$ (or on any twisted image cones) through fixed point localization. Under the decomposition of $\widetilde{\K}^G$ onto fixed points described above, for any point $\f\in\K^G\subset\widetilde{\K}^G$, 
\[
\f(q) = \sum_{\alpha\in\Fix_G(X)} \f_\alpha(q) \cdot \phi^\alpha, \quad \quad \text{ with } \ \f_\alpha(q) = i^*_\alpha \f(q) \in \widetilde{\K}^{pt}
\]
where $i_\alpha: \alpha \rightarrow X$ the inclusion map from the fixed point $\alpha$ to $X$. The criterion is exactly on these ``linear combination coefficients'' $\f_\alpha(q)$.

\begin{lemma}[Recursive criterion for $\L^X_G$] \label{RecursiveCriterion}
A point $\f(q)\in\K^G$ lives in $\L^X_G$ if and only if $\f_\alpha(q)$ satisfies the following two conditions for any $\alpha\in\Fix_G(X)$:
\begin{itemize}
    \item[(i)] $\f_\alpha(q)$ lives in $\widetilde{\L}^{pt}\subset \widetilde{\K}^{pt}$.
    \item[(ii)] For any $G$-character $\lambda$, the residue of $\f_\alpha(q)$ at any root of $\lambda$ satisfies the recursive formula
    \[
    \Res_{q = \lambda^{1/m}}\f_\alpha(q)\frac{dq}{q} = - \frac{Q^{mD}}{m}\frac{\Eu(T_\alpha M)}{\Eu(T_\psi \M_{0,2}(X;mD))} \f_\beta(\lambda^{1/m})
    \]
    if $\lambda$ is a tangent character at $\alpha$, and is zero otherwise. Here $\beta\in\Fix_G(X)$ is the unique fixed point such that the closure of the 1-dim $G$-orbit $\alpha\beta$ has tangent character $\lambda$ at $\alpha$, $D$ is the homological degree of $\overline{\alpha\beta}$, and $(\psi: \C P^1\rightarrow \overline{\alpha\beta}) \in \Fix_G(\M_{0,2}(X;mD))$ is the $m$-sheet covering ramified at the two marked points sent to $\alpha$ and $\beta$ respectively. 
\end{itemize}
\end{lemma}
Note that the poles at $q = \lambda^{1/m}$ for tangent character $\lambda$ come from the correlators in $\J^X_G$ and thus go to $\K^G_-$ under the polarization of $\K^G$, but they may only come from the input and thus are regarded as in $\widetilde{\K}^{pt}_+$ when we specialize to each $\widetilde{\K}^{pt}$. For twisted big $\J$-functions, Lemma \ref{RecursiveCriterion} is still correct, with Condition (i) unchanged and the recursive coefficient in Condition (ii) altered according to the twisting. For larger $G$-equivariant rational loop spaces like $\K^{G,T,S^T_+}$, Lemma \ref{RecursiveCriterion} is still correct, with $\widetilde{\K}^{pt}$ correspondingly extended.

\par The original version of the recursive criterion is proved in \cite{Givental:perm2} for $\L^X_G$ in the classical $G$-equivariant loop space, where only Laurent polynomials are allowed as inputs. However, it is not hard to see that it is correct for the rational $G$-equivariant loop spaces defined above as well. Indeed, in the proof in \cite{Givental:perm2},
\begin{itemize}
    \item that any point living in $\L^X_G$ satisfies the two conditions follows from fixed point localization on the moduli spaces;
    \item that any point in $\K^G$ satisfying the two conditions lives in $\L^X_G$ follows from an induction on the Novikov variables assuming knowledge of $\widetilde{\L}^{pt}$.
\end{itemize}
Both ingredients still hold over the rational loop spaces, so the proof there carries over to our case \emph{verbatim}.

\par Using Lemma \ref{RecursiveCriterion}, we may reduce the proof of Theorem \ref{TheoremEquivariantQuantumSerre} to verifying the following two facts
\begin{itemize}
    \item[(i)] Given any $J^2$ on $\widetilde{L}^{pt}$, $J^1$ as is given in the theorem is still a point on $\widetilde{L}^{pt}$.
    \item[(ii)] Given any $J^2$ satisfying the recursive formulae of $\L_G^{X,(\Eu^{-1},E^\vee,l+1)}$
    \[
    \Res_{q = \lambda^{1/m}} \f_\alpha(q) \cdot \frac{dq}{q} = -\frac{Q^{mD}}{m} \frac{\Eu(T_\alpha M)}{\Eu(T_\psi \M_{0,2}(X;mD))} \frac{\Eu(E^\vee)|_\alpha}{\Eu(E^\vee_{0,2,mD})|_\psi} \frac{\det^{-(l+1)}(E^\vee_{0,2,mD})|_\psi}{\det^{-(l+1)}(E^\vee)|_\alpha} \cdot \f_\beta(\lambda^{1/m}),
    \]
    $J^1$ as is given in the theorem satisfies the recursive formulae of $\L_G^{X,(\Eu,E,l)}$
    \[
    \Res_{q = \lambda^{1/m}} \f_\alpha(q) \cdot \frac{dq}{q} = -\frac{Q^{mD}}{m} \frac{\Eu(T_\alpha M)}{\Eu(T_\psi \M_{0,2}(X;mD))} \frac{\Eu(E_{0,2,mD})|_\psi}{\Eu(E)|_\alpha} \frac{\det^{-l}(E_{0,2,mD})|_\psi}{\det^{-l}(E)|_\alpha} \cdot \f_\beta(\lambda^{1/m}).
    \]
\end{itemize}
Condition (i) is true because the discrepancy between $J^1$ and $J^2$ comes essentially from the specialization at $a=b=1$ of (the ratio of) the hyper-geometric factors produced by $D_1$ and $D_2$ (see Section \ref{ProofQuantumSerre}), but such hypergeometric factors are known to preserve $\widetilde{L}^{pt}$ according to \cite{Givental:perm4}. Condition (ii), on the other hand, may be checked through direct computation. In fact, for each summand $f_j$ in the decomposition of $E$, $f_j$ restricts topologically to $\mathcal{O}(\langle c_1(f_j),D \rangle)$ on the orbit $\overline{\alpha\beta} = \C P^1$, and its restriction on the two end points satisfy 
\[
f_j|_\alpha = \lambda^{\langle c_1(f_j),D \rangle} \cdot f_j|_\beta.
\]
In this way, both $E_{0,2,mD}|_\psi$ and $E^\vee_{0,2,mD}|_\psi$ may be explicitly written out in terms of $E|_\alpha$ and $\lambda$, and thus all extra terms in the recursive coefficients of $\L_G^{X,(\Eu,E,l)}$ and of $\L_G^{X,(\Eu^{-1},E^\vee,l+1)}$ that come from the twistings. It is then not hard to see that the discrepancy between the recursive coefficients of $\L_G^{X,(\Eu^{-1},E^\vee,l+1)}$ and $\L_G^{X,(\Eu,E,l)}$ amounts exactly to the extra factors in $J^1$ compared to $J^2$ evaluated at the regular point $q = \lambda^{1/m}$. This completes our proof of the Theorem. We omit the computational details here.

\par Setting all $G$-equivariant parameters back to $1$, we obtain the non-$G$-equivariant KqSD in Section \ref{StatementQuantumSerre}. In this way, the above argument provides us with an alternative proof of Theorem \ref{TheoremQuantumSerre} for GKM varieties.

\par We end this section by remarking that $G$-equivariant analogues of Corollary \ref{Cor2} and \ref{Cor3} and Theorem \ref{1A} and \ref{1B} are still correct, and may be proved using Lemma \ref{RecursiveCriterion} as above. We only mention the following two key points, and the rest of the proof needs not to be changed:
\begin{itemize}
    \item The recursion coefficients of $\L_G^{X,(\Eu^{-1},E^\vee,l+1)}$ and $\L_G^{X,(\Eu,E,l)}$ may both be written either in terms of $Q_i'$ or in terms of $Q_i''$.
    \item For any element $J = J(q) = \sum_d Q^d J_d(q)$ in the $G$-equivariant loop space, multiplying each $J_d(q)$ by the factor $(1-\mu \cdot f_j \cdot q^{-\langle c_1(f_j), d\rangle})$ does not change the recursive formula that $J$ satisfies. In fact, it suffices to realize that for any fixed degree $d$,
    \[
    \left. 1-\mu \cdot f_j|_{\alpha} \cdot q^{-\langle c_1(f_j), d+mD\rangle} \middle.\right|_{q = \lambda^{1/m}} = \left. 1-\mu \cdot f_j|_{\beta} \cdot q^{-\langle c_1(f_j), d\rangle} \middle.\right|_{q = \lambda^{1/m}}.
    \]
    LHS and RHS of the equality above are exactly the extra factors that would appear respectively in LHS and RHS of the recursive formula, caused by the change $J_d(q)\mapsto J_d(q) \cdot (1-\mu \cdot f_j \cdot q^{-\langle c_1(f_j), d\rangle})$. Such change turns exactly the $J^2$ in Theorem \ref{TheoremQuantumSerre} to the one in Theorem \ref{1A}. The same is true for the factor $(1-\mu^{-1} \cdot f_j^{-1} \cdot q^{\langle c_1(f_j), d\rangle})$, which turns the $J^2$ in Theorem \ref{TheoremQuantumSerre} to the one in Theorem \ref{1B}.
\end{itemize}

\par We state only the $G$-equivariant versions of Theorem \ref{1A} and \ref{1B}.
\begin{manualtheorem}{3A}[$G$-equivariant KqSD] \label{2A}
For any point 
\[
J^1 = \sum_d \prod_{i=1}^n (Q_i')^{d_i} \cdot J^1_d \quad \in \L_G^{X,(\Eu,E,l)} \subset \K^G[[\mu]]^{Q'}
\]
with $J^1_d$ independent of $Q$ (or equivalently, $Q'$), we have
\[
J^2 = \Eu(\mu^{-1} E) \cdot \sum_d \ \prod_{i=1}^n \left( Q_i' \cdot \left(-q^{l}\right)^{c_1(E)_i}\right)^{d_i} \cdot J^1_d
\quad \in \L_G^{X,(\Eu^{-1},E^\vee,l+1)} \subset \K^G[[\mu]]^{Q'}.
\]
Here $\K^G[[\mu]]^{Q'}$ is the $G$-equivariant counterpart of $\K[[\mu]]^{Q'}$; $\Eu(\mu^{-1} E)\in K_G(X)[[\mu]]$ is the $G\times\C^\times_{\mu}$-equivariant Euler class of $E$ and satisfies $\Eu(\mu^{-1} E) \equiv 1 \ (\mod \ \mu)$.
\end{manualtheorem}
\begin{manualtheorem}{3B}[$G$-equivariant KqSD] \label{2B}
For any point 
\[
J^1 = \sum_d \prod_{i=1}^n (Q_i'')^{d_i} \cdot J^1_d \quad \in \L_G^{X,(\Eu,E,l)} \subset \K^G[[\mu^{-1}]]^{Q''}
\]
with $J^1_d$ independent of $Q$ (or equivalently, $Q''$), we have
\[
J^2 = \Eu(\mu E^\vee) \cdot \sum_d \ \prod_{i=1}^n \left( Q_i'' \cdot \left(-q^{(l+1)}\right)^{c_1(E)_i}\right)^{d_i} \cdot J^1_d
\quad \in \L_G^{X,(\Eu^{-1},E^\vee,l+1)} \subset \K^G[[\mu^{-1}]]^{Q''}.
\]
Here $\K^G[[\mu^{-1}]]^{Q''}$ is the $G$-equivariant counterpart of $\K[[\mu^{-1}]]^{Q''}$; $\Eu(\mu E^\vee)\in K_G(X)[[\mu^{-1}]]$ is the $G\times\C^\times_{\mu}$-equivariant Euler class of $E^\vee$ and satisfies $\Eu(\mu E^\vee) \equiv 1 \ (\mod \ \mu^{-1})$.
\end{manualtheorem}

\subsection{Over non-abelian GIT quotients} \label{NonAbelQuantumSerre}
\par The $G$-equivariant argument above allows us to generalize the (non-equivariant) KqSD to certain cases where the vector bundle $E$ in consideration is not necessarily generated by line bundles in $K(X)$, through the so-called idea of abelian/non-abelian correspondence. 

\par Let $X = V // H = V^s(H)/H$ be the GIT quotient of a vector space $V$ by a reductive group $H$, where $V^s(H)$ refers to semistable locus with respect to a given stability condition. We assume that semistable is equivalent to stable in our case, and that $H$ acts freely on $V^s(H)$. We denote by $S\subset H$ a maximal torus of $H$, by $V^s(S)$ the stable locus of the $S$-action, and assume that $S$ acts freely on it. In this case, we have the following diagram of quotients 
\begin{equation*}
\xymatrix{
	V^s(H)/S\ \ar@{^(->}^{\iota\quad\ }[r]\ar^q[d] & \ Y = V^s(S)/S\\
	X = V^s(H)/H &
}
\end{equation*}
where $\iota$ is an open embedding and $q$ is free quotient by $H/S$. We call $Y = V^s(S)/S$ the abelian quotient associated to $X$. It is a toric variety by definition. Moreover, for any vector bundle $E$ over $X$, we may find vector bundle $F$ over $Y$ such that $\iota^* F = q^* E$. We call such $F$ a lifting of $E$ over $Y$. Note that since $Y$ is a toric variety, the lifting $F$ is always primitive, of which the KqSD is known.

\par The idea of abelian/non-abelian correspondence, first introduced to quantum cohomology through the papers \cite{BCK1} and \cite{BCK2}, is to relate big $\J$-functions of $X$ to those of $Y$. For quantum K-theory, the idea is justified in \cite{Yan:flag} for the case where $X$ is a partial flag variety. More precisely, we treat the partial flag variety $X = \text{Flag}(v_1,\cdots, v_n; N)$ as a GIT quotient by
\[
X = V // H = \Hom(\mathbb{C}^{v_1},\mathbb{C}^{v_2})\oplus\cdots\oplus\Hom(\mathbb{C}^{v_n},\mathbb{C}^N) // GL(v_1)\times\cdots\times GL(v_n).
\]
Both $X$ and its associated abelian quotient $Y$ admit natural action by $G = (\C^\times)^N$ induced from the standard action of $G$ on $\mathbb{C}^N$. Then, it is proven that the image cone $\L^X_G$ of the $G$-equivariant big $\J$-function of $X$ may be recovered as the Weyl-group-invariant part of the image cone $\L^{tw,Y}_{\widetilde{G}}$ of the $(\Eu,\mathfrak{h}/\mathfrak{s})$-twisted $\widetilde{G}$-equivariant big $\J$-function of $Y$ under specialization of Novikov variables (as there are more of them for $Y$ compared to $X$). Here $\widetilde{G}$ is an enlarged torus acting on $Y$ with isolated fixed points and 1-dim orbits, and $\mathfrak{h}/\mathfrak{s}$ denotes the bundle over $Y$ associated to the $S$-representation $\mathfrak{h}/\mathfrak{s}$, where $\mathfrak{h}$ and $\mathfrak{s}$ are the Lie algebrae of $H$ and $S$ respectively. Moreover, twisting on the big $\J$-function of $X$ that is defined by a vector bundle $E$, either of the invertible type like $(\Eu,E)$ and $(\Eu^{-1},E)$, or in the form of level structure like $(E,l)$, may be recovered in the same way from the corresponding twisting on the $(\Eu,\mathfrak{h}/\mathfrak{s})$-twisted big $\J$-function of $Y$ but defined by the lifting vector bundle $F$ of $E$, thus either of the invertible type like $(\Eu,F)$ and $(\Eu^{-1},F)$, or in the form of level structure like $(F,l)$.

\par The $\widetilde{G}$-equivariant KqSD holds for $F$ over $Y$ as the entire $K(Y)$ is primitive. Taking the Weyl-group-invariant part, specializing the Novikov variables, and restricting $\widetilde{G}$ back to $G$, we may then prove the $G$-equivariant KqSD for any vector bundle $E$ over $X$, not necessarily primitive itself. For simplicity of notation, we state the result below only for $V_j\ (j=1,2,\cdots,n)$, the tautological bundles over the flag variety $X$. For each fixed $j$, we denote by $P_{js}\ (s=1,\cdots,v_j)$ the K-theoretic Chern roots of $V_j$. That is to say, exterior powers of $V_j$ are expressed formally as elementary symmetric polynomials in $P_{js}$. We take the Novikov variables $Q_1,\cdots,Q_n$ as associated to the basis of $H^2(X)$ given by $c_1(V_1),\cdots,c_1(V_n)$.
\begin{theorem}[KqSD for flag varieties]
Let $V_j\ (j=1,2,\cdots,n)$ be the $j$-th tautological bundle over the partial flag variety $X = \text{Flag}(v_1,\cdots, v_n; N)$. For any point 
\[
J^1 = \sum_{d = (d_{is}),d_{is}\geq 0} \ \prod_{i=1}^n Q_i^{\sum_{s=1}^{v_i}d_{is}} \cdot J^1_d \quad \in \L_G^{X,(\Eu,V_j,l)} \subset \K^G(\mu)
\]
with $J^1_d$ independent of $Q_i\ (i=1,2,\cdots,n)$, we have $J^2 \in \L_G^{X,(\Eu^{-1},V_j^\vee,l+1)} \subset \K^G(\mu)$ where 
\begin{align*}
J^2 & = \sum_{d} \ \prod_{i=1}^n \left( Q_i \cdot \left(-q^{-l\cdot \delta_{i,j}}\right)\right)^{\sum_{s=1}^{v_i}d_{is}} \cdot \prod_{s=1}^{v_j}\frac{1-\mu\cdot P_{js}^{-1}}{1-\mu \cdot P_{js}^{-1} \cdot q^{-d_{js}}} \cdot J^1_d \\
& = \sum_d \ \prod_{i=1}^n \left( Q_i \cdot \left(-q^{-(l+1)\cdot \delta_{i,j}}\right)\right)^{\sum_{s=1}^{v_i}d_{is}} \cdot \prod_{s=1}^{v_j}\frac{1-\mu^{-1} \cdot P_{js}}{1-\mu^{-1} \cdot P_{js} \cdot q^{d_{js}}} \cdot J^1_d.
\end{align*}
\end{theorem}
\begin{manualtheorem}{4A}[KqSD for flag varieties] \label{3A}
For any point 
\[
J^1 = \sum_{d = (d_{is}),d_{is}\geq 0} \ \prod_{i=1}^n (Q_i')^{\sum_{s=1}^{v_i}d_{is}} \cdot J^1_d \quad \in \L_G^{X,(\Eu,V_j,l)} \subset \K^G[[\mu]]^{Q'}
\]
with $J^1_d$ independent of $Q$ (or equivalently, $Q'$), we have
\[
J^2 = \Eu(\mu^{-1}V_j) \cdot \sum_{d} \ \prod_{i=1}^n \left( Q_i' \cdot \left(-q^{-l\cdot \delta_{i,j}}\right)\right)^{\sum_{s=1}^{v_i}d_{is}} \cdot J^1_d \quad \in \L_G^{X,(\Eu^{-1},V_j^\vee,l+1)} \subset \K^G[[\mu]]^{Q'}.
\]
\end{manualtheorem}
\begin{manualtheorem}{4B}[KqSD for flag varieties] \label{3B}
For any point 
\[
J^1 = \sum_{d = (d_{is}),d_{is}\geq 0} \ \prod_{i=1}^n (Q_i'')^{\sum_{s=1}^{v_i}d_{is}} \cdot J^1_d \quad \in \L_G^{X,(\Eu,E,l)} \subset \K[[\mu^{-1}]]^{Q''}
\]
with $J^1_d$ independent of $Q$ (or equivalently, $Q''$), we have
\[
J^2 = \Eu(\mu V_j^\vee) \cdot \sum_d \ \prod_{i=1}^n \left( Q_i'' \cdot \left(-q^{-(l+1)\cdot \delta_{i,j}}\right)\right)^{\sum_{s=1}^{v_i}d_{is}} \cdot J^1_d
\quad \in \L_G^{X,(\Eu^{-1},E^\vee,l+1)} \subset \K^G[[\mu^{-1}]]^{Q''}.
\]
\end{manualtheorem}
The degree term takes the form $d = (d_{is})_{i=1,s=1}^{n\quad v_i}$ because the summation is originally over curve classes on $Y$ with Novikov variables $\{Q_{is}\}_{i=1,s=1}^{n\quad v_i}$, and is reduced to $X$ only after the specialization of $Q_{is}\mapsto Q_i, \ \forall i,s$.

\par The above formulae obviously admit non-equivariant limits. Setting equivariant parameters of $G$ back to $1$, we obtain the non-torus-equivariant KqSD for the flag variety $X$. The formulae for $J^2$ will not change. This recovers Corollary \ref{CorollaryA} and \ref{CorollaryB} in the Introduction, which hold respectively in $\K[[\mu]]$ and in $\K[[\mu^{-1}]]$.

\bibliographystyle{abbrv}
\bibliography{notes}

$\,$\
\noindent
\textsc{Sorbonne Université and Université Paris Cité, CNRS, IMJ-PRG, F-75005 Paris, France}

\textit{e-mail address:} \href{mailto:xiaohanyan@imj-prg.fr}{xiaohanyan@imj-prg.fr}

\end{document}